\documentclass[fontsize=12pt,a4paper,headings=normal,
twoside=false,leqno,parskip=half-,abstract=true]{scrartcl}
\usepackage[utf8]{inputenc}
\usepackage{CJKutf8} 
\usepackage[greek,english]{babel}
\usepackage{csquotes}
\usepackage[official]{eurosym}
\usepackage[hyphens]{url}
\usepackage{hyperref}

\hypersetup{
 pdftitle={complex planar ODEs},
 pdfauthor={Bernold Fiedler},
 colorlinks=true,
 linkcolor=blue,
 citecolor=blue,
 filecolor=blue,
 urlcolor=blue} 

 \usepackage{afterpage}

\usepackage{graphicx}
\usepackage{wrapfig} 
\usepackage[format=plain,labelfont=bf,font=small]{caption}
\usepackage{xcolor}
\usepackage[arrow, matrix, curve]{xy}
\usepackage{float}

\usepackage{caption}
\usepackage{tabulary}
\usepackage{array}
\newcolumntype{N}[1]{>{\centering\arraybackslash}m{#1}}

\usepackage{amsmath,amsthm}
\swapnumbers 
\usepackage{amssymb, eurosym} 

\makeatletter
\newcommand{\tpitchfork}{%
  \vbox{
    \baselineskip\z@skip
    \lineskip-.52ex
    \lineskiplimit\maxdimen
    \m@th
    \ialign{##\crcr\hidewidth\smash{$-$}\hidewidth\crcr$\pitchfork$\crcr}
  }%
}
\makeatother
\usepackage{latexsym}
\usepackage{enumerate}

\usepackage[notref,notcite,color,final 
]{showkeys}

\definecolor{refkey}{rgb}{1,0,0}
\definecolor{labelkey}{rgb}{1,0,0}

\usepackage{cancel}

\usepackage{tikz}

  \mathchardef\ordinarycolon\mathcode`\:
  \mathcode`\:=\string"8000
  \begingroup \catcode`\:=\active
    \gdef:{\mathrel{\mathop\ordinarycolon}}
  \endgroup

\theoremstyle{plain}
\newtheorem{thm}{Theorem}[section]

\newtheorem{prop}[thm]{Proposition}
\newtheorem{cor}[thm]{Corollary}

\newtheorem{defi}[thm]{Definition}
\newtheorem{rem}[thm]{Remark}

\newcommand\mi{\mathrm{i}}
\newcommand\eps{\varepsilon}
\renewcommand\theta{\vartheta}
\renewcommand\rho{\varrho}
\renewcommand\phi{\varphi}
\renewcommand\Re{\mathrm{Re}\,}
\renewcommand\Im{\mathrm{Im}\,}

\newcommand{\N}{\mathbb{N}}
\newcommand{\Z}{\mathbb{Z}}
\newcommand{\Q}{\mathbb{Q}}
\newcommand{\R}{\mathbb{R}}
\newcommand{\C}{\mathbb{C}}
\renewcommand{\S}{\mathbb{S}}

\def\fx{\mathfrak{x}}
\def\fy{\mathfrak{y}}
\def\fF{\mathfrak{F}}
\def\ff{\mathfrak{f}}
\def\fg{\mathfrak{g}}
\def\fh{\mathfrak{h}}
\def\fH{\mathfrak{H}}
\def\fL{\mathfrak{L}}
\def\fw{\mathfrak{w}}

\def\bxi{\boldsymbol{\xi}}
\def\bz{\boldsymbol{\zeta}}

\def\CP{\C\mathrm{P}}
\def\Co{\mathbb{C}\setminus\{0\}}

\def\cR{\mathcal{R}}

\def\cC{\mathcal{C}}
\def\cD{\mathcal{D}}

\def\cT{\mathcal{T}}

\def\ss2{\mbox{\footnotesize $\sqrt{2}$}}

\begin{document}

\title{\LARGE{Beyond real blow-up:\\ 
\medskip
Masuda detours\\
\medskip
and complex holonomy
}}
\vspace{1cm}
{\subtitle{}
	\vspace{1ex}
	{}}\vspace{1ex}

\author{
 \\
\emph{Dedicated to my teacher Willi Jäger} \\ 
\emph{on the occasion of his 85th birthday}\\
{~}\\
Bernold Fiedler*
\\
\vspace{2cm}}

\date{\small{version of \today}}
\maketitle
\thispagestyle{empty}

\vfill

*\\
Institut für Mathematik\\
Freie Universität Berlin\\
Arnimallee 3\\ 
14195 Berlin, Germany


\newpage
\pagestyle{plain}
\pagenumbering{roman}
\setcounter{page}{1}

\begin{abstract}
\noindent
For real $\mathbf{b}$, consider quadratic heat, or Fujita, equations like
\begin{equation*}
\mathbf{w}_t=\mathbf{w}_{\bxi\bxi} + \mathbf{b}(\bxi)\,\mathbf{w}^2  
\end{equation*}
on $\bxi\in(0,\pi)$ with Neumann boundary conditions.
For $\mathbf{b}$=1, pioneering work by Kyûya Masuda in the 1980s aimed to circumvent PDE blow-up, which occurs in finite real time, by a detour which ventures through complex time.
Naive projection onto the first two Galerkin modes $\mathbf{w}=x+y \cos\bxi$ leads to a crude ODE caricature.
As in the PDE, spatially homogeneous solutions $y=0\neq x\in\R$ starting at $x_0$ blow up at finite real time $t=T=1/x_0$\,.
We aim for ODE ``linearization at infinity''.

\smallskip\noindent
In the spirit of Masuda, we extend real analytic solutions of polynomial planar ODEs to complex time $t$, and to real 4-dimensional $(x,y)\in\C^2$, to circumvent a real blow-up singularity at $t=T$.
More specifically we study the resulting complex foliations, in projective compactifications like $u=1/x,\ z=y/x$, including their holonomy at blow-up $u=0$.
We obtain analytic linearizations, at blow-up equilibria of Poincaré and Siegel type, based on spectral nonresonance.
Passage to imaginary time is an essential source of quasiperiodicity, even at real spectrum.
For iterated Masuda detours in the crude ODE caricature, we discuss the consequences of rational periodic nonresonance, and of irrational quasiperiodic nonresonance of Diophantine type.
Correspondingly, our results distinguish and highlight finite complex branching of solutions, at blow-up, versus almost-closure of irrationally quasiperiodic type.

\smallskip\noindent
Since iterated complex time loops are not feasible, for parabolic PDEs, our PDE-motivated approach is currently limited to ODEs.
On the other hand, all ODE results of the present paper exactly embed into certain PDEs of parabolic type, which possess a PDE-invariant Galerkin subspace.
We conclude with some comments on global aspects, PDEs, discretizations, and other applications.
\end{abstract}

{\small{
\tableofcontents
\listoffigures
}}


\newpage
\pagenumbering{arabic}
\setcounter{page}{1}

\section{Dedication and acknowledgment}\label{Ack} 

\numberwithin{equation}{section}
\numberwithin{figure}{section}
\numberwithin{table}{section}

This paper is dedicated to my dear teacher Willi Jäger on the occasion of his 85th birthday.
Vast is the range of topics and fields which he charismatically initiated in his lively group.
Boundless, almost, were the freedoms enjoyed by those who had the privilege to join.
Sparse and never limiting, but pointed and passionately emphatic, was his paternal advice, prodding us on.
For example: \emph{``Aber machen Sie doch endlich mal was \emph{Quasiperiodisches}, nicht wieder dieses periodische Zeug!''}
[``But do try something \emph{quasiperiodic}, for once, not that periodic stuff again!'']
Emphasis his, not mine.
And time flies.
With a delay of some 45 years, oh well, I at least try to deliver: see the Siegel domains of section \ref{SieLin}, and the irrationally winding holonomies and foliations in most examples.

\emph{``Aber Herr Fiedler -- Sie müssen doch vom \emph{Problem} ausgehen und nicht von der \emph{Methode}!''}
[``\ldots -- Thou shalt not start from the [mathematical] \emph{method}, but from the \emph{problem}!'']
The paternal commandments, repeated on many an appropriate occasion, kept ringing in my ears. 
Standard PDE methods in the problem of real blow-up are tricky comparison functions, Sobolev-like estimates galore, not to forget sophisticated integration by parts, and variational methods.
At least in the simplistically planar ODE variant, I attempt to pursue complex methods of complex dynamics, instead.
Prioritizing the problem over either method might perhaps excuse my patent lack of expertise, in both.

I had much help. Vassili Gelfreich very patiently explained his profound and motivating work on exponential splitting of homoclinic orbits. Anatoly Neishtadt influenced me with lucid hints on adiabatic elimination. 
Tibor Szabó helped along in combinatorics, with instant precision.
Karsten Matthies, Carlos Rocha, Jürgen Scheurle, Hannes Stuke, Nicolas Vassena, and the late colleagues Marek Fila and Claudia Wulff, have provided lasting interest and motivation in real and complex times.
And the tormented sighs of the most knowledgeable \emph{and} patient referees, poor souls, diligently encouraged many interspersed comments, in this long and meandering contribution.

\section{Nontechnical overview}\label{Over}

This paper explores blow-up at finite time $t\!=\!T$ for ODEs (ordinary differential equations).
Our perspective, however, is motivated by pioneering work of Masuda  on blow-up in a certain quadratic parabolic PDE (partial differential equation) \cite{Masuda1, Masuda2}.
Although most real solutions encounter complete blow-up, in finite real time, Masuda showed how to circumvent ``complete'' blow-up via detours in \emph{complex time} $t\in\mathbb{C}\setminus\{T\}$, which recombine at real times $t$ after blow-up.
Mostly, however, the two complex continuations fail to recombine and coincide after detours via the upper and lower complex half planes of $t$, and after blow-up. 
They even fail to turn real again, upon return of $t$ to the real axis.
See section \ref{Masuda} for a few more details.

Encouraged, as well as humbled, by such PDE motivation, we address the Masuda problem of complex-time continuation beyond real-time blow-up in a much simpler ODE setting
\begin{equation}
\label{ODEFt}
\dot{\mathbf{x}}(t)=\mathbf{F}(\mathbf{x}(t))
\end{equation}
of (vector) polynomials $\mathbf{F}$ of degree $m$, and with solutions $\mathbf{x}(t)\in\mathbb{C}^N$.
For illustration we will include a case study of ODE results for our crude 2-mode Galerkin caricatures  of Masuda type PDEs in sections \ref{C2}, \ref{Fuji}, and \ref{Con}.
More generally, and more accurately, section \ref{PDElift} will embed our ODE results into a parabolic PDE context.

Our restriction to ODEs \eqref{ODEFt} eliminates issues like solutions of parabolic PDEs in backwards complex time cones; see remark \ref{remdisclaimer}.
For example, we may consider a local closed time loop 
\begin{equation}
\label{gammat}
\begin{aligned}
   \gamma^t\,{:}\quad\ \S^1&\rightarrow\C\setminus\{T\}\,,   \\
     s\ &\mapsto\  t(s)\,,
\end{aligned}
\end{equation}
of nonzero finite winding number $\fw(\gamma^t)\neq0$ around an isolated blow-up time $T$.
We then solve the underlying ODE \eqref{ODEFt} along such a loop, i.e. along the real parametrization $s$ of complex $t$ with (not necessarily minimal) period $2\pi$.
The annotation of the curve $\gamma^t$ by \emph{superscript} $t$ indicates, that the values $\gamma^t(s):=t(s)$ refer to time $t$ rather than, say, the $t$-dependent ODE solution $\mathbf{x}(t)$ or some components of $\mathbf{x}$.
In fact, the ODE solution $\mathbf{x}(t)$ may, or may not, be periodic upon closure of the complex time loop $s\mapsto t\in\gamma^t$.
In section \ref{LoopStar}, we talk about \emph{blow-up loops} when the solution $s\mapsto x(t(s))$ closes up over the loop $\gamma^t$, i.e. $x(t(2\pi))=x(t(0))$.
In other words, the \emph{discrepancy} 
\begin{equation}
\label{discrepancy}
\mathbf{x}(t(s))\Big|_{s=0}^{s=2\pi}
\end{equation}
of the solution $\mathbf{x}(t(s))$ vanishes, along the complex time loop $\gamma^t$.
See definition \ref{defloop}.
In definition \ref{defstar} we then identify \emph{blow-up stars} associated to minimal blow-up loops of winding number $\fw(\gamma^t)\neq0$.
These are configurations of a number $|\fw(\gamma^t)|$ of radial blow-up branches $\mathbf{x}(t)$ along real times $0>t\nearrow T=0$ \emph{before blow-up}, alternating with $|\fw(\gamma^t)|$ blow-down branches $\mathbf{x}(t)$  along real times $0=T\swarrow t>0$ \emph{after blow-up}.
See figures \ref{figm34}, \ref{figm234} for illustrations of winding numbers $\fw(\gamma^t)=m-1=1,2,3$.

Our main theoretical results on blow-up loops and blow-up stars, based on ``linearization at infinity'', are only formulated in section \ref{B}.
After an introduction of the main concepts, in section \ref{Int}, the example sections \ref{Ric}-\ref{RecC2} gently illustrate these concepts.
Sections \ref{C2m} and \ref{Fuji}, \ref{Con} then apply the results of section \ref{B} to two classes of applications: $m$-homogeneous polynomial systems, and the motivating Galerkin caricature of Masuda's PDE mentioned in the abstract.
In section \ref{Ham} on polynomial Hamiltonian systems, our third class of applications, the linearization approach of section \ref{B} fails.
Instead, we relate the energy levels of blow-up solutions to the very classical language of complex algebraic curves.

For simplicity of presentation, we focus on the complex planar (i.e. real 4-dimensional) polynomial ODEs \eqref{ODEFt} with $\mathbf{x}=(x,y)\in\mathbb{C}^2$.
We rephrase blow-up as a polynomial ODE in terms of suitable \emph{projective coordinates} like $u\!=\!1/x,\ z\!=\!y/x$, on which we comment below; see also \eqref{fuz}.
In particular, we rescale complex time $t$ such that it takes infinite rescaled time $t_1$ to reach blow-up.
More precisely, blow-up to $x\!=\!\infty$ occurs towards a finite equilibrium $(u,z)=(0,e)$, in projective coordinates.
This makes $(u,z)$ a quite intuitive choice to study blow-up.
It is therefore natural to base our analysis of complex blow-up on diagonal linearization at blow-up equilibria $(0,e)$.

Our first result, theorem \ref{thmMs} in section \ref{Ls}, addresses blow-up loops of $\mathbf{x}(t)$, alias $(u(t),z(t))$, in complex time $t$, and their associated blow-up stars in real time $t$, which occur within a complex one-dimensional stable manifold of the blow-up equilibrium.
Therefore we summarize global results on the complex one-dimensional case in section \ref{Scalar}.

When both eigenvalues $\lambda_1,\lambda_2\neq0$ of the linearization at $(u,z)=(0,e)$ become involved, results depend on their \emph{spectral quotient} $\lambda=\lambda_1/\lambda_2\in\Co$.
\emph{Our main result}, theorem \ref{thmloop}, identifies blow-up loops and blow-up stars for positive rational spectral quotients $\lambda=n_1/n_2>0$.
In original time $t$ and projective blow-up coordinates $u,z$ we obtain closed loops $\gamma^t$ and $\gamma^u$ of winding numbers $(m-1)n_1$ and $n_1$\,, respectively; see \eqref{wt}, \eqref{wu}.
Theorem \ref{thmloop} also reformulates this result in terms of the Riemann surface $\cR$ generated by the blow-up solution.
Branching of $\cR$ at the blow-up equilibrium ensues.
At real irrational spectral quotients $\lambda\in\R\setminus\Q$\,, quasiperiodicity sneaks in when we pass to imaginary time.

Linearization in the saddle case of negative spectral quotients $\lambda<0$ usually requires $\lambda$ to be sufficiently irrational; see theorem \ref{thmsadlin}.
Irrationality requirements are expressed in terms of Diophantine conditions \eqref{Dio} or \eqref{Bryuno} which are familiar from KAM theory.
In that case, however, linearization leads to quasiperiodic behavior of the solutions discrepancies
\begin{equation}
\label{discrepancyn}
\mathbf{x}(t(s))\Big|_{s=0}^{s=2n\pi}
\end{equation}
upon n-fold iterations of the closed time loop $s\mapsto\gamma^t(s)$.
On the other hand, the presence of one-dimensional stable and unstable manifolds admits recourse to complex blow-up analysis along the lines of our previous scalar-minded results \cite{FiedlerShilnikov,FiedlerYamaguti}.

To our knowledge, all these results are new, in the context of blow-up.
The techniques, however, are heavily based on the highly recommended book \cite{Ilya} by Ilyashenko and Yakovenko, which still does not seem to have entered the canon of general knowledge, particularly in PDE circles.
Mainly, we employ the concepts of complex foliations, holo\-nomy, and analytic linearization developed there.
We will reference the pertinent parts throughout the paper.
These topics arise as follows.

The introduction of projective coordinates like $(u,z)$ comes at a price.
It often replaces ODEs like \eqref{ODEFt} by scaled ODEs 
\begin{equation}
\label{ODErhoF0}
 \dot{\mathbf{x}} =  1/\rho(\mathbf{x})\cdot\mathbf{F}(\mathbf{x})
\end{equation}
which involve a complex scalar factor $\rho=\rho(\mathbf{x})$.
See \eqref{ODErhoF}, \eqref{ODEuz}, \eqref{ODEvw} below.
The passage from eigenvalues $\lambda_1,\lambda_2\neq0$ to spectral quotients $\lambda=\lambda_1/\lambda_2\in\Co$ in our main results is another example.

Real factors $\rho=\rho(\mathbf{x})\neq 0$ just amount to a rescaling of real time $t$.
For real $\rho$ and t, indeed, the rescaling just stretches the tangent vectors of solutions by a factor $\rho$.
Time orbits, as sets, are preserved.
The nonzero factors $\rho$ are sometimes called (real) \emph{Euler multipliers}.
For real planar vector fields, this rescaling is sometimes subsumed in the language of differential forms, specifically of real \emph{Pfaffian 1-forms} with kernel generated by $\mathbf{F}$; see for example \cite{Amann, Hartman, Lang}.
Real Euler multipliers $\rho=\rho(\mathbf{x})$ thus rescale real time, but affect neither orbits nor the kernel of Pfaffians.

\emph{Complex Euler multipliers} $\rho=\rho(\mathbf{x})\neq 0$ in ODE \eqref{ODErhoF0}, quite analogously, just amount to a complex rescaling of complex time $t$.
When looking for blow-up loops $\gamma^t$ in \emph{complex time} $t$, such rescalings are quite admissible.
When looking for blow-up stars, however, i.e. for blow-up behavior in \emph{real time}, complex Euler multipliers are fraught with peril: they do not preserve the real time axis.
The blow-up stars therefore require additional analysis to recover the fate of real time under its complex rescaling.
See sections \ref{LoopStar}, \ref{PoiMas}.

The appropriate language for complex Euler multipliers $\rho(\mathbf{x})\in\C\setminus\{0\}$ in ODE \eqref{ODErhoF0}, then, are \emph{complex foliations} by complex one-dimensional \emph{leafs}.
The leafs are just the nonstationary solution orbits, in complex time.
These are obtained by a standard Picard-Lindelöf contraction mapping argument, or a holomorphic variant of the implicit function theorem.
From a real point of view, nonvanishing complex holomorphic vector fields $\mathbf{F}(\mathbf{x})$  identify the leafs as real 2-dimensional surfaces tangent to the real 2-dimensional complex span of $\mathbf{F}$ at $\mathbf{x}$. 
The Frobenius integrability condition, which amounts to the Cauchy-Riemann equations for $\mathbf{F}$, provides surfaces of ODE solutions in complex time, locally.
The analogue of real Pfaffian 1-forms are complex 1-forms $\omega$ with kernel $\mathbf{F}$.
We still call $\omega$ (complex) Pfaffian.
We summarize these language basics in section \ref{Basic}, mostly following \cite{Ilya}.
Courage, poor soul: this is not as hard as it may at first sound.
For a gentle introduction, we do include background and many examples.

\emph{Holonomy} addresses the behavior of foliations near zeros, alias equilibria, of $\mathbf{F}$ or its projective variants.
The concept is therefore suited particularly well to study Masuda detours around blow-up.
In section \ref{LinC2} we first address linear vector fields $\mathbf{F}$ of degree $m\!=\!1$, where complex time rescalings are not required because $\rho=1$.
After a direct approach by projective dimension reduction, in section \ref{LinProj}, we introduce linear holonomies in section \ref{LinHol}.
Sections \ref{z=e}-\ref{RatFol} then discuss linear ODEs ``at infinity'', i.e. in projective coordinates $(u,z)$, for general degrees $m$.
In particular we illustrate the roles of rational and irrational spectral quotients $\lambda$, at blow-up $(u,z)=(0,e)$ for blow-up loops and blow-up stars. See remark \ref{remirr} for a summary.

Section \ref{RecC2} on reciprocally linear ODEs illustrates the conceptual power of complex foliations.
Blow-up is marked here by vanishing denominators of the rational, non-polynomial vector field \eqref{ODErhoF0}.
An obvious complex Euler multiplier $\rho$ relates this case to the linear case of section \ref{LinC2}.
Results which were based on complex foliations therefore persist, verbatim.
It essentially remains to discuss the real-time effects and vanishing points of the complex multiplier $\rho$. 

After these conceptual preparations and illustrations, we formulate and prove our main results in section \ref{B}. 
Our proofs are based on classical results on analytic linearization at equilibria, as summarized in \cite{Ilya}.
Analytic linearization requires spectral nonresonance in Poincaré domains and, with additional Diophantine constraints, in Siegel domains.
We have adapted these results to our projective coordinates.

As announced above, we present case studies of three types of applications of our results.
\begin{enumerate}
  \item Section \ref{C2m} discusses blow-up for $m$-homogeneous foliations.
This is a generalization of the purely linear case $m\!=\!1$ of section \ref{LinC2} and involves nontrivial complex Euler multipliers $\rho=(1/x)^{m-1}$.
  \item Section \ref{Ham} addresses ODE blow-up for Hamiltonian ODEs and the second order pendulum with polynomial potential.
Our results are based on the classical relation of energy surfaces to complex algebraic curves.
We resort to algebraic curves, in that section, because analytic ``linearization at infinity'' as in section \ref{B} fails to apply.
The simplest interesting case of cubic potentials relates our results back to classical Weierstrass elliptic curves.
  \item Section \ref{Fuji}, at last, returns to the crude 2-mode Galerkin caricature announced as our PDE motivation.
Section \ref{Con} summarize our conclusions on this example, in less technical language.
The restriction to just a 2-mode ODE approximation of a relatively simple PDE blow-up exposes the huge gap which remains on our journey towards a detailed geometric understanding of concepts like complex blow-up loops and real-time blow-up stars in any full PDE setting.
\end{enumerate}
In section \ref{PDElift}, however, we construct a class of parabolic PDE examples where this gap is partially bridged; see variant \eqref{PDEf} of \eqref{PDEb}.
We adapt an old idea from \cite{FSandstede}, on real flow-embedding, to our complex setting.
Specifically, we embed any complex planar ODE vector field \eqref{ODEFt} into a complex scalar parabolic PDE \eqref{PDEf} on the circle, such that the ODE solutions $\mathbf{x}(t)$ define error-free PDE solutions on a complex two-dimensional, PDE-invariant Fourier subspace.
Since polynomial ODEs lift to polynomial PDEs, this lifts all ODE blow-up results of our present paper to that class of parabolic PDEs.

We conclude, in section \ref{Com}, with a discussion of some of the relevant literature, from the complex analysis side. 
In terms of applications, we mention seemingly unrelated real ODE systems which become close cousins when taking a complex variables approach.
Likewise, certain systems of real PDEs are just scalar complex PDEs in disguise.
We particularly recommend the intriguing aspects which link PDEs of different type when passing from real to imaginary ``time''.
For illustration in a parabolic PDE setting, we include some previous results, which point at a rather curious link between the real bounded solutions of real global attractors, in real time, and blow-up in imaginary time, i.e.\ for certain variants of Schrödinger type PDEs.

\section{Introduction: background and setup} \label{Int}

Let us substantiate the preliminary overview and the prerequisites on a technically more precise level.
In section \ref{RBlow-up} we address ODE blow-up in a standard real setting, in particular for real time. 
The scalar Riccati equation $\dot x=x^2$ will serve as our first and most elementary example for complex extensions, including complex time.
Section \ref{Masuda} pays homage to Masuda's original quadratic parabolic PDE setting, and summarizes his results from the 1980s on detours around real-time blow-up via complex time \cite{Masuda1,Masuda2}.
For our much more detailed ODE analysis, we currently have to pass to a crude 2-mode Galerkin caricature in the complex plane $\mathbf{x}=(x,y)\in\mathbb{C}^2$; see section \ref{C2}.
Therefore we provide some basic background on complex ODEs, complex foliations, and complex Euler multipliers, in section \ref{Basic}.
Guided by the intuition of section \ref{Over}, we also compactify $\mathbf{x}=(x,y)\in\mathbb{C}^2$ by the projective coordinates $[\fx]=[\xi\!:\!\eta\!:\!\zeta]\in\CP^2$; see \eqref{xyz}.
Projective coordinates play a crucial role for our main theoretical results in section \ref{B}: the description of ODE blow-up in complex time via blow-up loops and blow-up stars.
Section \ref{Fol} collects the coordinate-invariant complex foliations and the specific complex Euler multipliers associated with coordinate changes of ODE \eqref{ODEFt} in projective coordinates.
In section \ref{LoopStar} we arrive at the precise definitions of blow-up loops and blow-up stars, in terms of our projective coordinates.
Based on the background material developed here, the outline \ref{Out} again lists the remaining sections, albeit from a more technical perspective than in the nontechnical overview of section \ref{Over}.

\subsection{Blow-up in real time}\label{RBlow-up}

Solutions $\mathbf{w}=\mathbf{w}(t,\bxi)$ of semilinear, quasilinear, or nonlinear partial differential equations (PDEs) may become unbounded in finite real time $t\nearrow T<+\infty$.
This phenomenon is called \emph{blow-up}.
Blow-up comes in several technical flavors.
For example, we may explore blow-up of the spatial profiles $\bxi\mapsto \mathbf{w}(t,\bxi)$ in norm $|\mathbf{w}(t,\cdot)|$, in various Banach spaces $X$.
Or we may ask for pointwise and single-point blow-up of $|\mathbf{w}(t,\bxi_0)|$, for $t\nearrow T$ and at fixed $\bxi_0$\,, or for self-similar blow-up in unbounded domains, for the relations among all these concepts, and so on.
The term \emph{blow-down} is sometimes used to describe the same phenomenon in reverse, finite backward time $t\searrow T>-\infty$.
For many shades of blow-up in parabolic real PDEs see the most diligent monograph \cite{Quittner}, and the near-exhaustive list of references there.
The main methods of such traditional and very established blow-up analysis are sophisticated applications of variational techniques, related Sobolev type estimates, and comparison principles.

Eschewing spatial aspects, despite their nobly intuitive and geometric appeal, blow-up and blow-down also occurs in finite-dimensional ODE settings 
\begin{equation}
\label{ODEF}
\dot{\mathbf{x}}=\mathbf{F}(\mathbf{x})\,,
\end{equation}
traditionally with $\mathbf{x}\in \mathbf{X}=\R^N$.
For vector fields $\mathbf{F}$ which are locally Lipschitz, the maximal interval of existence $t\in(T_-,T_+)$ of any solution $\mathbf{x}(t)$ is always characterized by blow-down or blow-up $|\mathbf{x}(t)|\rightarrow\infty$, in case $T_-$ or $T_+$ are finite, respectively.
In case the vector field $\mathbf{F}$ is analytic, entire, meromorphic, rational, or just polynomial, solutions naturally extend to \emph{complex arguments} $\mathbf{x}\in \mathbf{X}=\C^N\cong\R^{2N}$, and to complex time.
We extend blow-up terminology to include solutions of ODEs \eqref{ODEF} which reach singularities $|\mathbf{F}|=\infty$ in finite complex time $t\rightarrow T\in\C$.

For ODEs, it makes sense to think of isolated singularities of solutions $t\mapsto\mathbf{x}(t)$ at $t\!=\!T$.
A simple paradigm is the homogeneously quadratic scalar \emph{Riccati equation}
\begin{equation}
\label{Ric0}
\dot x = x^2, \quad\textrm{with} \quad x(t)=1/(-t+1/x_0)\,,
\end{equation}
for nonzero initial conditions $x(0)=x_0\neq 0$.
Note analyticity in time $t\in\C$, except for a simple pole at $t\!=\!T\!=\!1/x_0$\,.
For real $x_0$\,, and real-time blow-up at $t\!=\!T$, the blow-up may be circumnavigated by arbitrarily small semi-circular detours of $t\in\C$, which venture into the upper or lower complex half-plane.
Either complex detour defines one and the same real-time continuation of the real blow-up solution $x(t)\in\C$, which becomes real and finite again, \emph{after} real-time blow-up.
In real time, all complex solutions become homoclinic to the equilibrium $x\!=\!0$, i.e.\ $x(t)\rightarrow 0$ for real $t\rightarrow\pm\infty$.
On the Riemann sphere $x\in\widehat{\C}=\C\cup\{\infty\}$, this includes the real blow-up-down solution through $x\!=\!\infty$.
See section \ref{Ric} for further details.

More than forty years ago, Kyûya Masuda was the first to leave the trodden path and explore such \emph{Masuda detours}, around blow-up at $t\!=\!T$ via complex time -- in a PDE setting.
See his pioneering work \cite{Masuda1,Masuda2}.

\subsection{Circumventing blow-up by Masuda detours in complex time}\label{Masuda}

For specific PDE motivation, consider the quadratic heat (or Fujita) equation
\begin{equation}
\label{PDEw}
\mathbf{w}_t=\mathbf{w}_{\bxi\bxi}+\mathbf{w}^2\,,
\end{equation}
on the real interval $0<\bxi<\pi$, under Neumann boundary conditions $u_{\bxi}=0$ at $\bxi =0,\pi$.
Integration over $\bxi$ implies that the spatially homogeneous solution $\mathbf{w}=0$ is the only real equilibrium $\mathbf{w}_t=0$.
Blow-up at finite real time $t=T>0$ occurs, for example, for positive initial conditions $\mathbf{w}(0,\bxi)=\mathbf{w}_0(\bxi)>0$.
The blow-up is \emph{complete}: any attempt to continue the blow-up solution by monotone approximation from below results in $\mathbf{w}(t,\bxi)\equiv\infty$, for all $t>T$ and all $\bxi$\,. 
See section II.27 in \cite{Quittner}.
Notably, this notion is based on comparison principles and is therefore limited to real solutions.
For sign-changing real variants, see for example \cite{FiedlerMatano}.

Spatially homogeneous solutions $\mathbf{w}(t,\bxi)= x(t)$ of PDE \eqref{PDEw} satisfy the Riccati ODE \eqref{Ric0}.
For almost homogeneous real initial conditions $ \mathbf{w}_0$\,, Masuda was then able to circumnavigate blow-up, at real time $t=T(\mathbf{w}_0)$, via a sectorial detour venturing into the upper or lower half-plane of complex $t$.
Notably, Masuda cautioned, the two detours via positive and negative imaginary parts of $t$ agree in their real-time overlap after blow-up, if \emph{and only if} $\mathbf{w}_0(\bxi)= x_0$ is spatially homogeneous.

Standard ODE techniques to circumvent isolated blow-up at $t\!=\!T$ then involve the analysis of complex \emph{holonomy} maps, i.e.\ of flows in complex time $t$ along simple or multiple loops $\gamma^t\subset\C\setminus\{T\}$ around the blow-up singularity at $t\!=\!T$; see \eqref{gammat} and  \cite{Ilya}.
For PDE \eqref{PDEw}, however, variational and comparison principles as in \cite{Quittner} both fail, as soon as we venture into the complex time domain.
Masuda's approach was therefore based on the fact that the local solution semigroup of PDE \eqref{PDEw} is sectorial, and hence analytic, in forward complex sectors $\arg t<\pi/2$.
See also \cite{Stukediss, Stukearxiv}, with an emphasis on heteroclinic orbits.

\begin{rem}\label{remdisclaimer} \emph{\textbf{Disclaimer.}}\quad
A fundamental PDE obstacle to a deeper complex analysis of the Masuda detours, therefore, is that solutions may not extend into backward complex sectors of $t$.
In particular, we cannot combine upper and lower half-circle detours to cover a single or iterated loops $\gamma^t$ around the singularity at $t\!=\!T$ and, thereby, define an infinite-dimensional version of holonomy.
As a first attempt, our paper will therefore explore the holonomy approach of \cite{Ilya} in the simplest interesting ODE caricature of the Masuda PDE paradigm, only: a Galerkin projection onto the first two Fourier modes.
We return to this topic in section \ref{LoopStar}, from the viewpoint of complex foliations.
\end{rem}

\subsection{A crude Galerkin caricature in $\C^2$}\label{C2}

To gain flexibility, our ODE caricature of complex-time Masuda detours will refer to the slightly more general quadratic PDE
\begin{equation}
\label{PDEb}
\mathbf{w}_t=\mathbf{w}_{\bxi\bxi}+\mathbf{b}(\bxi)\mathbf{w}^2\
\end{equation}
for $0<\bxi<\pi$, again under Neumann boundary condition.
Specifically, consider spatial dependence $\mathbf{b}(\bxi)$ on finite real Fourier sums:
\begin{equation}
\label{bj}
\mathbf{b}(\bxi)=b_0+\sum_k b_k\cos(k\bxi)\,.
\end{equation}
Without loss of generality, nonzero spatial averages $b_0$ of $\mathbf{b}$ allow us to rescale $\mathbf{w}$ such that $b_0=1$.
For \emph{symmetric} $\mathbf{b}(\pi-\bxi)=\mathbf{b}(\bxi)$, Fourier coefficients $b_k$ vanish at odd $k$.

For utmost simplification of $\mathbf{w}$, we only consider the first two eigenmodes $1,\ \cos\bxi$ of $\partial_{\bxi}^2$\,, orthogonal under the $L^2$ scalar product $\langle\cdot,\cdot\rangle$ of the spatial average $\tfrac{1}{\pi}\int_0^\pi$.
Orthogonal Galerkin projection of \eqref{PDEb} onto these eigenmodes, with coefficients $\mathbf{w}=x+y\cos\bxi$\,, yields the crude ODE caricature
\begin{equation}
\label{ODEb}
\begin{aligned}
 \dot x &=  \qquad \quad\ x^2 + b_1xy + \tfrac{1}{4}a y^2\,,  \\
 \dot y &=-y + b_1x^2 + axy  + \tfrac{1}{4}(3b_1+b_3)y^2\,,
\end{aligned}
\end{equation}
with parameter $a:=2+b_2$. For symmetric $\mathbf{b}$, where $b_1=b_3=0$, this simplifies to
\begin{equation}
\label{ODEa}
\begin{aligned}
 \dot x &=  \  x^2 + \tfrac{1}{4}ay^2\,,\qquad\qquad\qquad\qquad\quad\ \  \\
 \dot y &=y(-1 +  ax)\,.
\end{aligned}
\end{equation}
Nonzero ``spatially homogeneous'' solutions $y=0\neq x$ of the symmetric caricature \eqref{ODEa} satisfy the Riccati ODE \eqref{Ric0} and blow up in finite time $t=T=1/x_0$\,.
Since we are interested in Masuda detours around singularities, in complex time, we will consider \eqref{ODEa} for $(x,y)\in\C^2$, and study the compactified flow on two-dimensional complex projective space $\C P^2$.
Alternatively, we will also discuss antisymmetric modifications of $b_2=-2$, i.e.\ $a\!=\!0$, of \eqref{ODEb}. 
Our approach will follow the analysis of complex foliations in \cite{Ilya}, in parts.
See sections \ref{Fuji}, \ref{Con}.
For comments on higher-dimensional variants see section \ref{Dim}.

\subsection{Complex basics}\label{Basic}

For the convenience of the reader, we collect some basic background terminology and notation.
For later use, and to remain specific, we also introduce and discuss projective variables in $\CP^2$. 
Projective compactification of $\C^2\subset\CP^2$ will be prerequisite for the formulation of our main results on blow-up loops $\gamma^t\subset\C\setminus\{T\}$, in complex time, around blow-up time $T$.
See section \ref{LoopStar} and definition \ref{defloop}.

Let us start from ODE \eqref{ODEFt}, \eqref{ODEF} for some given locally analytic vector field $\mathbf{F}(\mathbf{x})$ in an open complex domain $\mathbf{x}\in\cD\subseteq \mathbf{X}:=\C^N$.
Here and below, the terms \emph{``analytic''} and \emph{``analyticity''} refer to local expansions by convergent power series.
\emph{``Holomorphy''} refers to complex differentiability and, therefore, Cauchy-Riemann equations.
For real differentiable functions the two notions coincide, by the Cauchy-Goursat theorem.
\emph{``Entire''} functions are globally analytic, e.g. for all complex time arguments $t\in\mathbb{C}$.
For general background on complex analysis, continuation of analytic functions, and Riemann surfaces, we refer to the textbooks \cite{Forster, Jost, Lamotke} and the standard references there.
More specifically, see \cite{Ilya} for details involving ODE solutions and flows in complex time.

For given initial condition $\mathbf{x}(0)=\mathbf{x}_0$ with local ODE solution $\mathbf{x}(t)$, let $\Phi^t(\mathbf{x}_0):=\mathbf{x}(t)$ define the local solution flow of the ODE.
Then the flow property 
\begin{equation}
\label{flow}
\Phi^{t_2}\circ\Phi^{t_1}= \Phi^{t_1+t_2}, \qquad \Phi^0=\mathrm{Id},
\end{equation}
holds, for any argument $\mathbf{x}_0$ and, locally, for all $t_1,t_2\in\mathbb{C}$ such that the  closed parallelogram spanned  by $t_1,t_2$ in $\mathbb{C}$ is contained in the domain of existence of the local semiflow $\Phi^t(\mathbf{x}_0)$.
This follows from the ODE and Cauchy's theorem.
Note how $\Phi^{t_1}$ and $\Phi^{t_2}$ commute.
For example, the local flow $\Phi^t$ in real time $t=t_1=r$ commutes with the local flow in imaginary  time $t=t_2=\mi s$.
For an application to parabolic versus Schrödinger type PDEs, see section \ref{PDE} below.

The local flow map $(t,\mathbf{x}_0)\mapsto\Phi^t(\mathbf{x}_0)$ is holomorphic in all variables.
For fixed times of small $|t|$, the flow maps $\Phi^t$ define local biholomorphism with holomorphic local inverse $\Phi^{-t}$; see \eqref{flow}.
Locally, in the scalar case $N\!=\!1$ and for fixed nonstationary $\mathbf{x}_0\in\cD,\ \mathbf{F}(\mathbf{x}_0)\neq0$, the ODE \eqref{ODEF} implies that solutions $t\mapsto \Phi^t(\mathbf{x}_0)$ are biholomorphically conformal.

\begin{defi}\label{defequi}
We use the following \emph{standard notions of equivalence} for ODE flows, based on conjugating bijections $\Psi$.
We speak of $C^0$ \emph{equivalence}, if $\Psi$ is a local or global homeomorphism.
\emph{Analytic equivalence} requires biholomorphic $\Psi$.
\emph{Flow equivalence}, of either type, requires flows $\Phi^t$ to be conjugated by $\Psi$, locally or globally, and in real or complex time.
\emph{Analytic flow, or orbit, linearization} asserts analytic flow, or orbit, equivalence, near an equilibrium, to the linearized ODE.
When the linearized ODE is in diagonal, semisimple form we also speak of \emph{analytic flow, or orbit, diagonalization}.
Local solutions, which are maximally extended in any open subset of interest, are called \emph{orbits}.
\emph{Orbit equivalence} only requires $\Psi$ to map orbits to orbits, as sets.
\end{defi}

For an example of global analytic equivalence on the Riemann sphere see the Riccati flows of section \ref{Ric}.
Global $C^0$ orbit equivalence in real time figures prominently in section \ref{Scalar}.
For a simple example see also the linear case of section \ref{LinProj}.
The local linearization and diagonalization results of section \ref{B} review analytic local flow equivalence to linear systems, in complex time and near equilibria.

Extending local solutions, maximally, we may encounter obstacles at finite complex times $T$.
Indeed, obstacles occur when $\lim |\mathbf{x}(t)|=\infty$ or $\liminf \mathrm{dist}(\mathbf{x}(t),\partial\cD)=0$, for $t\rightarrow T$.
We use the term \emph{blow-up} in either case, indiscriminately, when and where the boundary $\partial\cD$ denotes singularities of the vector field $\mathbf{F}$.

Following the Masuda paradigm, we are of course particularly interested in blow-up solutions $\lim |\mathbf{x}(t)|=\infty$ which become unbounded in finite time $t\rightarrow T\in\C$.
We will compactify $\mathbf{x}\in\C^N$ by \emph{homogeneous coordinates} in the complex $N$-manifold of complex projective space $\CP^N$.
As usual, the homogeneous coordinates $[\fx]=[\xi_1\!:\!\ldots\!:\!\xi_N\!:\!\zeta]\in\CP^N$ denote equivalence classes of $0\neq\fx=(\xi_1,\ldots,\xi_N,\zeta)\in\C^{N+1}$ under multiplication by a complex scalar $0\neq\sigma\in\C$.
For example, we represent $\mathbf{x}\in\C^N$ by the class $[\fx]$ of $\fx:=(\mathbf{x},1)$.
For $N\!=\!2$, as in \eqref{ODEb}, \eqref{ODEa}, homogeneous coordinates $[\xi\!:\!\eta\!:\!\zeta]$ compactify the complex plane $\mathbf{x}=(x,y)\in\C^2\subset\CP^2$ by the embedding $[\xi\!:\!\eta\!:\!\zeta]:=[x\!:\!y\!:\!1]\in\CP^2$.
Indeed we attach the Riemann sphere $[\xi\!:\!\eta\!:\!0]\in\CP^1\cong\widehat{\C}\cong \S^2$ at infinity $\zeta\!=\!0$.
More precisely, the compact complex analytic 2-manifold $\CP^2$ is given by the atlas of \emph{projective coordinates} $(x,y),(z,u),(v,w)\in\C^2$. 
For nonzero entries, the coordinates overlap as $[x\!:\!y\!:\!1]=[1\!:\!z\!:\!u]=[w\!:\!1\!:\!v]$.
Nonzero complementing coordinates $z,w$, for example, overlap as $z\!=\!1/w$.
In complete detail:
\begin{equation}
\label{xyz}
\begin{aligned}
    &x=1/u, &y=1/v,\quad &z=1/w,   \\
    &x=w/v, &y=z/u,\quad &z=y/x,   \\
    &w=x/y, &v=u/z,\quad &u=v/w.
\end{aligned}
\end{equation}
This represents infinity $[\xi\!:\!\eta\!:\!0]$ by $u\!=\!1/x\!=\!0$ or by $v\!=\!1/y\!=\!0$.

Another complication for maximal extensions of $\mathbf{x}(t)$ in complex time $t$ arises along closed time-loops $\gamma^t\subset\C\setminus\{T\}$, when blow-up occurs for some $t\rightarrow T$ in the interior of $\gamma^t$.
See \eqref{gammat} for this notation.
Since the Cauchy theorem for $t\mapsto\mathbf{x}(t)$ may fail, in the interior, time-loops $\gamma^t$ may generate multi-valued solutions $\mathbf{x}(t)$.
For a more detailed discussion of the resulting Riemann surfaces $\mathbf{x}(t)\in\cR\subset\C^N$ of maximal solutions in the scalar case $N\!=\!1$, see for example \cite{FiedlerShilnikov}.
For $N\!=\!2$ see also section \ref{RatFol} and theorem \ref{thmloop}.

In the overview section \ref{Over}, we have  observed how complex Euler multipliers $\rho\neq0$ preserve nonstationary orbits under complex time, alias leafs of complex foliations.
We have also announced the language of Pfaffian differential 1-forms as most appropriate in the context of $N\!=\!2$, i.e. for ODEs in $\C^2$. 
For a summary of our setup see also remark \ref{remEuler} below.
We now begin to deliver.

Consider nonstationary $\mathbf{x}_0\in\C^N$\,, i.e.\ $\mathbf{F}(\mathbf{x}_0)\neq 0$.
Then the complex \emph{flow-box theorem} for ODE \eqref{ODEFt}, \eqref{ODEF}  asserts local analytic equivalence of the flow $\Phi^t$ to the translation flow of any nonvanishing constant vector field $\mathbf{F}\equiv\mathrm{const}$.
In particular, we obtain a local \emph{foliation} of the domain $\cD\setminus\mathbf{F}^{-1}(0)\subseteq\C^N$ by complex one-dimensional (and hence real two-dimensional) \emph{leafs} $\fL$\,: the nonstationary orbits of ODE \eqref{ODEFt} in the flow-box, for complex $t$. 
Maximal complex extension of a solution $\mathbf{x}(t)$ is equivalent to the maximal extension of the leaf $\fL\,$ of $\mathbf{x}_0$\,, along with its accompanying local foliation.
Globally, the leafs are immersed in $\cD\setminus\mathbf{F}^{-1}(0)\subseteq\C^N$ but, in general, not embedded.
Our assumption $\mathbf{F}(\mathbf{x}_0)\neq 0$ excludes equilibria $\mathbf{x}_0$ from leafs.
Note how leafs depend on the domain $\cD$ of definition of $\mathbf{F}$.
\emph{Local leafs}, for example, refer to (small) neighborhoods $\cD$ of isolated equilibria $\mathbf{F}(\mathbf{x}_0)= 0$.

The language of holomorphic differential $(N-1)$-forms $\omega$ is particularly well-suited to describe such foliations: structures which consist of  decompositions into ``parallel'' complex lines, alias $\ker\omega$, up to local biholomorphisms. 
The complex two-dimensional case $N\!=\!2$ with local foliations into ODE orbits, in complex time, is our main concern in the present paper.
A suitable language is then provided by differential 1-forms, which we called Pfaffian forms in section \ref{Over}.
In complex coordinates $\mathbf{x}=(x,y)\in\C^2$ consider the complex differential 1-form
\begin{equation}
\label{omega}
\omega:= -g\,dx + f\,dy,
\end{equation}
for holomorphic nonzero $(f,g)=(f,g)(x,y)\in\C^2\setminus\{0\}$.
The holomorphic differential 1-form $\omega$ is a linear form on planar (tangent) vector fields $\mathbf{F}$.
\emph{Pull-back} of $\omega$ under coordinate transformations of $(x,y)$ is defined by substitution inside $f,g$ and formal differentiation of the differential symbols $dx, dy$.

For example, consider the vector field $\mathbf{F}=(f,g)=f\partial_x +g\partial_y$ which describes the ODE
\begin{equation}
\label{ODEfg}
\begin{aligned}
 \dot x &= f(x,y)\,, \\
 \dot y &= g(x,y)\,.
\end{aligned}
\end{equation}
Then \eqref{omega} implies that the foliation of $\cD\setminus\mathbf{F}^{-1}(0)$  by $\mathbf{F}$ is given by $\ker\omega$, i.e.
\begin{equation}
\label{omegaF=0}
\omega(\mathbf{F})=0\,.
\end{equation}
Indeed $\omega(\mathbf{F})= -gf+fg=0$.
In general, the kernel of the duality pairing $\omega(\mathbf{F})$ is invariant with respect to simultaneous pull-back of the form $\omega$ \emph{and} the vector field $\mathbf{F}$ under the same coordinate transformation; see \cite{Amann,Lang}.
But $\omega(1/\rho\cdot\mathbf{F})=0$ holds, likewise, for any scalar holomorphic function $\rho=\rho(\mathbf{x})\in\Co$, which we called a complex Euler multiplier in section \ref{Over}.
Indeed, neither the foliation nor the kernel are affected by $\rho$.

Exact 1-forms $\omega$ correspond to a complex version of Hamiltonian vector fields, 
\begin{equation}
\label{dH}
\omega=dH=H_x\,dx+H_y\,dy\,,
\end{equation}
with respect to the standard symplectic 2-form $\Omega=dx\wedge dy$, as follows.
The standard definition of Hamiltonian vector fields $\mathbf{F}=\mathbf{F}^H$ associated to a Hamilton function $H$ and symplectic 2-forms $\Omega$ is
\begin{equation}
\label{FH}
\omega=\Omega(\mathbf{F}^H,\cdot)=dH\,;
\end{equation}
see \cite{AM, Lang}. 
Comparing with \eqref{omega}, \eqref{ODEfg} we obtain  
\begin{equation}
\label{ODEH}
\begin{aligned}
   \dot x&=f=H_y\,,\\
   \dot y&=g=-H_x\,.
\end{aligned}
\end{equation}
Note how $0=\omega(\mathbf{F}^H)=dH\,\mathbf{F}^H$ implies that the complex Hamiltonian $H$ is constant along leafs of the foliation associated to \eqref{ODEfg}.
In the polynomial case, this leads to algebraic leafs and the associated theory of compact Riemann surfaces.
See section \ref{Hamhom} for further discussion.

Let us return to scalar multiples $1/\rho\cdot\mathbf{F}$ of given ODEs \eqref{ODEfg}; compare \eqref{ODErhoF0}.
Local solutions at nonzero $\mathbf{F}$ describe local complex foliations, by the flow box theorem.
The complex vector field $\mathbf{F}$ describes the field of tangents.
We may interpret multiplication by an \emph{Euler multiplier} $\rho=\rho(\mathbf{x})\neq0$ as a rescaling of complex time, as follows. 
Consider ODE \eqref{ODEF} for $1/\rho\cdot \mathbf{F}$, extended by the trivial ODE $\dot t\!=\!1$:
\begin{equation}
\label{ODErhoF}
\begin{aligned}
 \dot{\mathbf{x}} &=  1/\rho(\mathbf{x})\cdot\mathbf{F}(\mathbf{x})\,, \\
 \dot t &=1\,.
\end{aligned}
\end{equation}
We have already noted invariance of the leafs, the foliation, and the Pfaffian kernel of the first equation, under such Euler multipliers.
After multiplication of both sides by $\rho$, we introduce a new complex time variable $\tau$ such that $dt=\rho d\tau$. Denoting derivatives with respect to time $\tau$ as $'=d/d\tau$, we obtain
\begin{equation}
\label{ODEFrho}
\begin{aligned}
 \mathbf{x}' &= \mathbf{F}(\mathbf{x})\,, \\
 t' &=\rho(\mathbf{x})\,.
\end{aligned}
\end{equation}
Playing ODEs \eqref{ODErhoF} and \eqref{ODEFrho} against each other recommends Euler multipliers $\rho$ \emph{with zeros} as a powerful tool of blow-up analysis.
Zeros of $\rho$, for example, may cancel zeros of $\mathbf{F}$, and thus extend the domain of definition of the foliation \eqref{ODEFrho} beyond $\cD\setminus\mathbf{F}^{-1}(0)$.
Solving \eqref{ODEFrho} in time $\tau$, tells us precisely how time $t$ and nonlinearly rescaled time $\tau$ are related.
Zeros of $\rho$ may lead to finite times $t$ along unbounded paths of $\tau\in\C$.
At nonzero $\mathbf{F}$, this describes finite-time blow-up in $t$, at singularities of $1/\rho\cdot\mathbf{F}$ generated by the denominator $\rho$\,.
Zeros of $\rho$ may also ``outsource'' zero denominators of $\mathbf{F}=(f,g)$ in ODE \eqref{ODEfg}; see section \ref{RecC2}.
This allows us to holomorphically extend nontrivial holomorphic foliations, such that the remaining singularity locus of $\mathbf{F}=0$ is of complex codimension (at least) two, after outsourcing.
See \cite{Ilya}, section 2D.
See the proof of proposition \ref{proploop}, and sections \ref{LinC2}, \ref{RecC2}, for first examples.
In the projective compactification $\CP^2$ of $\C^2$, we will discuss many more examples involving blow-up of $x$ or $y$ in finite time.

For later use, let $dt_\ell=\rho_\ell d\tau$ denote the effect of several Euler multipliers $\rho_\ell$\,. 
Then
\begin{equation}
\label{rho12}
\frac{dt_{\ell_2}}{dt_{\ell_1}}=\rho_{\ell_2\ell1}:=\frac{\rho_{\ell_2}}{\rho_{\ell_1}}\,.
\end{equation}
Obvious consequences, for example, are $\rho_{31}=\rho_{32}\rho_{21}$ and $\rho_{12}=1/\rho_{21}$\,.

\subsection{Complex foliations in projective coordinates}\label{Fol}

The advantageous use of homogeneous projective coordinates $[\xi\!:\!\eta\!:\!\zeta]$, as in \eqref{xyz}, to extend polynomial ODEs from $\C^2$ to compact $\CP^2$ is fraught with peril.
Let us abbreviate $\fx:=(\xi,\eta,\zeta)\in\C^3\setminus\{0\}$ and let $\fF{:}\ \fx\mapsto\fF(\fx)\in\C^3$ be any nontrivial \emph{homogeneous polynomial} of degree $m$, i.e.
\begin{equation}
\label{fFm}
\fF(\sigma\fx)=\sigma^m\fF(\fx),
\end{equation}
for all $\fx\in\C^3$ and all $\sigma\in\C$.
Then $\fF$ induces a projective map $[\fF]{:}\ \CP^2\rightarrow \CP^2$.
However, neither the homogeneous ODE vector field
\begin{equation}
\label{ODEfF}
\fx'(\tau)=\fF(\fx)\,,
\end{equation} 
nor the ODE flow in time $\tau$, respect the proportions defining $\CP^2$, unless $m\!=\!1$.
We concede
\begin{equation}
\label{fxfy}
\fx(\tau) \textrm{ solves } \eqref{ODEfF} \quad\Leftrightarrow \quad \fy(\tau):=\sigma\,\fx(\sigma^{m-1}\tau) \textrm{ solves } \eqref{ODEfF}\,.
\end{equation}
However, this defines an ODE flow on the equivalences classes $[\fx]\in\CP^2$ in the scaling-equivariant linear case $m\!=\!1$, only. See also section \ref{LinProj}.

\begin{rem}\label{remEuler}
Complex Euler multipliers $\rho$ and complex foliations $\omega$ come to our rescue.
The next proposition shows how polynomial foliations $\omega(\mathbf{F})=0$ of $\C^2\setminus \mathbf{F}^{-1}(0)$ remain well-defined on the compactification $\CP^2$ of $\C^2$.
In other words, the ODE \eqref{ODEF} may not define a flow on $\CP^2$.
But it defines a foliation, with complex-time orbits as leafs $\fL$, in various local coordinates of $\CP^2$.
The foliations, in turn, are given by the kernel of a Pfaffian differential 1-form $\omega$.
The kernel is invariant under multiplication by nonzero Euler multipliers $\rho$.
The zeros of complex Euler multipliers will contribute crucially to extensions of the foliations, and to the detection of blow-up loops $\gamma^t$ in ODEs, along the Masuda PDE paradigm.
To study blow-up of \eqref{ODEfg} towards $u\!=\!1/x\!=\!0$ or $v\!=\!1/y\!=\!0$, we will therefore discuss many examples of this \emph{complex projective compactification} of foliations for polynomial ODEs, in the following sections.
\end{rem}

To prepare, let $m\geq1$ be the maximal degree of the complex polynomial $\mathbf{F}=(f,g)$.
Define the $m$-homogeneous polynomial $\fF=(\ff,\fg,\fh)$ in $\fx=(\xi,\eta,\zeta)$ by
\begin{align}
\label{ff}
   \ff(\xi,\eta,\zeta)&:=\xi\cdot\zeta^{m-1}+\zeta^m\cdot f(\xi/\zeta,\eta/\zeta)\,,   \\
\label{fg}
   \fg(\xi,\eta,\zeta)&:=\eta\cdot\zeta^{m-1}+\zeta^m\cdot g(\xi/\zeta,\eta/\zeta)\,,  \\  
\label{fh}
   \fh(\xi,\eta,\zeta)&:=\zeta\cdot\zeta^{m-1}\,.
\end{align}
For example, the sum $f(x,y):=\sum_{j,k} f_{jk}\,x^jy^k$ over $0\leq j+k\leq m$ implies
\begin{equation}
\label{ffdetail}
\ff(\xi,\eta,\zeta)= \xi\zeta^{m-1}+\sum_{j,k} f_{jk}\,\xi^j\eta^k\zeta^{m-j-k}\,,
\end{equation}
and similarly for $\fg$.

\begin{rem}\label{short}
To avoid pedantic notation here and below, cluttered by mere time transformations, we will write $u$ for $\tilde u$ in expressions like $\tilde u(t_1):=u(\tau)$. 
For derivatives, similarly, we will write $\dot u(t_1)$, shorthand for $\tfrac{d}{dt_1}\tilde u(t_1)$, and $\dot u(\tau)$ for $u'(\tau)=\tfrac{d}{d\tau}u(\tau)$.
Indeed, the arguments $t_1$ and $\tau$ of $u$ already entail the full notational regalia. 
\end{rem}

\begin{figure}[t]
\centering \includegraphics[width=0.71\textwidth]{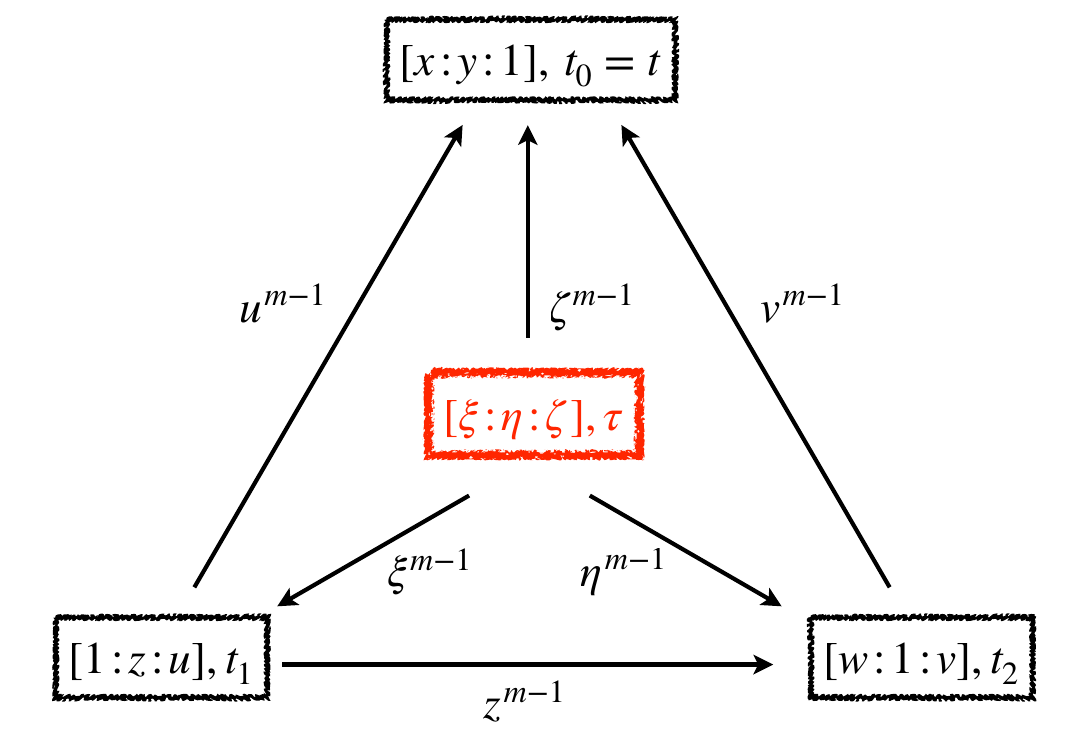}
\caption[Euler multipliers]{\emph{
Euler multipliers $\rho$ associated to homogeneous projective coordinates $[\xi\!:\!\eta\!:\!\zeta]=[x\!:\!y\!:\!1]=[1\!:\!z\!:\!u]=[w\!:\!1\!:\!v]$ of the projective compactification $\CP^2$ of $\C^2$.
Time variables in \eqref{ODEfF}, \eqref{ODExy}-\eqref{ODEvw} are $\tau$ and $t=t_0, t_1, t_2$\,, respectively.
Multipliers $\rho$ are indicated in the direction of arrows.
For example, $dt=\zeta^{m-1}d\tau=u^{m-1}dt_1$\,, and so on.
See \eqref{ttl} and \eqref{tt1}, \eqref{tt2}.
Note multiplicativity \eqref{rho12}, and hence commutativity, in the diagram.
}}
\label{figrho}
\end{figure}

\begin{prop}\label{propuvw}
Let the foliation $\omega(\mathbf{F})=0$ be defined on the nonequilibrium set of $\mathbf{F}=(f,g)$, by \eqref{omega}, \eqref{ODEfg}. See \eqref{omegaF=0}.
Consider the homogeneous projective coordinates $[\fx]=[\xi\!:\!\eta\!:\!\zeta]=[x\!:\!y\!:\!1]=[1\!:\!z\!:\!u]=[w\!:\!1\!:\!v]$ on $\CP^2$; see \eqref{xyz}.\\
Then the foliation $\ker\omega$ extends to the nonequilibrium set of \eqref{ODEfF} on $\CP^2$, as follows.
With an Euler multiplier $\rho_0:=\zeta^{m-1}\neq0$, ODE \eqref{ODEfg} with time $t\!=\!t_0$ reads
\begin{equation}
\label{ODExy}
\begin{aligned}
    \dot x(t)=\rho_0^{-1}\dot x(\tau)&= f(x,y)\,,   \\
     \dot y(t)=\rho_0^{-1}\dot y(\tau)&= g(x,y)\,,
\end{aligned}
\end{equation}
in the above shorthand. The foliation is extended by the foliation of the polynomial ODE
\begin{equation}
\label{ODEuz}
\begin{aligned}
     \dot u(t_1)=\rho_1^{-1}\dot u(\tau)&= -u\,f_1(u,z)\,,   \\
    \dot z(t_1)=\rho_1^{-1}\dot z(\tau)&= -z\,f_1(u,z)+g_1(u,z)=:p(u,z)\,.
\end{aligned}
\end{equation}
with Euler multiplier $\rho_1:=\xi^{m-1}$.
Analogously, an Euler multiplier $\rho_2:=\eta^{m-1}$ leads to the extension of the foliation by the polynomial ODE
\begin{equation}
\label{ODEvw}
\begin{aligned}
     \dot v(t_2)=\,\rho_2^{-1}\dot v(\tau)\,&=-v\,g_2(v,w)\,,   \\
     \dot w(t_2)=\rho_2^{-1}\dot w(\tau)&=-w\,g_2(v,w)+f_2(v,w)=:q(v,w)\,.  
\end{aligned}
\end{equation}
Here we have used the following transformations, which act linearly on polynomials $f,g$\,:
\begin{equation}
\label{f12}
f_1(u,z) :=u^mf(1/u,z/u)\,,\quad f_2(v,w):= v^mf(w/v,1/v)\,,
\end{equation}
and analogously for $g$.
Similarly
\begin{equation}
\label{f21}
f_1(u,z)=z^mf_2(u/z,1/z)\,,\quad f_2(v,w)=w^mf_1(v/w,1/w), \ \textrm{etc}.
\end{equation}
The time variables $t\!=\!t_\ell$ in  \eqref{ODExy}-\eqref{ODEvw} with labels $\ell=0,1,2$, respectively, are related to $\tau$ in \eqref{ODEfF} by
\begin{equation}
\label{ttl}
dt_\ell=\rho_\ell\, d\tau\,.
\end{equation}
See \eqref{ODEFrho}, \eqref{rho12} and figure \ref{figrho}.
\end{prop}

\begin{proof}
Consider ODE \eqref{ODEfF} for the $m$-homogeneous vector field $\fF$ defined by \eqref{ff}-\eqref{fh}.
We first verify that ODE \eqref{ODEfF} describes the prescribed original vector field \eqref{ODExy}, alias \eqref{ODEfg}, in time $t_0=t$.
Indeed,  $[\xi\!:\!\eta\!:\!\zeta]=[x\!:\!y\!:\!1]$ implies $x=\xi/\zeta$ and $y=\eta/\zeta$.
Therefore \eqref{ODEfF} implies
\begin{equation}
\label{xpf}
\begin{aligned}
\dot x(\tau)&=\zeta^{-1}\left(\dot \xi(\tau)-x \dot \zeta(\tau)\right)=\zeta^{-1}\left(\ff(\xi,\eta,\zeta)-x\fh(\xi,\eta,\zeta)\right)=\\
&=\zeta^{m-1}\left(x+f(x,y)-x\right)= \rho_0\,f(x,y)=\rho_0\,\dot x(t)\,.
\end{aligned}
\end{equation}
Replacing $x,\xi,\ff,f$ by $y,\eta,\fg,g$ completes the proof of claim \eqref{ODExy}.
The nonzero Euler multiplier $\rho_0=\zeta^{m-1}$ satisfies \eqref{ttl}, for $\ell\!=\!0$.

We verify the second equation of claim \eqref{ODEuz} next, for $[\xi\!:\!\eta\!:\!\zeta]=[1\!:\!z\!:\!u]$, i.e.\ $u=\zeta/\xi=1/x$ and $z=\eta/\xi=y/x$. 
Therefore \eqref{ODEfF} implies
\begin{equation}
\label{zpf}
\begin{aligned}
\dot z(\tau)&=\xi^{-1}\left(\dot \eta(\tau)-z \dot \xi(\tau)\right)=
\xi^{-1}\left(\fg(\xi,\eta,\zeta)-z\ff(\xi,\eta,\zeta)\right)=\\
&=\xi^{m-1}\left(\fg(1,z,u)-z\ff(1,z,u)\right)=\\
&=\rho_1\left(zu^{m-1}{+}u^mg(1/u,z/u)-z\left(u^{m-1}{+}u^mf(1/u,z/u)\right)\right)=\\
&= \rho_1 \left(g_1(u,z)-zf_1(u,z)\right)=\rho_1\,\dot z(t_1)\,.
\end{aligned}
\end{equation}
All remaining claims of \eqref{ODEuz}, \eqref{ODEvw} are completely analogous.

Up to different Euler multipliers $\rho_\ell$ in \eqref{ttl}, the three ODEs \eqref{ODExy}-\eqref{ODEvw} are versions of one and the same ODE \eqref{ODEfF}.
Therefore they define one and the same foliation $\ker\,\omega$ on any intersection of their domains.
This proves the proposition.
\end{proof}

Our main interest are blow-up solutions $(x(t),y(t))$ of ODE \eqref{ODExy}, which reach $x\!=\!\infty$ or $y\!=\!\infty$ in finite original time $t\!=\!t_0$\,.
Projectively, blow-up is resolved as the Riemann sphere $\widehat{\C}=\CP^1$ of coordinates $[\xi\!:\!\eta\!:\!0]$, i.e.\  by projective coordinates $u\!=\!0$ or $v\!=\!0$ in \eqref{ODEuz} and \eqref{ODEvw}, respectively.
We therefore collect the relevant time transformations.
\begin{cor}\label{tt}
In the setting of proposition \ref{propuvw}, the time transformations from $t_1$ in \eqref{ODEuz} and $t_2$ in \eqref{ODEvw} to $t\!=\!t_0$ in \eqref{ODExy} are given by the Euler multipliers
\begin{align}
\label{tt1}
    dt&=u^{m-1}dt_1\,,   \\
\label{tt2}
    dt&=v^{m-1}dt_2\,.
\end{align}
\end{cor}
\begin{proof}
Since $u=\zeta/\xi$ and $v=\zeta/\eta$, relations \eqref{ttl} with the rule \eqref{rho12} of compositions imply $dt=\rho_0\rho_1^{-1}\,dt_1=(\zeta/\xi)^{m-1}\,dt_1=u^{m-1}\,dt_1$\,. 
See figure \ref{figrho}.
Similarly, $dt=\rho_0\rho_2^{-1}\,dt_2=(\zeta/\eta)^{m-1}\,dt_2=v^{m-1}\,dt_2$\,, as claimed.
\end{proof}

\subsection{Blow-up loops and blow-up stars in complex and real time}\label{LoopStar}

We return to Masuda detours around blow-up, in the general setting of polynomial systems \eqref{ODExy} and in the language of foliations.
Slightly more generally than Masuda, we start from a complex reference solution $(x,y)(t)\in\C^2$ which blows up, say, for $t\rightarrow T=0$.
Our compactification in $\CP^2$ refers to blow-up coordinates $u\!=\!1/x,\ z\!=\!y/x$ and $v\!=\!1/y,\ w\!=\!x/y$ in rescaled times \eqref{ttl}, as in the previous section \ref{Fol}.
Due to invariance of infinity, represented in $\CP^2$ by $u\!=\!0$ or $v\!=\!0$, blow-up requires convergence to equilibrium, $(u,z)\rightarrow(0,e)$ or $(v,w)\rightarrow(0,e')=(0,1/e)$, respectively.
Skipping an analogous discussion of $v,w$, we may therefore proceed in coordinates $(u,z)$, assuming $z\!=\!y/x$ remains finite as we approach blow-up.
Let $(u,z)\in\fL\subset\CP^2\setminus\{u\!=\!0\}$ denote the local leaf of the blow-up solution, i.e.\ the orbit in complex time.
Then any Masuda detour must remain in that leaf $\fL$.

Consider any local circumvention of the blow-up time $t\!=\!T\!=\!0$, by a closed curve $\gamma^t(s):=t(s)\in\Co,\ 0\leq s\leq 2\pi$, as in \eqref{gammat}.
We assume nonzero winding number $\fw(\gamma^t)$ around $T\!=\!0$.
Without loss of generality, we may choose left winding $\fw(\gamma^t)>0$.
Let the curve $\gamma^{uz}\subset\fL$ be defined by the solution $\gamma^{uz}(s):=(u(t(s)),z(t(s)))$ of \eqref{ODEuz}, \eqref{tt1} along the curve $\gamma^t$ in the local blow-up leaf $\fL$.
See remarks \ref{remEuler}, \ref{short}, and corollary \ref{tt}, for passage from time $t_1$ in $(u(t_1),z(t_1))$ to original time $t\!=\!t_0$\,, and vice versa.

\begin{defi}\label{defloop}
We call the above pair $(\gamma^t,\,\gamma^{uz})$ a \emph{blow-up loop}, if the lifted solution curve $\gamma^{uz}\subset\fL$ closes up, with vanishing discrepancy as defined in \eqref{discrepancy}.
In other words,
\begin{equation}
\label{discrepancyuz}
\gamma^{uz}(s)\Big|_{s=0}^{s=2\pi}=0\,.
\end{equation} 
We call a blow-up loop \emph{minimal}, if the winding number $\fw(\gamma^t)>0$ is minimal among all blow-up loops in the local leaf $\fL$.
\end{defi}

Minimal blow-up loops are only determined up to free homotopy of the closed loop $\gamma^t$, locally in a small punctured disk around $t\!=\!T$.
This follows from local commutativity \eqref{flow} of ODE flows in complex time.

The loop $\gamma^t$ may be self-intersecting.
For example, we may always cover any blow-up loop $(\gamma^t,\gamma^{uz})$ repeatedly, by periodic extension of the parametrization $s\mapsto\gamma^t(s)$ to $0\leq s\leq 2n\pi$.
Minimality excludes such iterated loops.

By invariance of $\{u\!=\!0\}$, the $z$-plane $u\!=\!0$ cannot intersect $\gamma^{uz}$.
Therefore the projections $\fL^u,\,\gamma^u$ of $\fL,\,\gamma^{uz}$ onto $u$ are never zero: $\gamma^u\subset\fL^u\subseteq\Co $.
For $\iota\in\{t,u\}$, we abbreviate winding numbers as $\fw_\iota:=\fw(\gamma^\iota)$.

\begin{prop}\label{proploop}
Fix polynomial degree $m\geq2$.
Let  $(u(t),z(t))$ denote a blow-up solution of \eqref{ODEuz}, which converges towards an equilibrium $(u,z)=(0,e)$ along a path $t\rightarrow T=0$.
Assume a blow-up loop $(\gamma^t,\gamma^{uz})$ in complex time $\gamma^t(s)=t(s)\neq T$ around blow-up time $T$, sufficiently close to $(t,u,z)=(T,0,e)$.
See definition \ref{defloop}.
Also assume nonvanishing semisimple eigenvalues $\lambda_1,\,\lambda_2$ of the linearization at $(0,e)$.\\
Then the integer winding numbers $\fw_\iota\in\Z$ of the loop $t(s)\in\gamma^t\subset\C\setminus\{T\}$ and of the projected solution loop $u(t(s))\in\gamma^u\subset\Co,\ 0\leq s\leq 2\pi$, are related by
\begin{equation}
\label{wtu}
\fw_t=(m-1)\,\fw_u\,.
\end{equation}
In particular, nonzero winding of either curve implies nonzero winding of the other.\\
The same statement holds, analogously, for the windings $\fw_t$ and $\fw_v$ of blow-up loops $(\gamma^t,\gamma^{vw})$ in \eqref{ODEvw}.
\end{prop}

\begin{proof}
Up to substitution of $u,z$ by $v,w$, it is sufficient to prove the claim for the blow-up loop $(\gamma^t,\gamma^{uz})$.
For convenience and without loss of generality, we may shift $z$ in \eqref{ODEuz} such that $e\!=\!0$ with diagonal linearization.
We obtain
\begin{equation}
\label{uzdiag}
\begin{aligned}
     \dot u(t_1)&=u(\lambda_1+\ldots)\,,   \\
    \dot z(t_1)&=z(\lambda_2+\ldots)+u^2(c_2+\ldots)\,,
\end{aligned}
\end{equation}
with some complex coefficient $c_2$\,, and up to higher order terms in $u,z$.
The Euler multiplier $\rho(u)=u^{m-1}$ in \eqref{tt1}, \eqref{uzdiag} then provides a local expansion
\begin{equation}
\label{dtdunondeg}
dt=u^{m-1}dt_1=u^{m-1}\,\frac{dt_1}{du}\,du = \tfrac{1}{\lambda_1}\,u^{m-2}\,(1+\ldots)\,du\,.
\end{equation} 
Integration along the local blow-up path $t\rightarrow T$ yields
\begin{equation}
\label{tuexpan}
t-T= \tfrac{1}{\lambda_1(m-1)}\, u^{m-1}(1+\ldots)\,,
\end{equation}
with higher order terms in both $u$ and $z$, viz.~uniformly small perturbations depending on $t-T$ along the local path. 
By \eqref{flow}, local flows in complex time commute.
Therefore \eqref{tuexpan} extends to the loop $\gamma^t$.
But this only shows that the curve $\gamma^t$ performs \emph{almost} $m-1$ complete windings around blow-up time $T$, upon a single loop of small $\gamma^u\subset\Co$.
Since we have \emph{assumed} a discrepancy-free blow-up loop $(\gamma^t,\gamma^{uz})$, on the other hand, both loops must close, at some finite integer winding numbers.
Therefore expansion \eqref{tuexpan} is sufficient to prove claim \eqref{wtu}.
\end{proof}

Our main results, already surveyed in the overview of section \ref{Over} and presented in section \ref{B}, establish the existence of blow-up loops.
The following sections discuss many specific examples.
For the excluded case $m\!=\!1$, at $x\!=\!y\!=\!0$, see the planar linear case of section \ref{LinC2}.
For slight generalizations of the explicit scalar Riccati case $m\!=\!2$ in \eqref{Ric0}, see also section \ref{Ric} and reduction \eqref{Riclin} to the scalar linear case.
In section \ref{RecC2}, we discuss a reciprocally linear system.
We address the general $m$-homogeneous case, in section \ref{C2m}, and the Hamiltonian case of plane complex algebraic curves in section \ref{Ham}.
Section \ref{Fuji} returns to a case study for the crude ODE caricatures \eqref{ODEb}, \eqref{ODEa} of our original Masuda example with $m\!=\!2$: the quadratic Fujita PDEs \eqref{PDEw}, \eqref{PDEb}.
See section \ref{Con} for a less technical summary.

Caution is in place, however, concerning continuation after ``complete'' blow-up -- even in real time.
In the sections to follow, we encounter many examples where the blow-up loop $(\gamma^t,\,\gamma^{uz})$ closes up nicely,  without any discrepancies, perhaps after some iterated cycling of $\gamma^t$ around $t\!=\!T$. 
Any continuation beyond blow-up, however, may require complex $u$ if we insist on real time $t$.
Or complex time, if we insist on real $u$.
For explicit cubic examples, $m\!=\!3$, see figures \ref{figm34}(b) and \ref{figm234}(center).

The phenomenon of blow-up loops $(\gamma^t,\,\gamma^{uz})$ with winding numbers $\fw_t = (m\!-\!1)\fw_u>1$ is closely related to nonuniqueness of complex continuation after ``complete'' blow-up, in real time $t$.
See figure \ref{figm234} for illustration of $\fw_u=1$ and degrees $m=2,3,4$.
Indeed, a circular multiple loop $\gamma^t$ performs at least $\fw_t$ crossings of the positive real time axis $t>T=0$, after blow-up time $T$,  before the loop $\gamma^{uz}$ first closes up again.
In a minimal blow-up loop, we heuristically also expect $\fw_t$ distinct complex-valued blow-down continuations $(u(t),z(t))$, in the same local leaf $\omega(\fL)=0$, \emph{after} blow-up and in real time $t>0$.
Arriving at blow-up from a particular real-time solution $0>t\nearrow T=0$, we may think of the sparkling fireworks radiating from a rocket after, well, ``blow-up''.
For real times $0>t\nearrow T$ \emph{preceding} blow-up time $T\!=\!0$, we likewise expect $\fw_t-1$ distinct further complex-valued companions $(u(t),z(t))\in\fL$, in the same leaf $\fL$ as our original blow-up solution, which all collide and blow up at $t\!=\!T$, synchronously.

\begin{defi}\label{defstar}
Let $(\gamma^t,\gamma^{uz})$ be a minimal blow-up loop in the local leaf $\fL$, with winding number $0<\fw_t:=\fw(\gamma^t)\in\Z$.
A  \emph{real-time blow-up star} is the collection of all $\fw_t$ complex blow-up branches of solutions, and all $\fw_t$ blow-down continuations, in real time $\pm (t-T) \searrow 0$ respectively, which alternate in the leaf $\fL$.\\
The stars include, and identify with, their companions of blow-up loops $(\gamma^t,\gamma^{vw})$ in blow-up coordinates $v,w$.
\end{defi}

The original results \cite{Masuda1, Masuda2} by Masuda observe absence of blow-up loops for PDE perturbations of the spatially homogeneous Riccati case $\mathbf{w}(t,\bxi)=x(t)$.
The ODE Riccati case of quadratic scalar nonlinearities itself features unique real continuation, due to simple winding $\fw(\gamma^t)=1$ of its blow-up loop; see the explicit solution \eqref{Ric0} and section \ref{Ric}.
In section \ref{Scalar}, we encounter general winding $\fw(\gamma^t)=m-1$ for scalar nonlinearities $\dot x=P(x)$ of polynomial degree $m$.
Already reciprocally linear $\C^2$ cases provide real-time blow-up stars with any finite number of branches, at suitable resonances; see section \ref{RecC2}.

In all these cases, blow-up loops reflect the local discrepancy-free behavior of complex foliations $\omega\!=\!0$ and their complex leafs $\fL$ near blow-up equilibria $(u,z)=(0,e)$.
We therefore study these complex foliations in some detail, in the linear example of section \ref{LinC2} and in the central section \ref{B} on analytic ``linearization at infinity''.
In particular, see theorem \ref{thmloop} in section \ref{PoiMas} on the relation among minimal blow-up loops, blow-up stars, the branching of Riemann surfaces, and meromorphically branched blow-up. 
In conclusion, it is the notion of complex foliations which decides on the success, or failure, of discrepancy-free blow-up loops in the Masuda sense -- and on the complex or real continuation of solutions, in real time, by their associated real-time blow-up stars.

\subsection{Technical outline} \label{Out}

For a nontechnical overview recall section \ref{Over}.
The previous section \ref{LoopStar} has introduced blow-up loops and blow-up stars as our central concepts for the analysis of ODE blow-up in the complex plane; see definitions \ref{defloop} and \ref{defstar}.
Methodically, our approach is largely based on \cite{Ilya}. 
For ample illustration, but also for preparation of proofs, the eight sections \ref{Ric}-\ref{RecC2} and \ref{C2m}-\ref{Con} are dedicated to specific classes of examples.
In section \ref{Ric}, and as a warm-up towards blow-up, we revisit autonomous scalar quadratic Riccati ODEs in the complex domain.
Our previous work on scalar complex polynomial ODEs, which remains relevant within complex one-dimensional separatrices of blow-up saddles, is reviewed in section \ref{Scalar}.
As an indispensable preparation for our analysis of nonlinear blow-up by ``linearization at infinity'', section \ref{LinC2} studies autonomous linear ODEs in $\C^2$ and $\CP^2$.
A first application of complex Euler multipliers, and the associated language of complex foliations, is blow-up in reciprocally linear ODEs and their reduction to the linear case; see section \ref{RecC2}.
Our main results, in section \ref{B}, study blow-up loops and blow-up stars near blow-up equilibria, in projective coordinates.
See theorem \ref{thmloop}.
The results amount to, and rigorously justify, analytic ``linearization at infinity''.
In particular, linearization at nonresonance provides a precise quantitative description of the actual Masuda discrepancies \eqref{discrepancy} ``after'' blow-up, and their vanishing along blow-up loops.
Section \ref{C2m} addresses blow-up in $m$-homogeneous polynomial ODEs.
Section \ref{Ham} comments on the general case of polynomial Hamiltonian ODEs, where the linearization approach of section \ref{B} fails.
Instead we resort to complex energy surfaces, as algebraic leafs, i.e. as compact Riemann surfaces.
The quadratic Weierstrass and the cubic Duffing pendulum are simple examples.
In section \ref{Fuji} we return to the original PDE motivation of sections \ref{RBlow-up}-\ref{C2}\,: crude two-mode Galerkin caricatures \eqref{ODEb}, \eqref{ODEa} of the Masuda paradigm.
Section \ref{Fuji} follows the linearization approach of section \ref{B} and studies these ODEs ``at infinity'' $[\xi\!:\!\eta\!:\!0]\in\CP^1\subset\CP^2$, via the projective compactifications \eqref{ODExG}-\eqref{ODEwG2} in $\CP^2$.
We highlight and compare spatially symmetric heterogeneous quadratic coefficients $\mathbf{b}(\xi)$ in \eqref{PDEb}, and certain spatially antisymmetric modifications.
The conclusions of section \ref{Fuji} on Masuda detours, blow-up loops and real-time blow-up stars are summarized, once again, in section \ref{Con}.
Section \ref{PDElift} admits a broader class of PDE variants of \eqref{PDEb}: 
complex parabolic reaction-drift-diffusion equations on the circle.
This class encompasses all polynomial examples of complex planar ODE blow-up of our present paper and lifts them to a specific PDE context.
We conclude with comments on previous literature and further perspectives of our sketch, as summarized at the beginning of section \ref{Com}.

\section{Example: Riccati equations}\label{Ric}

\begin{figure}[t]
\centering \includegraphics[width=0.86\textwidth]{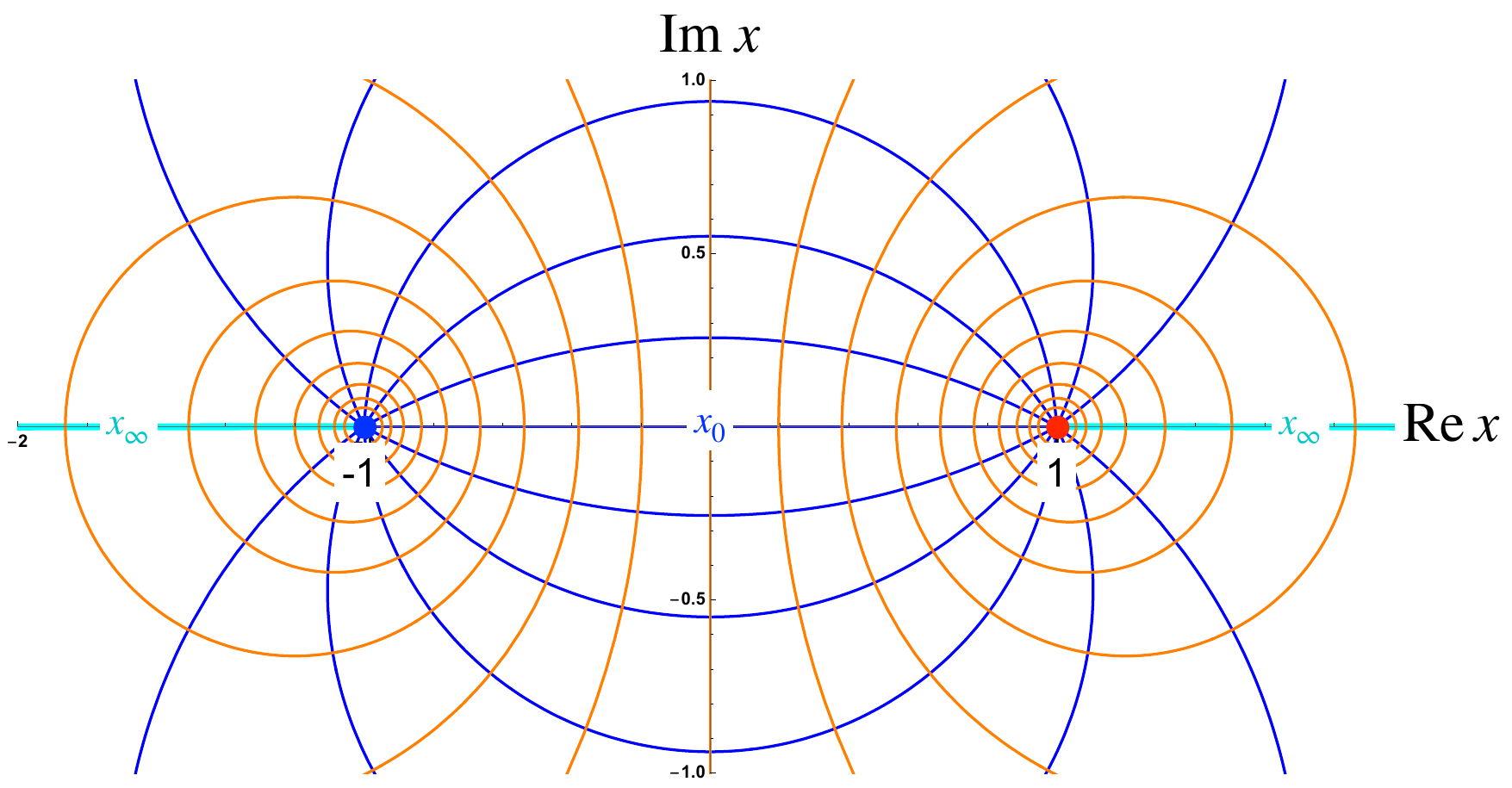}
\caption[The complex Riccati equation]{\emph{
Phase portrait of the quadratic complex Riccati equation $\,\dot x=x^2-1$, i.e.\ \eqref{ODERic} with  $a=1,\ e_1=+1,\ e_2=-1$.
In real time $t$, the equilibrium $e_1=+1$ is a source (red) and $e_2=-1$ (blue) is a sink.
All other real-time orbits (blue) are heteroclinic from $+1$ to $-1$; compare \eqref{het}.
They foliate the cylinder leaf $\fL=\widehat{\C}\setminus\{\pm 1\}$ into circular segments with midpoints on the imaginary axis.
The blow-up and blow-down real-time orbits $x_\infty\subset\C^2$ (cyan) are parts of a single heteroclinic orbit through $x\!=\!\infty$, on the Riemann sphere $x\in\widehat{\C}$.
Hence blow-up can be circumnavigated by a blow-up loop.
Note that the definition and selection of the cyan orbits is not invariant under Möbius transformations.
\\
Nonstationary orbits in imaginary time (orange), in contrast, provide an iso-periodic foliation of the complex leaf $\fL$ with minimal period $\pi$.
The periodic orbits are nested around $\pm 1$.
On the Riemann sphere, this includes the iso-periodic orbit on the imaginary axis, which experiences $\C^2$ blow-up and blow-down in finite imaginary time.
Blow-up and blow-down in imaginary time are related to real heteroclinicity of $x_0$ in real time, just as in the PDE context of \cite{FiedlerFila}; see section \ref{PDE}.
Since the flow $t\mapsto\Phi^t(x_0)$ is conformal, for any fixed $x_0\in\widehat{\C}\setminus\{\pm1\}$, the blue and orange circle families are mutually orthogonal.
}\cite{FiedlerShilnikov}}
\label{figRiccati}
\end{figure}

We have mentioned the purely quadratic Riccati ODE \eqref{Ric0} as the spatially homogeneous version of the Masuda PDE paradigm \eqref{PDEw}.
All nonstationary real-time orbits $x(t)$ in the Riemann sphere $\widehat{\mathbb{C}}$ are homoclinic to the algebraically double equilibrium $x=e_1=e_2=0$, for $|t|\rightarrow\infty$.
In $\C\setminus\R$, all trajectories are circles, tangent to the real axis at the origin.\!\!~
\footnote{If you never did: draw and verify this, e.g., based on the explicit solutions \eqref{Ric0}.
In class, it may serve as a nice example for lower-semicontinuity of blow-up time $T\leq+\infty$. }
We now illustrate some of our concepts for the ``general'' quadratic Riccati equation 
\begin{equation}
\label{ODERic}
\dot{x}=a(x-e_1)(x-e_2)
\end{equation}
with distinct equilibria $e_1, e_2\in\C$ and some coefficient $a\in\Co$.
The real-time case $e_1=0<a,\ e_2=K,\ x\geq0$ has acquired some textbook fame in elementary population models and economics with ``carrying capacity'' $K>0$.
See figure \ref{figRiccati} for an illustration of complex orbits $x(t)$ with $a=e_1=1,\ e_2=-1$, in real and imaginary time $t$.

The group of biholomorphic automorphisms of the Riemann sphere $\widehat{\mathbb{C}}$, alias the \emph{Möbius group}, consists of the fractional linear  \emph{Möbius transformations}
\begin{equation}
\label{Mobius}
x\mapsto u:=\frac{\mathfrak{a}x+\mathfrak{b}}{\mathfrak{c}x+\mathfrak{d}}\quad\mathrm{with}\quad
\begin{pmatrix}
    \mathfrak{a}  &  \mathfrak{b}  \\
    \mathfrak{c}  &  \mathfrak{d}
\end{pmatrix}
\in \mathrm{PGL}(2,\mathbb{C})\cong \mathrm{SL}(2,\mathbb{C})/\{\pm \mathrm{Id}\}\,.
\end{equation} 
The group is generated by \emph{shifts} $u=x+\mathfrak{b}$, \emph{complex scalings} $\mathfrak{b}=\mathfrak{c}=0$ of $x$ by nonzero $a=\mathfrak{a}^2=r\exp(\mi\phi)$, and the involutary \emph{inversion} $u\!=\!1/x$.

Let us map $x\!=\!\infty$ to $u\!=\!0$ by the inversion $u\!=\!1/x$, for example.
Then
\begin{equation}
\label{udot=1}
\dot u=-a(1-e_1u)(1-e_2u)=-a+\ldots\neq 0\,,
\end{equation}
near $u\!=\!0$.
In particular, $t=-a^{-1}u+\ldots$ is a local biholomorphism: any small complex loop $\gamma^u$ of $u$ around $u\!=\!0$ becomes a loop $\gamma^t$ of $t$ around the blow-up time $t\!=\!T\!=\!0$, and vice versa.
``Blow-up'' at $x\!=\!\infty$ can therefore be circumnavigated by an explicit blow-up loop $(\gamma^t,\gamma^u)$, as in definition \ref{defloop}.
Moreover, the resulting global flow $\Phi^t$ extends to the Riemann sphere $x\in\widehat{\C}=\CP^1$.
Since $\Phi^t$ consists of automorphisms, the flow property \eqref{flow} identifies $\Phi^t$ as a 1-parameter subgroup of Möbius transformations -- without any further calculation.

A Möbius map of the distinct equilibria $x=e_1, e_2$ to $u=0,\infty$, respectively, preserves the linearizations at equilibria and therefore transforms \eqref{ODERic} to the linear variant
\begin{equation}
\label{Riclin}
\dot u= a(e_1-e_2) u\,,
\end{equation}
now with the equilibrium $x\!=\!e_2$ at $u\!=\!\infty$.
The associated exponential global flow $\Phi^t$ therefore identifies all 1-parameter subgroups of the Möbius group, up to conjugation, with the exception of the purely quadratic case \eqref{Ric0}, viz.~the limiting case $e_1=e_2=0$ of \eqref{Ric0}.

Möbius transformations preserve the set of lines and circle segments, jointly. 
Orbit equivalence of $\dot u=u$ to $\dot x=x^2-1$, by a Möbius transformation, therefore establishes the circular segments and circles of figure \ref{figRiccati}.
We illustrate the flows in real time (blue), and in imaginary time (orange). 
In real time, we obtain a unique source $e_1=+1$ (red) and sink $e_2=-1$ (blue), but no poles. 
All other real-time orbits $x(t)$ (also blue) are \emph{heteroclinic} from $+1$ to $-1$, i.e. $\lim x(t)=\mp 1$ for $t\rightarrow\pm\infty$.
In symbols:
\begin{equation}
\label{het}
x(t)\,{:}\quad +1\leadsto -1\,.
\end{equation}
On the real axis in the Riemann sphere $\widehat{\C}$, we encounter the heteroclinic blow-up and blow-down orbit $x_\infty$ through $x\!=\!\infty$ (cyan).

In imaginary time $t$, the two equilibria become \emph{Lyapunov centers}: the eigenvalue $\lambda$ of the linearization is purely imaginary and, hence, their neighborhoods are foliated by periodic orbits of minimal imaginary period $2\pi/\lambda$ in $t$.
In our case, all nonstationary orbits (orange) are iso-periodic of minimal imaginary period $\pm\mi\pi$.
They foliate the cylinder $\widehat{\C}\setminus\{e_1\,,\ e_2\}$.
Because the flow $t\mapsto\Phi^t(x_0)$ of Möbius transformations is conformal, for any fixed $x_0$ in the single complex leaf $\fL=\widehat{\C}\setminus\{e_1,e_2\}$, the heteroclinic and iso-periodic foliations in real and imaginary time $t$ are mutually orthogonal, just as the real and imaginary time directions themselves are.

\section{Example: Scalar polynomial ODEs in real and complex time}\label{Scalar}

Our next example are scalar complex polynomial ODEs of degree $m\geq 2$,
\begin{equation}
\label{ODEP}
\dot x(t)=P(x)=\sum_{j=0}^m \,p_{m-j}\,x^j=p_0\,(x-e_1)\cdot \ldots \cdot (x-e_m)\,.
\end{equation}
For later reference in sections \ref{z=e}, \ref{Ls}, \ref{C2m}, \ref{Hamm}, and \ref{Con}, we summarize our previous results on the global dynamics of \eqref{ODEP} in the Riemann sphere $\widehat{\C}$.
See \cite{FiedlerShilnikov, FiedlerYamaguti} and the references there for complete details and related results.
Our results hold for generic complex polynomials $P$. 
Explicit solution by separation of variables, it turned out, does \emph{not} tell it all:
geometric issues like the blow-up loop and blow-up star of the degenerate equilibrium ``at infinity'', as well as the global geometry of the associated foliations by branched Riemann surfaces, or the mere combinatorial counting \eqref{countpol1} below, of all generic cases, would all be missing.

\begin{figure}[t!]
\centering \includegraphics[width=0.9\textwidth]{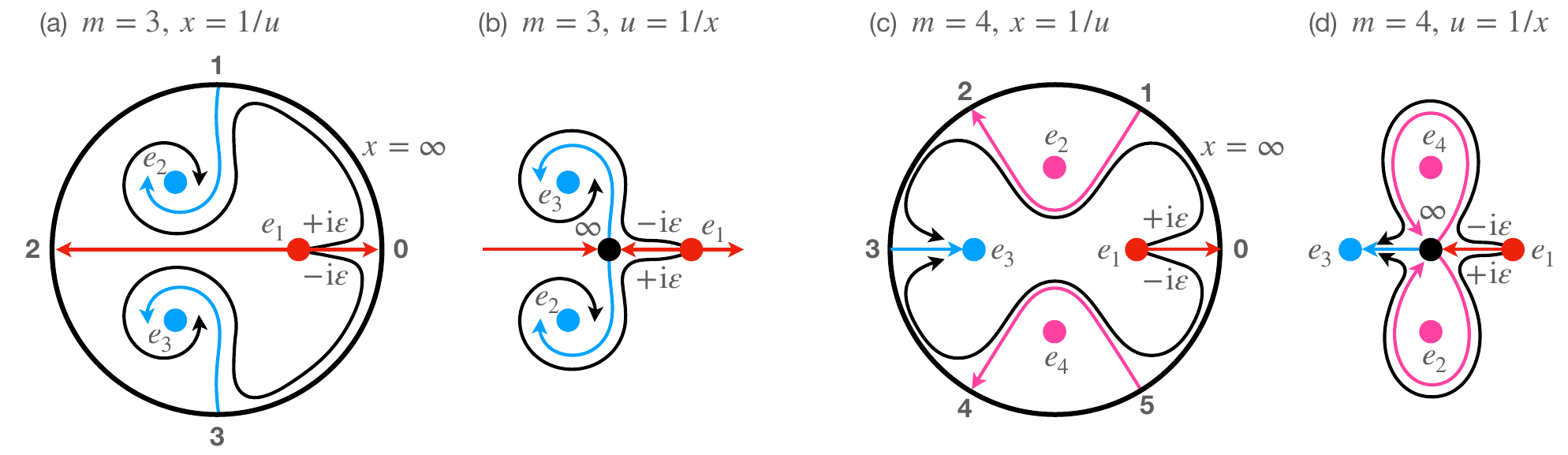}
\caption[Cyclotomic ODEs of degrees 3 and 4]{\emph{
Schematic phase portraits, in real time, of complex-valued ODEs \eqref{ODEP}, for cyclotomic polynomial vector fields $\dot{x}=P(x):=x^m-1$.
The equilibria $x\!=\!e_j$ are the $m$-th roots of unity.
Blue indicates blow-down orbits and their sink targets.
Blow-up orbits, and the sources from which they emanate in reverse time, for $t\searrow -\infty$, are red.
For $m\!=\!3$ see $x$ in (a), and $u\!=\!1/x$ in (b).
Similarly, $x$ in (c) and $u\!=\!1/x$ in (d) refer to $m\!=\!4$.
The invariant \emph{Poincaré circle} $r\!=\!|u|\!=\!0,\ \phi\in\mathbb{S}^1$ of the Poincaré compactification of $x$, by a circle $u=r\exp(\mi\phi)$ ``at infinity'', is marked black in (a), (c).
In (b), (d), a black dot marks the one-point compactification $x\!=\!\infty$ by $u\!=\!1/x\!=\!0$ in the Riemann sphere $\widehat{\C}\cong\S^2$.
Note the real-time blow-up stars of $m-1$ red blow-up orbits towards $u\!=\!0$, which alternate with $m-1$ blue blow-down orbits around $u\!=\!0$, in (b), (d).
In (a), (c) these become separatrices associated to saddles on the circle $u\!=\!0$ ``at infinity''.
Saddles are marked by integers $\mathbf{0},\ldots, \mathbf{2}m\mathbf{-3}$.
Two Lyapunov centers (purple) of purely imaginary $P'(e_j)=\pm4\,\mi$ occur in (c), (d), by nongeneric degeneracy.
They are surrounded by foliated counter-rotating nests of periodic orbits with minimal period $\pi/2$.
Each nest, in turn, is bounded by a homoclinic orbit (also purple) of conflated blow-down-up orbits; see (d).
In (c), these become heteroclinic orbits between boundary saddles on the circle $u\!=\!0$ ``at infinity''.
}\cite{FiedlerYamaguti}}
\label{figm34}
\end{figure}

\begin{figure}[t] 
\centering \includegraphics[width=0.9\textwidth]{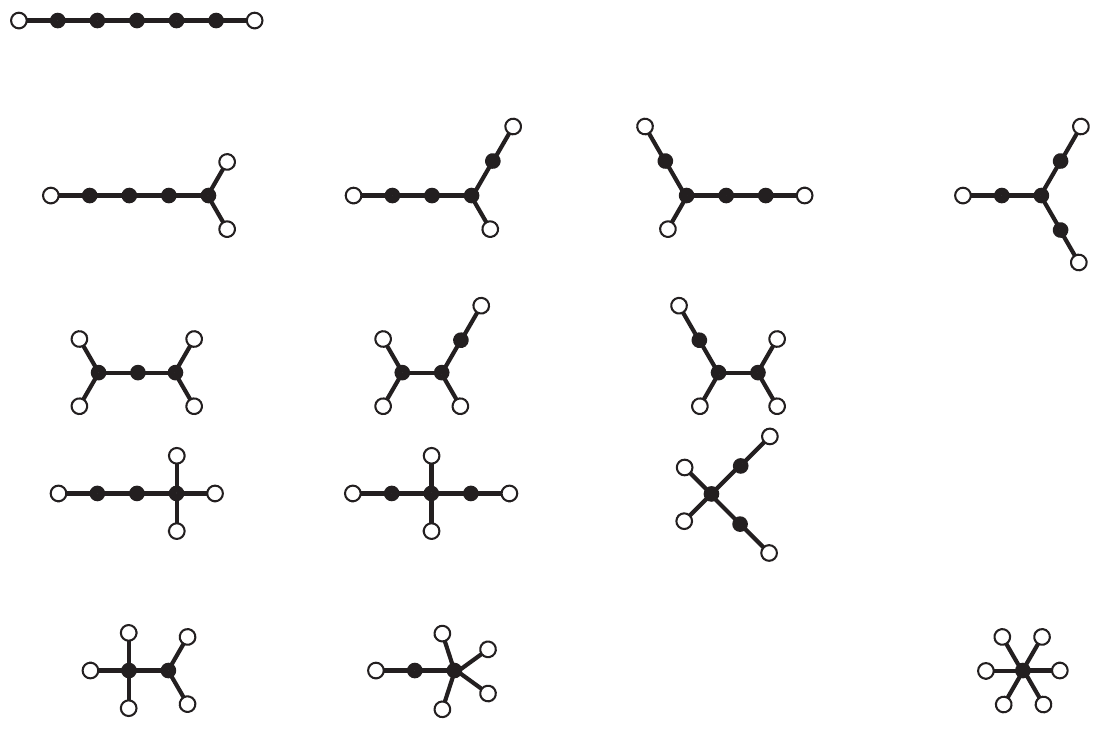}
\caption[Generic polynomials of degree 7, and their 14 trees on the Riemann sphere]{\emph{
The $14$ planar trees $\mathcal{T}$, alias reduced connection graphs $\mathcal{C}^*$, with $m\!=\!7$ vertices, up to orientation preserving equivalence.
Note the two pairs which are mirror-symmetric to each other, but not mirror-symmetric, individually.
Terminal vertices, of edge degree one, are distinguished by circles. 
Solid dots mark all other vertices.
The trees are enumerated by increasing number $2,\ldots,6$ of terminal vertices.
Each tree possesses two bi-colorations which mark the vertices as red sources and blue sinks, alternatingly.
These also defines edge orientations, from sources to sinks.
Color swaps correspond to time reversal.
}\cite{FiedlerYamaguti}}
\label{fig14trees}
\end{figure}

To fix the degree, we assume $p_0\neq0$, viz.~$p_0=1$ without loss of generality.
Analogously to section \ref{Fol}, we compactify blow-up of \eqref{ODEP} by projective coordinates $[\xi\!:\!\zeta]=[x\!:\!1]=[1\!:\!u]\in\CP^1=\widehat{\C}$, i.e.~$u\!=\!1/x$.
As in \eqref{ODEuz} and \eqref{tt1} we obtain
\begin{equation}
\label{ODEu}
\dot u(t_1)= -u\,P_1(u)\,
\end{equation}
with $P_1(u):=u^mP(1/u)=p_0+\ldots+p_mu^m$\ , and $dt=u^{m-1}dt_1$\,.
Note the linearization with eigenvalue $\lambda=-p_0\neq0$ at the equilibrium $u=e_0:=0$, alias $x\!=\!\infty$.
See sections \ref{Over} and \ref{Int} for terminology.

The ODE flow near $x\!=\!\infty$, alias $u\!=\!0$ in blow-up coordinate $u:\!=\!1/x$, illustrates  discrepancy-free blow-up loops $(\gamma^t,\gamma^u)$, in complex time, and their associated blow-up star, in real time.
Indeed, theorems \ref{thmMs}, \ref{thmloop} below, with dummy variables $z\!=\!0$, will identify local minimal blow-up loops $\gamma^t,\gamma^u \subset \Co$, defined by simple closed loops $\fw(\gamma^u)=1$; see definition \ref{defloop}.
Then proposition \ref{proploop} determines the winding number $\fw(\gamma^t)=m-1$ of the associated blow-up star.
In fact, $m-1$ complex blow-up branches (red), for real $0>t\nearrow T=0$, alternate with their continuations by $m-1$ blow-down branches (blue), for real $0<t\searrow T$ after blow-up time $T$.
See definition \ref{defstar} and figure \ref{figm34}, (b), (c) for illustration.
In sections \ref{z=e}, \ref{Ls} below, we discuss the adaptation of these scalar considerations to the dynamics in stable leafs of blow-up; see in particular the schematic illustration in figure \ref{figm234}.

Let us assume hyperbolicity of all finite equilibria $x\!=\!e_j$\,, in real time, i.e.\ $\Re P'(e_j)\neq0$.
In particular, all zeros $e_j$ are simple.
In real time $t$, the equilibrium $x\!=\!e_j$ is a \emph{sink}, for $\Re P'(e_j)<0$, and a \emph{source} if $\Re P'(e_j)>0$.

Under further genericity assumptions on the coefficients $p_j$ of complex polynomials $P(x)$ of degree $m$, we have studied the resulting global flows of ODE \eqref{ODEP} on the Riemann sphere $\widehat{\C}\cong\S^2$.
In real time $t$, we have obtained $1\,{:}\,1\ C^0$ orbit equivalence to generic real-time Morse flows on the real 2-sphere $\S^2$, with a total of $m$ nondegenerate saddles and sinks $e_j$ and a single pole of order $m-2$ at infinity.
See definition \ref{defequi} for equivalence terminology.
See figure \ref{figm34} for cyclotomic examples of degrees $m=3,4$, where sources and blow-up orbits are colored red.
Blue color indicates sinks and blow-down orbits towards them.

$C^0$ orbit equivalence classes of real-time flows can be represented by (reduced) \emph{connection graphs} $\cC^*$ between source and sink equilibria. 
Vertices are the equilibria $e_j$\,.
Single directed edges represent the existence of heteroclinic orbits $x(t)$ from sources, for $t\searrow -\infty$, to sinks for $t\nearrow\infty$; compare \eqref{het}.
Omitted blow-up orbits to the pole at infinity, from source equilibria, and blow-down orbits, towards sink equilibria, separate basins of attraction in forward and backward time, respectively.
The resulting connection graphs $\cC^*$ are planar trees $\cT$ of $m$ vertices -- \emph{all of them}.
See figure \ref{fig14trees} for an illustration of degree $m\!=\!7$.

The $1\,{:}\,1$ representation of $C^0$ orbit equivalence classes by planar trees enables combinatorial counting of the possible configurations.
The precise classification is up to planar homeomorphisms which preserve planar orientation.
They also preserve orbits and edges, respectively, but are allowed to reverse all edge directions simultaneously.
Equivalence to dual circle diagrams of $m-1$ unlabeled disjoint chords (handshakes) among $2(m-1)$ equidistantly spaced vertices on the circle allows the following counts $A_m$ of equivalence classes, due to \cite{countpol1, countpol2}:
\begin{equation}
\label{countpol1}
\begin{aligned}
A_m\,=\,\tfrac{1}{2(m-1)m}\, \tbinom{2(m-1)}{m-1}\,&+\,\tfrac{1}{4(m-1)}\tbinom{m}{m/2}\,+\,\tfrac{1}{m-1}\,\boldsymbol{\phi}(m-1) \,+\\
                  &+\,\tfrac{1}{2(m-1)}\,\sum_{k=2}^{m-2}\,\tbinom{2k}{k}\, \boldsymbol{\phi}(\tfrac{m-1}{k}) \,. 
\end{aligned}
\end{equation}
See also \cite{oeispol}.
Here $\boldsymbol{\phi}$ denotes the Euler totient count of coprime integers, and the sum only runs over proper divisors $k$ of $m-1$.
For odd $m$, the second summand on the right is omitted.
Explicit counts $A_m$ for $2\leq m\leq 16$ are
\begin{equation}
\label{countpol}
1, 1, 2, 3, 6, 14, 34, 95, 280, 854, 2694, 8714, 28640, 95640, 323396.
\end{equation}

\section{Example: Linear foliations and linear holonomy in $\C^2$}\label{LinC2}

Our main results on blow-up loops $(\gamma^t,\gamma^{uz})$ and their associated blow-up stars  are based on linearization at blow-up equilibria, in projective compactifications.
See definitions \ref{defloop} and \ref{defstar}, propositions \ref{propuvw}, \ref{proploop}, and sections \ref{z=e}-\ref{RatFol}, \ref{B} below.
As our first example of a system in complex time, we therefore study the complex foliations of the linear ODE
\begin{equation}
\label{xylin}
\begin{aligned}
    \dot x(t) &=\alpha x+\beta y\,,   \\
    \dot y(t) &=\gamma x+\delta y\,,
\end{aligned}
\end{equation}
on $\C^2$, with constant complex coefficients. 
ODE \eqref{xydiag} diagonalizes \eqref{xylin}.
By linearity, ODE \eqref{xydiag} \emph{does} define a flow on complex projective space $[x\!:\!y]\in\CP^1=\widehat{\C}$, which we explore in section \ref{LinProj}.
Following \cite{Ilya}, alternatively, section \ref{LinHol} studies periodicity versus quasiperiodicity of the associated holonomies at $x\!=\!y\!=\!0$.
The remaining sections \ref{z=e}-\ref{RatFol} assume diagonal ``linearity at infinity'', i.e.~of ODE \eqref{ODEuz} in projective coordinates $u\!=\!1/x,\ z\!=\!y/x$ and in rescaled time $t_1$; see \eqref{uzlin}.
Only Euler multipliers $\rho=u^{m-1}$, however, allow us to pass to original time $dt=\rho dt_1$\,; see \eqref{tt1} and the blow-up stars which arise in figure \ref{figm234}.
The trivial stable leaf $z\!=\!0$ is addressed in section \ref{z=e}.
The general leafs depend on irrational versus rational spectral quotients, as described in sections \ref{IrrFol} and \ref{RatFol}, respectively.

A Pfaffian form $\omega$ of the foliation $\omega\!=\!0$ associated to \eqref{xylin} is
\begin{equation}
\label{omegalin}
\omega=-(\gamma x+\delta y)\,dx+(\alpha x+\beta y)\,dy\,;
\end{equation}
compare \eqref{omega}-\eqref{ODEfg}.
Assuming nonvanishing determinant and discriminant of the coefficient matrix, a linear substitution diagonalizes \eqref{xylin}:
\begin{equation}
\label{xydiag}
\begin{aligned}
    \dot x(t) &= \lambda_1\, x\,,   \\
    \dot y(t) &= \lambda_2\, y\,,
\end{aligned}
\end{equation}
with complex eigenvalues $\lambda_1,\,\lambda_2$ and associated Pfaffian $\omega=-\lambda_2dx+\lambda_1dy$.
Finite-time blow-up does not occur in the linear case, of course.
Still, we discuss the foliation $\omega\!=\!0$ in some detail, as an elementary example of independent interest and for later perusal in projective compactifications \eqref{ODEuz}, which \emph{do} involve blow-up in finite time.
In particular, we illustrate the concept of holonomy maps in the linear case.

\subsection{Projective dimension reduction}\label{LinProj}

Equivariance of ODE \eqref{xylin} under self-similar scaling reduces dimension.
For any constant scaling factor $\sigma\in\Co$, \emph{scaling equivariance} observes that the scaled orbit $(\sigma x(t),\sigma y(t))$ is a solution whenever $(x(t),y(t))$ is.
In diagonalized variables \eqref{xydiag}, then, $z:\!=\!y/x$ is scaling invariant and represents orbits under the scaling group of $\sigma$.
Since nonvanishing discriminant implies $\lambda_1\neq\lambda_2$\,, the resulting flow on the Riemann sphere $z\in\widehat{\C}=\CP^1$ for the linear ODE
\begin{equation}
\label{linz}
\dot z=(-\dot x/x+\dot y/y)z=(-\lambda_1+\lambda_2)z
\end{equation}
defines a single cylindrical leaf $\fL\,{:}\ z\in\Co$.
By equivariance, \eqref{linz} is autonomous in the scaling invariant $z$. 
For dynamics see \eqref{Riclin} and, in other coordinates, the Riccati equation of figure \ref{figRiccati}.
The foliation $\omega\!=\!0$ for $(x,y)\in\CP^2$ is then parametrized by the scaling parameter $\sigma$, with holomorphically equivalent flows \eqref{linz} on most leafs.
The only exceptions are the two trivial leafs $\fL_0,\,\fL_\infty$ at $z=0,\infty$ which are invariant under scaling by $\sigma$.
These are readily identified as the nonequilibrium parts of the invariant ``manifolds'' given by the coordinate axes, alias eigenspaces $\{y\!=\!0\}$ and $\{x\!=\!0\}$ of the nonzero eigenvalues $\lambda_1$ and $\lambda_2$\,, respectively.
The trivial leafs are cylindrical, as well, because the equilibria $x\!=\!y\!=\!0$ and infinity in $\widehat{\C}=\CP^1$  have to be omitted, in either case.
All cylindrical leafs are topologically and biholomorphically equivalent.
Since leafs consist of single orbits, in complex time, they are trivially orbit equivalent.
With few trivial exceptions, the three types of leafs fail to be analytically flow equivalent, due to their complex time periods $2\pi\mi/\lambda_\iota$ in $x,y$.
Even if we allow for rescaling of time by  nonzero real factors, only real collinear nonzero eigenvalues $\lambda_1,\,\lambda_2$\,, and $\lambda_3:=\lambda_2-\lambda_1$ qualify for leaf equivalence.

\subsection{Linear holonomy}\label{LinHol}

Alternatively to scaling invariance, we may follow \cite{Ilya} and describe the foliation $\omega\!=\!0$ associated to the diagonal linear system \eqref{xydiag} by a \emph{holonomy map} $h$ of (vertical) $x$ over (horizontal) $y$, as follows.
Fix $y=y_0\neq0$ and let $x_0\in\C$ label the leafs of the local foliation over $y\!=\!y_0$\,.
We call $\{y\!=\!y_0\}$ a \emph{transverse section} to the foliation $\omega\!=\!0$.
Indeed, the leaf $\fL$ of $(x_0,y_0)$ is transverse to the section $\{y\!=\!y_0\}$, at $(x_0,y_0)$, because $\dot y=\lambda_2y_0\neq0$ there.
Now let $s\mapsto y(s)=\gamma^y(s),\ 0\leq s\leq2\pi$ parametrize any closed differentiable Jordan curve $\gamma^y\subset\Co$, oriented positively around $y\!=\!0$ and starting at $\gamma(0)\!=\!y_0$\,, say.
By complex time-invariance of leafs, \eqref{xydiag} implies that the leaf $\fL$ of $(x_0,y_0)$ intersects the transverse section $\{y=y(s)\}$ at $(x(s),y(s))$, where 
\begin{equation}
\label{lingam}
x'= \frac{dx}{dy}\cdot \frac{dy}{ds}=\frac{\dot x}{\dot y}\cdot y'(s)= \frac{\lambda_1}{\lambda_2} \cdot \frac{y'(s)}{y(s)}\cdot x \,.
\end{equation}
Here $':=d/ds$, and $y'(s)/y(s)$ is given by the choice of $y(s)=\gamma^y(s)$.
Integrating $x'/x$ over the cycle $\gamma^y$, this implies
\begin{equation}
\label{xhol}
x(2\pi)=x_0\exp\left(\int_0^{2\pi} \frac{x'(s)}{x(s)}\,ds\right) =x_0\exp\left(\int_{\gamma^y}\, \frac{\lambda_1}{\lambda_2} \cdot \frac{dy}{y}\right)= x_0\exp\left(2\pi\mi\,\lambda_1/\lambda_2\right) \,, 
\end{equation}
by the residue theorem along $\gamma^y$.
This holds independently of our particular choice of $\gamma^y$ above.
The holonomy map is then defined as the return map
\begin{equation}
\label{linhol}
h\,{:}\quad x_0\mapsto x(2\pi)=x_0\exp\left(2\pi\mi\,\lambda_1/\lambda_2\right)\,,
\end{equation}
characterized by the \emph{spectral quotient} $\lambda:=\lambda_1/\lambda_2\ \mathrm{mod}\,\Z$\,.
Analogously to \eqref{discrepancy}, \eqref{discrepancyuz}, we call $x(2\pi)-x_0$ the \emph{discrepancy}, and $\exp\left(2\pi\mi\,\lambda_1/\lambda_2\right)$ the \emph{holonomy multiplier}, of the linear holonomy $h$.
Alternatively, we may simply consider \eqref{xydiag} in imaginary time direction $t=\mi s/\lambda_2$\,. 
This amounts to the particular choice $y(s)=y_0\exp(\mi s)$, for the circle $\gamma^y$.
The holonomy $h$ then follows by integration of $x'(s)=\mi\lambda\, x(s)$ along $0\leq s\leq2\pi$.
Since all nonzero orbits of $y'(s)=\mi\, y(s)$ are then periodic of the same period $2\pi$, we may also identify $h$ as a (remaining) Poincaré return map along the periodic orbit $\gamma^y$ in real time $s$.
In this case, the holonomy $h$ is often called the \emph{monodromy} of the periodic orbit $\gamma^y$, and the holonomy multiplier is a Floquet multiplier of $\gamma^y$ \cite{ArnoldEnc, Ilya}.

In coordinate invariant language, the holonomy map of $\gamma^y$ encodes the return map defined by the leafs of the foliation $\omega\!=\!0$, extended over $\gamma^y$.
We repeat: holonomy is defined by the foliation, and \emph{not} by the particular flows on its leafs.
This eliminates the role of Euler multipliers $\rho$\,.
Up to biholomorphic conjugacy, holonomy is also independent of the chosen cross section to the foliation and only depends on the free homotopy class of $\gamma^y\subset\Co$.
The holonomies of cycles $\gamma^y$ with winding number $n:=\fw(\gamma^y)\in\Z\cong\pi_1(\Co )$ around $y\!=\!0$ correspond to iterates $h^n$.
Note that trivial holonomy $h^n=\mathrm{Id}$ occurs if, and only if, $\lambda=\lambda_1/\lambda_2\in\Q\setminus\{0\}$ is rational with $n\lambda\in\Z$.
Irrational spectral quotients $\lambda$ lead to \emph{quasiperiodic holonomy iterates}, even in the case of \emph{real} $\lambda_1,\lambda_2$\,.

Depending on our precise set-up, however, the same linear system \eqref{xydiag} may also spawn other holonomy maps.
For example, the holonomy of (vertical) $y$ over positively oriented (horizontal) cycles $x(s)=\gamma^x(s)\in\Co$ around $x\!=\!0$ in the complex $x$-plane becomes
\begin{equation}
\label{linholy}
h\,{:}\quad y_0\mapsto y_1=\exp\left(2\pi\mi\,\lambda_2/\lambda_1\right) y_0\,.
\end{equation}
Note the inverse spectral quotient $\lambda_2/\lambda_1=1/\lambda$ in the holonomy multiplier, compared to $\lambda=\lambda_1/\lambda_2$ in \eqref{linhol}.

Blow-up coordinates $u:=1/x,\ z:=y/x$ provide the diagonal linear variant
\begin{equation}
\label{uzdiaglin}
\begin{aligned}
    \dot u(t_1)&= -\lambda_1\, u\,,   \\
    \dot z(t_1) &= (\lambda_2-\lambda_1)\, z
\end{aligned}
\end{equation}
of \eqref{uzdiag}; see \eqref{xyz}, \eqref{ODEuz}, and \eqref{f12}.
The ODE for $z\!=\!y/x$ has already appeared in our discussion \eqref{linz} of scaling.
We obtain spectral quotients $-\lambda_1/(\lambda_2-\lambda_1)=1/(1-1/\lambda)$ for the $u$-holonomy over positive cycles $\gamma^z$ in $z$, and the inverse $-1/\lambda\ \mathrm{mod}\,\Z$ for the $z$-holonomy over $\gamma^u$.
Passage from $x$ to $u\!=\!1/x$ reverses the sign for $\dot{u}$ in \eqref{uzdiaglin}, and the sign of the $y$- or $z$-holonomy over $x$.
The reason, of course, is that closed curves $\gamma^x$ around $x\!=\!0$ in the cylinder $\widehat{\C}\setminus\{0,\infty\}$ reverse orientation when considered as closed curves $\gamma^u$ around $u\!=\!0$ alias $x\!=\!\infty$ on the Riemann sphere.
Similarly, an orientation reversal of $\gamma$ leads to sign reversal of the spectral quotients which characterize holonomy in the various set-ups.

The alternative blow-up coordinates $v:=1/y,\ w:=x/y$ provide the diagonal linear variant
\begin{equation}
\label{vwdiag}
\begin{aligned}
    \dot v(t_2) &= -\lambda_2\, v\,,   \\
    \dot w(t_2) &= (\lambda_1-\lambda_2)\, w\,;
\end{aligned}
\end{equation}
see \eqref{xyz}, \eqref{ODEvw}, and \eqref{f12}.
For $v$-holonomy over $w$ and $w$-holonomy over $v$ we obtain mutually inverse spectral quotients $1/(1-\lambda)$ and $-\lambda\ \mathrm{mod}\,\Z$, respectively.
In total, this associates $\pm\lambda\,,\pm 1/(1-1/\lambda)\ \mathrm{mod}\,\Z$, and their inverses, to spectral quotients $\lambda$ of holonomies around the three equilibria $[0\!:\!0\!:\!1],\ [1\!:\!0\!:\!0],\ [0\!:\!1\!:\!0]$ of the projective compactification of the linear diagonal ODE \eqref{xydiag}.

The fixed point indices +1 of each equilibrium, by the way, add up to the Lefschetz number 
$\cal{L}$$(\Phi^t)=3$ of the associated flow $\Phi_t$\,, for fixed small $|t|$ which prevent nonstationary periodicity.
By homotopy to the identity map $\Phi^0=\mathrm{Id}$, this number justly coincides with the genus of $\CP^2$, i.e.\ with the alternating sum of its Betti numbers.

\subsection{The stable leaf $\{z=0\}$}\label{z=e}

To prepare for blow-up section \ref{B}, in the language set by proposition \ref{proploop}, we now assume total linearity of ODE \eqref{ODEuz}.
Linearity requires $f,g$ with constant $f_1$ and linear $g_1$ in \eqref{f12}, at any given degree $m$.
For $f(x,y):=-\lambda_1x^m$ and $g(x,y):=(\lambda_2-\lambda_1)x^{m-1}y$, we obtain the linear diagonal ODE
\begin{equation}
\label{uzlin}
\begin{aligned}
     \dot u(t_1)&=\lambda_1\,u\,,   \\
    \dot z(t_1)&=\lambda_2 \,z\,.
\end{aligned}
\end{equation}
In other words, we have just substituted $u,z$ for $x,y$ in \eqref{xydiag}.
We now assume $\lambda_1,\, \lambda_2$ are real with $\lambda_1<0$.

The simplest leafs are contained in the punctured invariant coordinate axes $\{u\!=\!0\}$ and $\{z\!=\!0\}$, each a cylinder $\Co$.
The invariant leaf of $\{u\!=\!0\}\subset\CP^2$ is an artifact of projective compactification at $x\!=\!1/u\!=\!\infty$.
We therefore discuss the trivial, but more interesting, \emph{stable leaf} $\fL$ of $\{z\!=\!e\}$, here for $z\!=\!e\!=\!0$.

\begin{figure}[t]
\centering \includegraphics[width=0.9\textwidth]{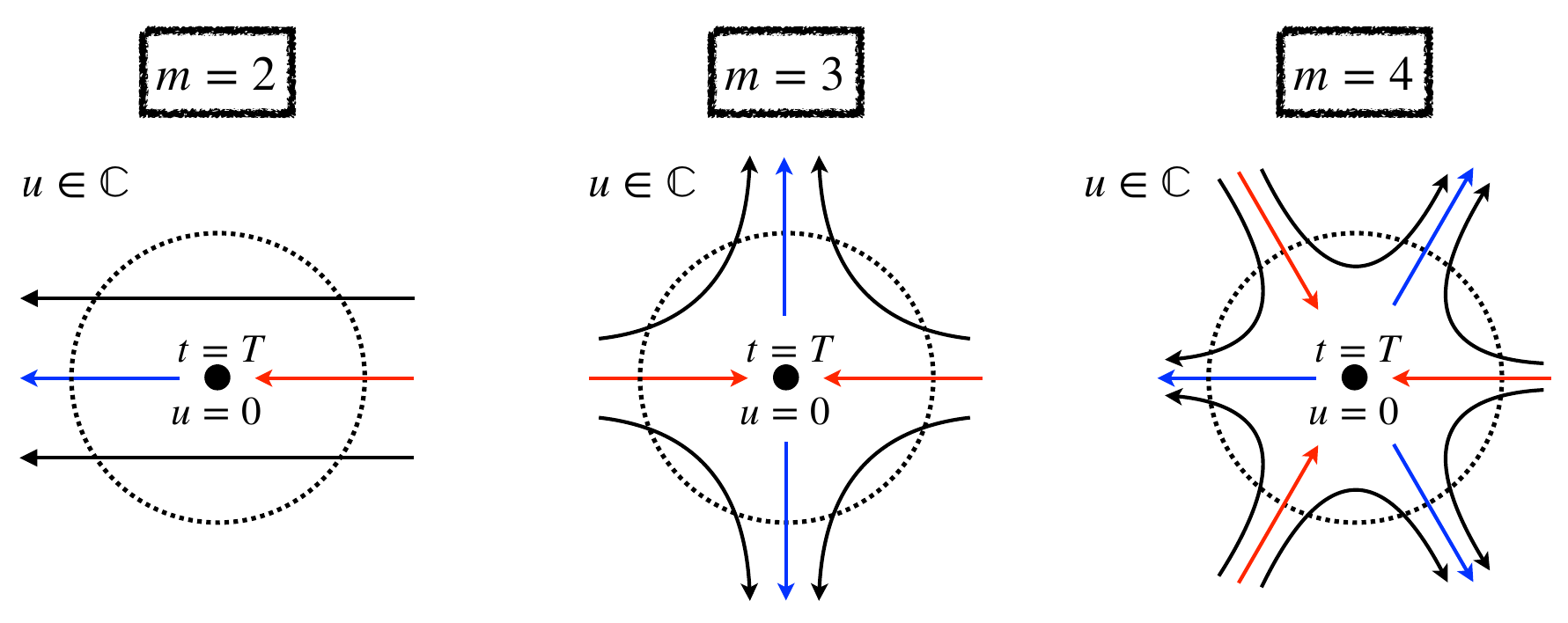}
\caption[Blow-up loops and real-time blow-up stars]{\emph{
Phase portraits of minimal blow-up loops and real-time blow-up stars near blow-up equilibrium $u\!=\!0$, in the complex plane $u\in\C$ of the stable leaf $\fL\subset\{z\!=\!0\}$ of $\lambda_1<0$ in the linear diagonal ODE \eqref{uzlin}. 
Time rescaling to original time $t\sim -u^{m-1}$ is by the Euler multiplier $\rho=u^{m-1}$ of \eqref{tt1}.
See \eqref{tu} for $m=2,3,4$.
The real-time blow-up stars consist of red blow-up and blue blow-down separatrix orbits.
Other orbits (black) follow parallel real time directions of original time $t\in\C$, with fixed nonzero imaginary parts $\Im t $.
Blow-up of real $0<x=1/u\nearrow+\infty$ occurs in finite real time $0>t\nearrow T=0$.
The blow-up loop $\gamma^{uz}=\gamma^u\times\{0\}$ (black dashed) arises from $m-1$ iterates $\gamma^t$ of a simple loop in $t$, and then without any resulting discrepancy. 
This illustrates the (minimal) winding numbers $\fw_u=1$ and $\fw_t=(m-1)\fw_u=m-1$ of $\gamma^u$ and $\gamma^t$ asserted in \eqref{wtu} of proposition \ref{proploop}.
The alternating $m-1$ red and $m-1$ blue separatrices of the hyperbolic sectors mark real-time blow-up stars.\\
Note the parallel flow, for the quadratic case $m\!=\!2$ (left).
Real-time blow-down $-\infty\swarrow x=1/u<0$ occurs from continuing $0=T\swarrow t>0$.
The continuation can be reached, discrepancy-free, by single complex Masuda detours in time along the upper or the lower half of the blow-up loop $\gamma^t$.
Also compare the cyan orbits in the Riccati case of figure \ref{figRiccati}, for a more global view.\\
The cubic case $m\!=\!3$ (center) requires two cycles through the complex loop $\gamma^t$.
Accordingly, we traverse the negative axis of original time $t<0$ twice, resulting in two red blow-up orbits along real $u\neq0$.
Indeed, cubic $\dot x=x^3$ blows up for positive \emph{and} negative $x\!=\!1/u$.
Any blue blow-down continuation to $0=T\swarrow t>0$, along the blow-up loop $\gamma^u$, requires purely imaginary $x,u$.
Moreover, there are two conflicting options for such continuations beyond blow-up.\\
For $m\!=\!4$ (right), as for all even $m$, the real time axis again allows for a continuation beyond blow-up of real $x=1/u>0$, by blow-down along real $u<0$.
However, triple passage through $t<0$ and $t>0$, each, along the triply iterated blow-up loop cycle $\gamma^t$ now leads to three blow-up separatrices, which alternate with three blow-down separatrices.
For further examples see figure \ref{figm34} (b), (d).
}}
\label{figm234}
\end{figure}

The stable leaf $\{z\!=\!0\}$ of $\dot u(t_1)=\lambda_1u(t_1)$ for the blow-up with $\lambda_1<0$ corresponds exactly to the linear variant \eqref{Riclin} of Riccati flows in $u$, with $\lambda_1:=a(e_1-e_2)$.
The flow looks regular, in the time derivative $\dot{}=d/dt_1$\,.
Figure \ref{figm234} of the cases $m=2,3,4$, however, shows that some caution is required in the real-time interpretation of this trivial result within the stable leaf $\fL$.
The intricacies come from the nonlinear effect of the Euler multiplier $\rho=u^{m-1}$ in corollary \ref{tt}. 
In the blow-up context of section \ref{LoopStar}, the Euler multiplier $\rho$ is required when passing from $t_1$ to original time $dt=u^{m-1}dt_1$\,; see \eqref{tt1}.
The relation \eqref{tuexpan} between $t$ and $u$ becomes exact, here, without any higher order terms:
\begin{equation}
\label{tu}
t-T=\tfrac{1}{\lambda_1(m-1)}\,u^{m-1}\,.
\end{equation} 
Therefore the orbit $\gamma^{uz}=\gamma^u\times\{0\}$ of the solution $(u(t),0)$ closes up, first and as a discrepancy-free \emph{minimal blow-up loop}, upon $m-1$ iterates $\gamma^t$ of a simple loop in $t$.
The winding numbers are $\fw_u=1$ and $\fw_t=m-1$; see \eqref{wtu}.

Blow-up occurs for real $0>t\nearrow T=0$ and real $0<x(t)=1/u(t_1)\nearrow+\infty$.
For even $m$, the blow-up connects to blow-down for real $0=T\swarrow t>0$ and real $-\infty\swarrow x=1/u<0$.
See the cases $m=2,4$ in figure \ref{figm234}.
The separatrix orbits of the associated real-time blow-up star are marked red and blue there.
After $\pm(m-1)/2$ cycles through a simple $t$-loop around $T\!=\!0$, the Masuda detour traverses half a $u$-loop $\pm \Im \gamma^u\geq0$ of the blow-up loop $(\gamma^t,\,\gamma^u)$, black dashed in figure \ref{figm234}, and therefore leads to such naive real continuation, without any discrepancy. 

For odd $m\geq3$, however, $\pm(m-1)/2$ cycles through $\gamma^t$ lead back to $t>0$, i.e.\ to blow-up again, in the blow-up loop.
See the case $m\!=\!3$ of figure \ref{figm234}, center.
Then there is no canonical choice among the $m-1$ mutually inconsistent blue blow-down orbits, all complex-valued, to continue real red blow-up, or any other of the $m-1$ red blow-up orbits.
Still, each available option of $\gamma^{uz}$ is reachable, consistently, through the $(m{-}1)$-fold blow-up loop along $\gamma^t$.
Small purely imaginary perturbations $t\pm\eps\mi$ of red real blow-up time $t$, which had been attempted so far \cite{Masuda1, Masuda2, Stukediss, Stukearxiv}, are tied to $\gamma^t$ itself, and therefore lead to completely inconsistent ``neighboring'' blue continuations for $m\geq3$, in spite of non-discrepancy over a full blow-up loop $\gamma^u$.
See the black orbits in figures \ref{figm234} and \ref{figm34} (b). 
Indeed, the above considerations arose in section \ref{Scalar} on the global dynamics of generic polynomial vector fields $\dot x=P(x)$ of degree $m$ on the Riemann sphere $\widehat{\C}=\CP^1$.
In particular, the above degenerate real-time phase portraits of $u$ at $u\!=\!1/x\!=\!0$ can also be described in terms of a degenerate blow-up equilibrium.
Locally, that degeneracy generates $2(m-1)$ hyperbolic sectors, in the sense of chapter VII.9 in \cite{Hartman}.
The standard case of a linearly hyperbolic equilibrium at infinity corresponds to $m\!=\!3$.
See also \cite{FiedlerShilnikov} for further discussion in terms of the Poincaré compactification of $\C$, rather than the projective one-point compactification by the Riemann sphere $\CP^1=\widehat{\C}$.

\subsection{Irrational spectral quotients}\label{IrrFol}

Consider the linear diagonal ODE \eqref{uzlin} with real irrational spectral quotients $0\neq\lambda=\lambda_1/\lambda_2\in\R\setminus\Q$ next.
As in \eqref{linholy},  we obtain the $z$-holonomy
\begin{equation}
\label{irrholy}
h\,{:}\quad z_1=\exp(2\pi\mi/\lambda)\, z_0
\end{equation}
at irrational angles $\lambda$, over any fixed $u_0\in\Co$.
In particular, the holonomy of the foliation is \emph{quasiperiodic} in $u$, over any $z_0\neq0$.
Therefore, the associated leafs 
\begin{equation}
\label{irrleaf}
\fL\, {:}\quad z=\exp(\log(u/u_0)/\lambda)\,z_0
\end{equation}
are infinitely branched over $u\in\Co$ by the Riemann surface of the logarithm.
They never quite close up to form a discrepancy-free blow-up loop, in the sense of definition \ref{defloop}.
Note density of each leaf $\fL$ in $\CP^2$.
So, solution trajectories $\gamma^{uz}$ in the sense of proposition \ref{proploop} close \emph{almost}: not quite, but with any required precision, under sufficiently high iterations of any simple circular Masuda loop $\gamma^t$.

In the saddle case $\lambda_1<0<\lambda_2$ of negative spectral quotients $\lambda<0$, the irrational leafs $\fL$ do not contain real blow-up solutions asymptotic to $(u,z)=(0,0)$.
The node case $\lambda>0$, in contrast, does contain real blow-up solutions, in each such leaf.
Then quasiperiodicity prevents discrepancy-free blow-up loops around blow-up, but also provides arbitrarily close approximations to holonomy multipliers 1, i.e. to their ever elusive non-discrepancy.

\subsection{Rational spectral quotients, and Riemann branching}\label{RatFol}

Consider the linear diagonal ODE \eqref{uzlin} with nonzero rational spectral quotients $0\neq\lambda=\lambda_1/\lambda_2=n_1/n_2\in\Q$ next, i.e.\ with coprime integers $n_1,n_2$\,.
Rescaling time by a constant real Euler multiplier, we may assume $\lambda_\iota=-n_\iota$ for both eigenvalues.
Here $n_1>0$, and $n_2\neq0$ possesses any sign.
By \eqref{irrholy}, the holonomy multiplier $\exp(2\pi\mi/\lambda)$ and the linear $z$-holonomy $h$ of nontrivial leafs $\fL$ become primitive $n_1$-th roots of unity:
\begin{equation}
\label{hroot}
h^{n_1}=\mathrm{Id}\,.
\end{equation}
All $z\in\Co$ are periodic under iteration of $h$, with minimal period $n_1$\,.

Holonomy traces leafs.
In time $t_1$ of \eqref{uzlin} and with $\theta:=\exp(-t_1)$, we have explicit exponential solutions
\begin{equation}
\label{uzsol}
\begin{aligned}
     u(t_1)&=u_0\exp(-n_1t_1)=u_0\,\theta^{n_1}\,,   \\
     z(t_1)&=z_0\exp(-n_2t_1)=z_0\,\theta^{n_2}\,.
     \end{aligned}
\end{equation}
In other words, any nontrivial leaf is algebraic:
\begin{equation}
\label{Ln1n2}
\fL^c\,{:}\quad u^{n_2}=c\,z^{n_1}\neq0\,
\end{equation}
with constant $c=z_0^{n_1}u_0^{-n_2}\in\Co$\,.
Note flow equivalence to the standard leaf $c\!=\!1$, by biholomorphic diagonal scaling of $u,z$.

Only for $n_2>0$, the leaf contains real-time orbits towards blow-up at $u\!=\!z\!=\!0$. 
Parametrization by $\theta\in\C$ then identifies the closure $\fL^c\cup\{0\}$ of any nontrivial leaf as a \emph{Riemann surface} $\cR$:
the generalized parabola with \emph{ramification index} $n_1>0$ over $u$ and ramification index $n_2>0$ over $z$.
Purely imaginary $t_1\!=\!-\mi s$ defines a blow-up loop $(\gamma^t,\,\gamma^{uz})$, which only closes after a minimal period of $2\pi$ in $s$, without discrepancy.
In particular, the blow-up loop is minimal, according to definition \ref{defloop}.
With resulting winding number $\fw(\gamma^u)=n_1$\,, proposition \ref{proploop} tells us that the blow-up loop requires $\fw(\gamma^t)=(m-1)n_1$ windings around blow-up time $t\!=\!T\!=\!0$ in original time $t$.

We determine the associated blow-up star next; see definition \ref{defstar}.
We pass to original time $t$ via $dt=\rho\,dt_1=u^{m-1}\,dt_1$\,, along the loop $u=u_0\,\theta^{n_1}$ with $t_1=-\mi s,\ \theta=\exp(\mi s)$.
With $\theta\!=\!0$ at blow-up $t\!=\!T\!=\!0$, we obtain
\begin{equation}
\label{tthetalin}
t-T=-\tfrac{1}{(m-1)n_1}u_0^{m-1}\,\theta^{(m-1)n_1}\,.
\end{equation}
This determines the blow-up star on $\fL^c$ to consist of $\fw_t=(m-1)n_1$ radial blow-up and blow-down lines, in the regular parametrization of the Riemann surface $\cR=\fL^c\cup\{(0,0)\}$ by $|\theta|\leq1$.

The winding of the above blow-up loops $\gamma^{uz}$ in the standard leaf $\fL^1$ of $c\!=\!1$ is best visualized within real standard spheres $\S^3=\{|u|^2+|z|^2=\eps^2\}$.
Then $|\theta|=1$ on the unit circle parametrize $\fL^1\cap\S^3$ as an $n_1\!:\!n_2$ torus knot in $\S^3$, unless at least one of $n_\iota\!=\!\pm1$; see for example \cite{Bries, Lamotke}.

Similarly knotted leafs, by the way, have been used by Art Winfree to construct initial conditions for knotted scroll wave filaments of spatially three-dimensional excitable media.
For their unknotting dynamics see \cite{FiedlerMantel, FiedlerScheel}, and the references there.
See \cite{KupitzHauser}, for a delicate experimental setup, and section 11 of \cite{Kupitz} for an experimental confirmation of our predictions on collisions of scroll waves, as a preliminary stage.

\begin{rem}\label{remnonss}
The special case $\lambda\!=\!1$ of vanishing discriminant, $\lambda_1=\lambda_2$\,, leads to algebraic multiplicity two and is not diagonalizable, in general.
Only if we \emph{assume semisimple eigenvalues}, alias geometric multiplicity two, viz.~diagonal linearization \eqref{xydiag}, our previous analysis applies.
This leads to trivial identity holonomies of $z$ over $u$ and of $w$ over $v$ by \eqref{uzdiaglin} and \eqref{vwdiag}, respectively.
\emph{Non-semisimple spectral quotients} $\lambda\!=\!1$, in contrast, fail to possess nontrivial blow-up loops, even in the linear case
\begin{equation}
\label{Jordan}
\begin{aligned}
    \dot u(t_1)&=-u\,,   \\
    \dot z(t_1)&=u-z\,.
\end{aligned}
\end{equation}
Indeed, any nontrivial blow-up loop $\gamma^u$ of $u(t_1)=u_0\exp (-t_1)\neq0$ requires some imaginary period $2\mi \pi k$ of $t_1$\,.
This prevents any periodicity of $z(t_1)=(z_0 + u_0\,t_1)\exp(-t_1)$.
Therefore we had to assume semisimple eigenvalues, in case $\lambda\!=\!1$.
\end{rem}

\begin{rem}\label{remirr}
In all linear holonomy constructions, we emphasize the role of rational versus real irrational spectral quotients $\lambda=\lambda_1/\lambda_2$\,.
For \emph{rational quotients},  iterated loops $\gamma^t$ of $t$ cause solutions $(u(t),z(t))$ to match up consistently to form discrepancy-free \emph{blow-up loops} $(\gamma^t,\gamma^{uz})$, i.e.\ they possess holonomy multipliers 1, after finitely iterated detour cycles.
This attests to leafs $\fL^c$ of the complex foliation as (finitely) \emph{branched Riemann surfaces} $\cR$ at $u\!=\!z\!=\!0$.
See also theorem \ref{thmloop} below.
For \emph{irrational quotients}, in contrast, the leafs $\fL^c$ never quite match up.
Instead, however, iterated Masuda detours $\gamma^t$ will generate \emph{quasiperiodic holonomy multipliers}.
For suitably chosen large iterates, these get arbitrarily close to the multiplier 1 of discrepancy-free matching.
\end{rem}

\section{Example: Reciprocally linear blow-up in $\C^2$}\label{RecC2}

Our third example in complex time are reciprocally linear ODEs of the form
\begin{equation}
\label{xyr}
\begin{aligned}
    \dot x(t) &=1/(\gamma x+\delta y)\,, \\
    \dot y(t) &=1/(\alpha x+\beta y)   
\end{aligned}
\end{equation}
on $\C^2$, with real coefficients. 
In this example, blow-up occurs towards the singularity $x\!=\!y\!=\!0$ where the vector field blows up, just as it did for $|(x,y)|\rightarrow\infty$ before.
In \eqref{Eulerr}, a suitable complex Euler multiplier $\rho$ ``outsources'' denominators and reduces the reciprocally linear ODE \eqref{xyr} to the linear case \eqref{xylin} of section \ref{LinC2}.
Since complex foliations are not affected by complex $\rho\neq0$, we only have to address $\rho\!=\!0$, and can invoke the previous section \ref{LinC2} anywhere else.
This demonstrates the conceptual power of the language of complex foliations and their holonomy, in the analysis of complex-time blow-up loops and real-time blow-up stars as introduced in definitions \ref{defloop} and \ref{defstar}.

We multiply the reciprocally linear ODE \eqref{xyr} by, and absorb, an Euler multiplier 
\begin{equation}
\label{Eulerr}
\rho(x,y):=(\alpha x+\beta y)(\gamma x+\delta y)\,.
\end{equation}
The complex Euler multiplier $\rho$ does not affect the associated foliation $\omega\!=\!0$, in the pole-free domain $\cD:=\{\rho\neq0\}\subset\C^2$ of ODE \eqref{xyr}, since $\dot x, \dot y$ do not vanish.
In fact,  the foliation $\omega\!=\!0$ of \eqref{xyr} on $\cD$ extends to the foliation of the linear ODE \eqref{xylin}, verbatim, on $\C^2\setminus\{(0,0)\}$.

As in section \ref{LinC2}, we assume nonvanishing determinant $\det$ and positive discriminant $d=\mathrm{tr}^2-4\det$.
This diagonalizes \eqref{xylin} to \eqref{xydiag}, with nonzero simple real eigenvalues $\lambda_1, \lambda_2$ and spectral quotient $\lambda$.
The interested reader may delight in the cases of complex conjugate or multiple eigenvalues, or of complex coefficients, along similar lines.

The discussion of linear holonomy in section \ref{LinC2} then applies verbatim, as far as the extended foliation $\omega\!=\!0$ is concerned.
The discussion of sections \ref{z=e}-\ref{RatFol} concerning the complex flows within leafs $\fL$, however, does involve the new Euler multiplier $\rho$ of \eqref{Eulerr} and requires further comment.
To be specific, but also remain sufficiently general, we consider the diagonalization \eqref{xydiag} and replace \eqref{Eulerr} by an accordingly transformed Euler multiplier 
\begin{equation}
\label{rhor}
\rho=(x-ay)(x-by)\,.
\end{equation}
Indeed, the diagonalization process still allows us to pick real $a,b$, arbitrarily.
We discuss distinct nonzero $a,b$, leaving more degenerate cases to the reader again.
For $\rho=xy$ see the discussion of \eqref{xyrp} below.

In view of real blow-up, we now
compare the systems
\begin{equation}
\label{ODExyr}
\begin{aligned}
     \dot x(t_1)=\rho\, \dot x(t)&= \lambda_1\,x\,,   \\
      \dot y(t_1)=\rho\, \dot y(t)&= \lambda_2\,y\,,
\end{aligned}
\end{equation}
in the shorthand notation of remark \ref{short}, here for $x,y$ in original time $t$ and rescaled time $dt=\rho\,dt_1$\,, respectively.
This is analogous to \eqref{ODExy} and \eqref{tt1}, but now with $\rho$ from \eqref{rhor}.

Assume $\lambda_1<0$, for blow-up.
We consider the trivial stable leaf $\fL\,{:}\ y=0\neq x$ first.
Then $\rho=x^2$ implies $x\dot x(t)= \lambda_1$ on $\fL$, i.e.\ $x^2(t)=2\lambda_1t$.
The resulting double branching of $x$ over $t$ allows us to peruse the discussion of proposition \ref{proploop} and section \ref{z=e} for $m\!=\!3$.
See figure \ref{figm234} (center).
In the associated blow-up star, consider real-valued blow-up (red), for $0>t\nearrow T=0$ and $\pm x(t)$ towards the singularity of ODE \eqref{ODExyr} at $x\!=\!0$.
Then real-time continuation for $T=0\swarrow t>0$ requires purely imaginary blow-down values of $x(t)$ (blue).
In conclusion, we obtain discrepancy-free \emph{minimal blow-up loops $(\gamma^t,\,\gamma^u)$ around blow-up in the stable leaf, which are doubly branched over $t$}.
The arguments for the trivial blow-up or blow-down leaf $\fL\,{:}\ x=0\neq y$ and real $\lambda_2\neq0$ are analogous.

Next we explore the general leaf $\fL^c$ of nonzero $x_0,\,y_0$ for discrepancy-free blow-up loops $(\gamma^t,\,\gamma^{xy})$.
In section \ref{IrrFol} we have already noticed how blow-up loops require rational spectral quotients $\lambda\in\Q\setminus\{0\}$, viz. coprime integer nonzero $\lambda_\iota=-n_\iota$\,. 
We therefore follow section \ref{RatFol}.
For blow-up, we again assume negative eigenvalues, ordered by $0<n_1<n_2$\,, without loss of generality.

Each nontrivial Riemann leaf $\fL^c\,{:}\  x^{n_2}=c\,y^{n_1}$ is then algebraic, with some complex $c\neq0$ and branching at $x\!=\!y\!=\!0$, as before.
Although each fixed leaf $\fL^c$ provides $2(n_2-n_1)$ nontrivial intersections with vanishing Euler multipliers $\rho\!=\!0$ of \eqref{rhor}, these cannot accumulate to the singularity $x\!=\!y\!=\!0$ of interest, within $\fL^c$.
There, the Euler multiplier satisfies $\rho=x^2+\ldots$\,, locally, up to ramified terms of higher order.
Locally and as before, we may therefore invoke proposition \ref{proploop} for $(x,y)$ with $m-1:=2$ and section \ref{RatFol} with winding number $\fw(\gamma^x)=n_1$.
This establishes discrepancy-free \emph{minimal blow-up loops  $(\gamma^t,\,\gamma^{xy})$ with winding number $\fw(\gamma^t)=2n_1$}\,.

More globally, it remains to explore the nontrivial intersections of the algebraic leafs $\fL^c$ with the singularities of ODE \eqref{ODExyr} along the vanishing lines $\rho\!=\!0$ of the Euler multiplier \eqref{rhor}.
By complex rescaling of $x,y$, it is sufficient to address the standardized case $a=1\neq b$ at the intersection $x\!=\!y\!=\!1$ in the standard leaf $\fL:=\fL^1$ of $c\!=\!1$.
The remaining $(n_2-n_1)$-th roots of unity $x\!=\!y$ can be treated analogously, along with the intersections of $\fL$ with $x\!=\!by$.
Let $t_1=t=T=0$ denote the blow-up time associated to the intersection $x\!=\!y\!=\!1$.
The flow \eqref{ODExyr} is regular in $t_1$\,, with $dt=\rho\, dt_1$\,.
Integration therefore provides expansions
\begin{equation}
\label{rhotr}
\begin{aligned}
\rho&=(x-y)(x-by)=(x-y)((1-b)+\ldots)=(1-b)(n_2-n_1)\,t_1 + \ldots\,;\\
t&=\tfrac{1}{2}(1-b)(n_2-n_1)\,t_1^2 + \ldots\,.
\end{aligned}
\end{equation}
As in the case $m\!=\!3$ of the stable leaf, this identifies discrepancy-free \emph{minimal blow-up loops $(\gamma^t,\,\gamma^{xy})$ around blow-up, here at the nontrivial singularities of $\fL^1\cap\{\rho\!=\!0\}$.}
Note winding numbers $\fw_x=1$ of $\gamma^x$ around $x\!=\!1$, and $\fw_t=2$ of $\gamma^t$ around $T\!=\!0$.

\section{Main results on linearization of real and complex blow-up}\label{B}

In this section, as inspired by Masuda, we develop our main results on circumventions, in complex time, of real-time blow-up for polynomial systems \eqref{ODExy} of degree $m$.
More specifically, our results aim for blow-up loops, in complex time, and their associated real-time blow-up stars.
See definitions \ref{defloop}, \ref{defstar} and our main result, theorem \ref{thmloop} in section \ref{PoiMas} below.
Our results are based on analytic ``linearization at infinity'', i.e.\ on analytic linearization of ODEs \eqref{ODEuz} or \eqref{ODEvw} at equilibria $(u,z)=(0,e)$ or $(v,w)=(0,1/e)$, in projective coordinates $u\!=\!1/x,\ z\!=\!y/x$ or $v\!=\!1/y,\ w\!=\!x/y$, and in the sense of definition \ref{defequi}.
See section \ref{Fol} and in particular proposition \ref{propuvw}.

Real-time blow-up then refers to blow-up $u\rightarrow 0$ for $0>t\nearrow T=0$ in original time $t$.
In particular, we aim for discrepancy-free blow-up loops $(\gamma^t,\,\gamma^{uz})$ around blow-up which close perfectly, possibly after loops $\gamma^t\subset\C\setminus\{T\}$ with several windings, locally around blow-up time $t\!=\!T\!=\!0$.
We distinguish rescaled time $t_1$ of $u,z$ from original time $t$ of $x,y$, with Euler multiplier $dt=\rho\,dt_1$\,.
See remark \ref{short} for shorthand notation.

Real-time blow-up requires convergence of $(u,z)$ to equilibria $u\!=\!0,\ z\!=\!e$ of \eqref{ODEuz}; see section \ref{LoopStar}.
We therefore assume real eigenvalues $\lambda_1<0\neq\lambda_2$\,, with spectral quotients $\lambda=\lambda_1/\lambda_2\neq0,1$.
Shifting $e$ to $z\!=\!0$ and linear diagonalization leads to local expansions
\begin{equation}
\label{fuz}
\begin{aligned}
     \dot u( t_1)&=f^u(u,z):=u(\lambda_1+\ldots)\,,   \\
     \dot z( t_1)&=f^z(u,z):=z(\lambda_2+\ldots)+u^2(c_2+\ldots)\,,
\end{aligned}
\end{equation}
with an Euler multiplier $\rho(u)=u^{m-1}$ of $dt=\rho\, d t_1$\,.
See \eqref{uzdiag} and \eqref{tt1}.

In section \ref{Ls} we assume a stable leaf $\fL\subset\{z\!=\!0\}$ of the foliation $\omega\!=\!0$ associated to \eqref{fuz}.
This allows us to proceed in analogy to linear section \ref{z=e}; see theorem \ref{thmMs}.
In sections \ref{ResLin}-\ref{SieLin}, our approach to blow-up loops will be based on local \emph{analytic diagonalization}: we seek biholomorphisms $(u,z)=\Psi(\tilde u,\tilde z)$ which transform nonresonant nonlinear systems \eqref{fuz} to their \emph{diagonal} flow or orbit equivalent linearization
\begin{equation}
\label{tuzdiag}
\begin{aligned}
    \dot {\tilde u}(t_1) &= \lambda_1\, \tilde u\,,   \\
    \dot {\tilde z}(t_1) &= \lambda_2\, \tilde z\,,
\end{aligned}
\end{equation}
locally near $u\!=\!z\!=\!0$.
See definition \ref{defequi} for the notion of local analytic flow equivalence.
Background on analytic linearization and the prerequisite spectral nonresonance conditions on the spectral quotient $\lambda:=\lambda_1/\lambda_2$ is provided in sections \ref{ResLin}, \ref{PoiLin}, and \ref{SieLin}, by references to \cite{Ilya}.
See in particular theorems \ref{thmPoi} and \ref{thmsadlin}.
Our main result, theorem \ref{thmloop} below, then establishes discrepancy-free blow-up loops in the nonresonant Poincaré cases of nodes $0<\lambda,\,1/\lambda\in\Q\setminus\N$.
This includes a description of leaf closures as Riemann surfaces $\cR$ with branching at blow-up.
In section \ref{SieLin} we comment on the quasiperiodic leaves which arise in the complementary Siegel domain of saddles $\lambda<0$, under Diophantine conditions on the irrationality of the spectral quotient $\lambda$.

Why, dear energetic reader, do we present our main results so late in the paper? 
The previous sections are an attempt to carefully build up and explain our unconventional approach, and its pitfalls, by several carefully chosen examples.
The examples may look trivial, at first - until you actually try, from scratch. 
The ``boring'' linear case of section \ref{LinC2}, in particular, outsources part of the proof for the nonlinear case.
Skip all that -- and enter at your own risk.

\subsection{Trivial stable leaf}\label{Ls}

\begin{thm}\label{thmMs}
Assume a trivial stable leaf $\fL\subset\{z\!=\!0\!\neq\!u\}$ exists, associated to the stable eigenvalue $\lambda_1<0$ of the linearization at the equilibrium of $u\!=\!z\!=\!0$ of \eqref{ODEuz}, \eqref{fuz}.\\
Then any local loop $\gamma^u$ of $u\in\fL\setminus\{0\}$ with winding number $\fw_u$=1 defines a minimal blow-up loop $(\gamma^t,\,\gamma^u)$ of \eqref{ODEuz}.
The winding number $\fw(\gamma^t)=m-1$ is determined by the degree $m\geq2$ of the polynomial $f^u(u,0)$.
\end{thm}
\begin{proof}
The stable leaf $\fL$ simplifies \eqref{fuz} to
\begin{equation}
\label{fuzs}
\begin{aligned}
     \dot u( t_1)&=f^u(u,z):=u(\lambda_1+\ldots)\,,   \\
     \dot z( t_1)&=f^z(u,z):=z(\lambda_2+\ldots)\,.
\end{aligned}
\end{equation}
The dynamics on the stable leaf is given by
\begin{equation}
\label{ufs}
	\dot u( t_1)=f^u(u,0)\,.
\end{equation}
Without loss of generality, we may invoke Poincaré and assert analytic flow linearization of \eqref{ufs} by a local scalar biholomorphism $\Psi(u)=u+\ldots$\,, to obtain
\begin{equation}
\label{uflin}
	\dot u( t_1)=\lambda_1\,u\,,
\end{equation}
albeit with an Euler multiplier $\rho(\Psi(u))=u^{m-1}(1+\ldots)$ instead of pure $\rho=u^{m-1}$.
See \cite{ArnoldODE,Ilya} for this scalar variant of analytic flow linearization; the vector case will be detailed in sections \ref{PoiLin}, \ref{PoiMas} below.
The local expansion \eqref{tuexpan} of original time $t$ remains valid, now in terms of $u$ alone.
In particular the loop $\gamma^t$ closes whenever $\gamma^u$ does.
This proves that $(\gamma^t,\,\gamma^u)$ is indeed a blow-up loop.
Hence we may invoke \eqref{wtu} of proposition \ref{proploop}, after all, to identify the minimal winding number $\fw(\gamma^t)=m-1$, as claimed.
This proves the theorem.
\end{proof}

Up to the flow linearizing biholomorphism $\Psi$, blow-up stars follow as in the discussion of the stable leaf for the linear case; see section \ref{z=e}.

\subsection{Resonances and analytic diagonalization}\label{ResLin}

Going back as far as the dissertation of Henri Poincaré, analytic flow diagonalization is based on \emph{integer nonresonance}.
See section I.5 of \cite{Ilya} and the references there, for more detailed accounts of the classical linearization results summarized here and in sections \ref{PoiLin}, \ref{SieLin}, \ref{Dim} below.
\begin{defi}\label{defPoiSie}
A complex pair $\boldsymbol{\lambda}:=(\lambda_1,\lambda_2)\in(\Co)^2$ belongs to the \emph{Poincaré domain}, if the closed straight line interval 
\begin{equation}
\label{interval}
[\lambda_1,\lambda_2]:=\{(1-t)\lambda_1+t\lambda_2\,|\, 0\leq t\leq1\} 
\end{equation}
between $\lambda_1$ and $\lambda_2$\,, alias the convex hull of $\{\lambda_1,\lambda_2\}$, does not contain zero.
If $0\in[\lambda_1,\lambda_2]$, in contrast, then $\boldsymbol{\lambda}$ belongs to the complementary \emph{Siegel domain}.\\
Let $\boldsymbol{\alpha}=(\alpha_1,\alpha_2)$ denote pairs of nonnegative integers $\alpha_\iota$\,, such that $|\boldsymbol{\alpha}|:=\alpha_1+\alpha_2\geq2$.
The pair $\boldsymbol{\lambda}$  is called \emph{resonant of order} $|\boldsymbol{\alpha}|$, if for some $\boldsymbol{\alpha}$ and $\iota=1,2$
\begin{equation}
\label{res}
\lambda_\iota=\boldsymbol{\alpha}\cdot\boldsymbol{\lambda}:=\alpha_1\lambda_1+\alpha_2\lambda_2\,.
\end{equation}
In absence of any resonance, the pair $\boldsymbol{\lambda}$ is called \emph{nonresonant}.
\end{defi}

Resonance in the Poincaré domain is easy to determine.

\begin{prop}\label{propPoires} 
Let $\boldsymbol{\lambda}:=(\lambda_1,\lambda_2)\in(\Co)^2$ and $\lambda:=\lambda_1/\lambda_2\in\Co$. 
Then 
\begin{enumerate}[(i)]
	\item $\boldsymbol{\lambda}$ is in the Poincaré domain, if and only if either $\lambda\in\C\setminus\R$ is nonreal, or else $\lambda>0$.
	\item Nonreal quotients $\lambda\in\C\setminus\R$ are always nonresonant.
	\item Real $\lambda>0$ are resonant, if and only if $\lambda$ or $1/\lambda$ is in $\N\setminus\{1\}$. 
\end{enumerate}
\end{prop}

\begin{proof}
Claim (i) follows from definition \ref{defPoiSie}.
The requirement $|\boldsymbol{\alpha}|\geq2$ implies claim (ii).
Indeed, nonreal spectral quotients $\lambda$ then always violate the real resonance condition \eqref{res}, because one of the two $\alpha_\iota$ has to vanish, while the other has to exceed 1.
The same condition prevents resonances for semisimple $\lambda\!=\!1$, and for a real collinear solution of resonance equation \eqref{res} with $\lambda>0$.
This proves claim (iii) and the proposition. 
\end{proof}

In the following sections, we will base our search for blow-up loops on local analytic diagonalization. 
Definition \ref{defequi} then asks for a local biholomorphism $(u,z)=\Psi(\tilde u,\tilde z)$:
\begin{equation}
\label{Psi}
\begin{aligned}
u&=\Psi^u(\tilde u,\tilde z):=\tilde u+\ldots\,,\\
z&=\Psi^z(\tilde u,\tilde z):=\tilde z+\ldots\,,
\end{aligned}
\end{equation}
which linearizes the flow of system \eqref{fuz} to become diagonal; see \eqref{tuzdiag}.
Any such diagonalization requires semisimple nonresonant spectrum $\boldsymbol{\lambda}=(\lambda_1,\lambda_2)$, in general.
In this context we speak of \emph{spectral pairs} $\boldsymbol{\lambda}$ with \emph{spectral quotient} $\lambda:=\lambda_1/\lambda_2\in\Co$.
For real spectrum $\boldsymbol{\lambda}$, the saddle case of $\lambda<0$ falls into the Siegel domain.
The Poincaré domain $\lambda>0$ leads to real-time stable or unstable nodes .
We summarize available results, separately, for the Poincaré node domain $\lambda>0$ and for the Siegel saddle domain $\lambda<0$.

\subsection{The Poincaré domain}\label{PoiLin}

We start with the Poincaré domain of spectral quotients $\lambda\neq0$ which are nonreal, or positive. See proposition \ref{propPoires}.
The following theorem is essentially due to Poincaré.
See \cite{ArnoldODE} and, for a complete proof, section I.5B in \cite{Ilya}. 

\begin{thm}\label{thmPoi}
Assume the spectral pair $\boldsymbol{\lambda}:=(\lambda_1,\lambda_2)$ of eigenvalues $\lambda_1,\lambda_2\in\Co$ is in the Poincaré domain and nonresonant.\\
Then \eqref{fuz} possesses a locally analytic flow linearization $\Psi$ at $u\!=\!z\!=\!0$; see \eqref{Psi}.
In the semisimple case, the linearization is diagonal.
\end{thm}

As discussed in section \ref{IrrFol}, already, and again in section \ref{SieLin} below for the Siegel domain, Poincaré linearization of blow-up ODE \eqref{ODEuz} in the case of irrational spectral quotients $0\neq\lambda\in\C\setminus\Q$ does not lead to complex circumvention of real-time blow-up.
Instead, we have to rely on discrepancy-free blow-up loops $(\gamma^t,\,\gamma^u)$ in the trivial stable leaf $\fL\subset \{z\!=\!0\}$ of $\lambda_1<0$.

\subsection{Minimal blow-up loops and blow-up stars in the Poincaré domain}\label{PoiMas}

Sections \ref{C2m}, \ref{Fuji}, and \ref{Con} essentially apply the following main result on blow-up loops and blow-up stars.
See definitions \ref{defloop} and \ref{defstar}.

\begin{thm}\label{thmloop}
Consider any equilibrium $(u,z)=(0,e)$ of the blow-up equation \eqref{ODEuz}, for $f,g$ of polynomial degree $m\geq2$.
For the linearization at $(0,e)$, we assume integer negative coprime eigenvalues $\lambda_1=-n_1,\ \lambda_2=-n_2$ with nonresonant, but rational, Poincaré spectral quotient $\lambda=n_1/n_2>0$.
In other words $0<\lambda,1/\lambda\in (\Q\setminus\Z)\cup\{1\}$.
In case $\lambda\!=\!1$ we assume $\lambda_1=\lambda_2=-1$ to be semisimple, i.e.\ without Jordan block of length 2. \\
Then there exist discrepancy-free minimal blow-up loops $(\gamma^t,\,\gamma^{uz})$, in blow-up variables $u\!=\!1/x,\ z\!=\!y/x$ of \eqref{ODEuz} and in original time $t$ of \eqref{ODExy}, locally around blow-up of $(x(t),y(t))$ at $t\!=\!T\!=\!0$ and $(u,z)=(0,e)$, with the following properties.

\begin{enumerate}[(i)]
  \item The local leaf $\fL$ of $\gamma^{uz}$ closes to a Riemann surface $\cR=\fL\cup\{0,e\}$ with branching point $(u,z)=(0,e)$.
  \item In a suitable regular parametrization $\Theta{:}\ \{|\theta|\leq1\}\rightarrow \cR$ by the complex unit disk, we obtain complex analytic expansions
  \begin{align}
\label{at}
   t-T & =a_t\,\theta^{(m-1)n_1}+\ldots\\
\label{au}
   u &= a_u\,\theta^{n_1}+\ldots 
\end{align}
with nonzero coefficients $a_t\,,\, a_u$\,, such that $\gamma^{uz}=\Theta(\S^1)$.
  \item The associated winding numbers of the projected closed loops $\gamma^t$ and $\gamma^u$ are
\begin{align}
\label{wt}
    \fw_t&=(m-1)\,n_1>0\,,   \\
\label{wu}
    \fw_u&=n_1>0\,.  
\end{align}
  \item The real-time blow-up star in the local leaf $\fL\subset\cR$ is parametrized, alternatingly, by the $\fw_t=(m-1)n_1$ real analytic blow-up and blow-down branches of $\theta$ in \eqref{at} such that $t\leq T$ and $t\geq T$, respectively.
\end{enumerate}
Analogous results hold for blow-up variables $v,w$ replacing $u,z$.
Reversing time, the statements remain true for blow-down at positive $\lambda_\iota=n_\iota$\,.
\end{thm}

\begin{proof}
We give an outline first.
We work in the normalized coordinates of \eqref{fuz} with ${\mathbf{F}}=(f^u,f^z)$, and $z\!=\!e$ shifted to $e\!=\!0$ with diagonalized linear part.
For the requirement of semisimple $\lambda\!=\!1$ see remark \ref{remnonss}.
Nonresonance proposition \ref{propPoires} allows us to invoke the Poincaré local analytic flow diagonalization \eqref{tuzdiag} for $(\tilde u,\tilde z):=\Psi^{-1}(u,z)$; see theorem \ref{thmPoi}.
In particular the resulting foliation $\omega\!=\!0$ becomes algebraic with nontrivial local leafs $\tilde\fL^c{:}\  \tilde u^{n_2}=c\,\tilde z^{n_1}\neq0,\ c\neq0$\,; compare \eqref{Ln1n2}.
Then section \ref{RatFol} applies, in time $t_1$\,.
With $\lambda_\iota=-n_\iota$\,, this shows that solution loops $\gamma^{\tilde u \tilde z}$ of $(\tilde u(t_1),\tilde z(t_1))\in\tilde\fL^c$ in imaginary time $t_1=-\mi s,\ 0\leq s\leq 2\pi$, close up after minimal period $s=2\pi$ and with $\fw_{\tilde u}=n_1$\,.
Consider the loop $\gamma^{uz}:=\Psi(\gamma^{\tilde u \tilde z})$ on the local leaf $\fL:=\Psi(\tilde\fL)$.
Direct integration, as in section \ref{RatFol} and proposition \ref{proploop}, will then establishes a minimal blow-up loop $(\gamma^t,\gamma^{uz})$ and an associated blow-up star with properties (i)--(iv).

Indeed, local analytic flow diagonalization by $\Psi$ as in \eqref{Psi} identifies the loop $\gamma^{uz}:=\Psi(\gamma^{\tilde u \tilde z})$ on the local leaf $\fL:=\Psi(\tilde\fL)$.
However, $\Psi=(\Psi^u,\Psi^z)$ modifies the Euler multiplier $\rho(u)=u^{m-1}$ associated to \eqref{ODEuz} to become
\begin{equation}
\label{rhotil}
\tilde \rho(\tilde u,\tilde z)=(\Psi^u)^{m-1}={\tilde u}^{m-1}(1+\ldots)\,.
\end{equation}
We implicitly claim here, and prove in \textbf{step 1} below, that $\tilde u$ can be factored out, in the first component $\Psi^u$ of the near-identity analytic flow linearization \eqref{Psi}:
\begin{equation}
\label{Psiu}
\Psi^u=\tilde u\cdot(1+\ldots)\,.
\end{equation}
Factorization preserves invariance of the $z$-axis $\{u\!=\!0\}$.
Locally, in other words,
\begin{equation}
\label{zaxis}
u\!=\!0 \ \Leftrightarrow\ \tilde u\!=\!0 \ \Leftrightarrow\ \tilde\rho\!=\!0\,.
\end{equation}
To establish $\gamma^t,\gamma^{uz}$ as a blow-up loop, we have to address the parametrization of the loop  $\gamma^{uz}:=\Psi(\gamma^{\tilde u \tilde z})$ by original time $dt=\tilde\rho\,dt_1$\,.
We already noticed how closed loops $\gamma^{\tilde u \tilde z}$, for the diagonally linear ODE \eqref{uzlin} in variables $(\tilde u,\tilde z)$, arise in time $t_1=-\mi s,\ 0\leq s\leq 2\pi$.
See  section \ref{RatFol} and in particular \eqref{uzsol}, i.e.~$\tilde u(t_1)=\tilde u_0\,\theta^{n_1},\ \tilde z(t_1)=\tilde z_0\,\theta^{n_2}$ in the \emph{regularizing parameter} $\theta:=\exp(-t_1)$ of our theorem.
We consider any nontrivial local leaf $\tilde\fL$, with small nonzero $\tilde u_0\,,\, \tilde z_0$\,, and $|\theta|\leq1$.
The loop $\gamma^{\tilde u \tilde z}$ in $t_1=-\mi s$ then possesses minimal period $s=2\pi$, only, and lifts to a closed loop $\gamma^{uz}$, by $\Psi$.
In \textbf{step 2} below, we show that the lifted Euler multiplier $dt=\tilde\rho\, dt_1$ of \eqref{rhotil} generates a closed loop $s\mapsto\gamma^t(s):=t(s)\in\Co$.
This shows that $(\gamma^t,\gamma^{uz})$ is a blow-up loop, up to steps 1 and 2 below.
We postpone minimality of $(\gamma^t,\gamma^{uz})$ to our proof of claim (iii).

\begin{description}
\item[Step 1: Proof of factorization claim \eqref{Psiu}.]\hfill\\[2mm]
Analytic flow linearization consists of formal diagonalization, by formal nonresonant normal form theory as in section 2.4 of \cite{Vdb}, and a proof of convergence, by majorants as in \cite{Ilya}, section I.5B.
To prove factorization \eqref{Psiu} it is therefore sufficient to scrutinize the formal normal form process, which proceeds by induction over the degrees $n$ of the successive transformations $\Psi$ and the successively transformed vector fields.
The transformed vector field $\tilde {\mathbf{F}}$ of ${\mathbf{F}}$, in variables $(\tilde u,\tilde z)$, is given by the pull-back
\begin{equation}
\label{Psi*F}
\tilde {\mathbf{F}} = \Psi^*{\mathbf{F}} := (D\Psi)^{-1}\cdot {\mathbf{F}}\circ\Psi\,,
\end{equation}
where $D\Psi$ denotes the Jacobi matrix.
Note that $\Psi^*{\mathbf{F}}$ is linear in ${\mathbf{F}}$.

To start the induction, at the linear level $n\!=\!1$,  we note that the linear part $A=D\mathbf{F}(0)$ of the blow-up equation \eqref{ODEuz} is lower triangular.
Under our spectral assumptions, linear diagonalization only requires linear lower triangular $\Psi$ of the form \eqref{Psiu}.

For the induction step from $n{-}1$ to $n$, we assume that some linearly nonresonant ${\mathbf{F}}$ has been linearized, up to degree $n{-}1$, by a transformation $\Psi$ of degree $n{-}1\geq1$.
For the $u$-components, we also assume the factorizations
\begin{align}
\label{fuu}
    f^u&=u\cdot f_1^u(u,z)\,,   \\
\label{Psiuu}
    \Psi^u&=\tilde u\cdot \Psi_1^u(\tilde u,\tilde z)\,,
\end{align}
with $\Psi_1^u=1+\ldots$\ .
Here \eqref{fuu} expresses the general assumption that ODE \eqref{fuz} preserves the invariant complex plane $u\!=\!0$, and \eqref{Psiuu} is the induction hypothesis \eqref{Psiu} on $\Psi$.

The induction step in the formal normal form procedure is based on the action of $\mathrm{ad} A\,\mathbf{G}:=A\mathbf{G}(\mathbf{x})-D\mathbf{G}(\mathbf{x})A\mathbf{x}$ on homogeneous vector polynomials $\mathbf{G}$ of degree ${n\geq2}$ in $\mathbf{x}=(u,z)$.
Let $\mathbf{e}_\iota$ denote the unit vectors in $\C^2$.
Here $A=D\mathbf{F}(0)=\mathrm{diag}\,[\lambda_1,\lambda_2]$ is the linear part of $\mathbf{F}$ which we have diagonalized from the start.
Then $\mathrm{ad} A$ is diagonalizable with eigenvalues $\lambda_\iota-\boldsymbol{\alpha}\cdot\boldsymbol{\lambda}
$ and eigenvectors $u^{\alpha_1}z^{\alpha_2}\mathbf{e}_\iota$ of polynomial degree $n=|\boldsymbol{\alpha}|$; compare the resonance condition \eqref{res}.
By our nonresonance assumption, $\mathrm{ad} A$ is invertible.
The normal form induction step consists of elimination of all $n$-th order terms $u^{\alpha_1}z^{\alpha_2}\mathbf{e}_\iota$ of ${\mathbf{F}}$ by corresponding terms 
\begin{equation}
\label{normal}
(\lambda_\iota-\boldsymbol{\alpha}\cdot\boldsymbol{\lambda})^{-1}\cdot u^{\alpha_1}z^{\alpha_2}\mathbf{e}_\iota
\end{equation}
in $\Psi$.
This step, by the way, is the source of small denominators in the Siegel domain.
By assumption \eqref{fuu}, the component $f^u$ does not contain any pure $z$-terms with $\alpha_1=0$, for $\iota\!=\!1$.
Therefore the normal form transformation $\Psi^u$ does not contain any pure $z$-terms, either.
This completes the induction step to degree $n$, for $\Psi^u$ in \eqref{Psiuu}.

It remains to show that pull-back $\Psi^*{\mathbf{F}}$ under $\Psi$ satisfying \eqref{Psiuu} preserves the factorization \eqref{fuu} of $f^u$. 
Indeed, $D\Psi$ and $D\Psi^{-1}$ in \eqref{Psi*F} are lower triangular, at $\tilde u\!=\!0$.
Therefore
\begin{equation}
\label{fuu=0}
\tilde f^{\tilde u} = (\Psi_{\tilde u}^{\tilde u})^{-1}\cdot f^u(\tilde u,\tilde z) =  \tilde u\cdot (\Psi_{\tilde u}^{\tilde u})^{-1}\,f_1^u(\tilde u,\tilde z)=0
\end{equation}
at $\tilde u\!=\!0$.
Here we used subscripts $\tilde u$ to indicate partial derivatives with respect to $\tilde u$.
Since \eqref{fuu=0} perpetuates \eqref{fuu} from $f^u$ to $\tilde f^{\tilde u}$, i.e.\ to degree $n$, this completes the induction step of the normal form procedure.
Convergence of the resulting power series for $\Psi$, by \cite{Ilya}, then proves claim \eqref{Psiu}.

\item[Step 2: Proof that $\gamma^t\subset\C\setminus\{0\}$ is a closed loop.]\hfill\\[2mm]
We have to evaluate the effect of the analytic linearization $\Psi$ on time $t\!=\!t(s)$ along $t_1\!=\!-\mi s,\ 0\leq s\leq 2\pi$, via the Euler multiplier $dt=\tilde\rho\, dt_1$ of \eqref{rhotil}.
Fix small $(\tilde u_0\,,\,\tilde z_0)\in\tilde\fL^c$.
Along $\tilde u(s)\!=\!\tilde u_0\, \theta^{n_1},\ \tilde z(s)\!=\!\tilde z_0\, \theta^{n_2}$  with $\theta\!=\!\exp(\mi s),\ d\theta\!=\!\mi\theta\,ds$, we integrate and substitute \eqref{Psiu}:
\begin{equation}
\label{tloop}
\begin{aligned}
   \int_{t(0)}^{t(2\pi)} dt &= \int_{t_1(0)}^{t_1(2\pi)} \tilde\rho\, dt_1=-\mi\int_0^{2\pi} \big(\Psi^u(\tilde u(s),\tilde z(s))\big)^{m-1}\,ds    \\
    &= -\oint \big(\Psi^u(\tilde u_0 \theta^{n_1},\tilde z_0 \theta^{n_2})\big)^{m-1}\,d\theta/\theta =\\
    &= -\oint \tilde u_0^{m-1}\, \theta^{(m-1)n_1}(1+\ldots)^{m-1}\,d\theta/\theta =0\,.
\end{aligned}
\end{equation}
Here we have applied Cauchy's theorem to the last contour integral over $\theta\!:=\!\exp(\mi s)\in\S^1$.
Since $(m-1)n_1\geq1$, indeed, the integrand is analytic in $|\theta|\leq1$, for fixed small $|\tilde u_0|,|\tilde z_0|$.
This completes \textbf{step 2}, and establishes the nontrivial blow-up loop $(\gamma^t,\gamma^{uz})$.
\end{description}

It remains to address claims (i)--(iv).
\begin{enumerate}[(i)]
  \item follows from section \ref{RatFol}, for $\tilde{\cR}:=\tilde \fL^c\cup\{(0,0)\}$ and $\cR:=\Psi(\tilde\cR)$.
  \item is based on the map $(u,z)=\Theta(\theta):=\Psi(\tilde u_0\,\theta^{n_1},\tilde z_0\,\theta^{n_2})$. 
  Insertion of factorization \eqref{Psiu} shows expansion \eqref{au}. 
  Expansion \eqref{at} follows from insertion of $\Theta(\theta)$ and \eqref{Psiu} in expansion \eqref{tuexpan} from the proof of proposition \ref{proploop}.
  \item follows from expansions \eqref{at}, \eqref{au} on the closed loops $\gamma^t,\gamma^u$. 
  The winding number \eqref{wt} does not depend on our choice of $(\tilde u_0\,,\,\tilde z_0)\in\tilde\fL^c$. 
  It does not even depend on the choice of any nontrivial small $|\tilde u_0|,|\tilde z_0|>0$. 
  Therefore the blow-up loop $(\gamma^t,\gamma^{uz})$ is minimal.
  \item follows from expansion \eqref{at} for complex $|\theta|\leq 1$.
  Indeed, the local biholomorphism $\tilde\theta:=\theta\cdot(1+\ldots)^{1/((m-1)n_1)}$ straightens the real analytic blow-up star on $\cR$ to become radial and alternating in $\tilde\theta$, just as in the linear case.
\end{enumerate}
This proves the theorem.
\end{proof}

In the setting of theorem \ref{thmloop}, we obtain expansions of $(t,u(t),z(t)$ with respect to $\theta$, viz.~the complex root $(t-T)^{1/\fw_t}$; see \eqref{at}.
We call these expansions the associated \emph{(finitely branched) meromorphic blow-up expansion} of $x\!=\!1/u,\,y\!=\!z/u$, in terms of $(t-T)^{1/\fw_t}$. 
In the setting of proposition \ref{proploop}, we also obtain an integer leading exponent $t-T\sim u^{m-1}+\ldots$, because $\fw_t=(m-1)\fw_u$\,.

Irrational spectral quotients $0<\lambda\not\in\Q$ are also nonresonant, of course.
Analytic flow linearization at nodes shows how open sets of initial conditions reach blow-up, in finite real time. 
Our discussion of the linear case in section \ref{IrrFol} detected how nontrivial leafs of the foliation $\omega\!=\!0$ keep winding, with quasiperiodic holonomy multipliers $\exp(2\pi\mi\lambda)$. 
See remark \ref{remirr}.
Although such quasiperiodic winding occurs arbitrarily near blow-up itself, it never quite closes up to form a discrepancy-free blow-up loop.
This also prevents finite meromorphic branching, and generates essential singularities instead.
Also in the nonlinear case, quasiperiodic holonomy allows us to achieve arbitrarily small discrepancy gaps, upon iterating the same local loop $\gamma^t\subset\C\setminus\{T\} $ around blow-up time $t\!=\!T\!=\!0$, over and over.
Indeed, the solutions $(u(t),z(t))$ along iterated loops $\gamma^t$ \emph{almost} close up -- with arbitrarily small remaining relative errors, alias discrepancies, in the projective coordinates $u\!=\!1/x,\ z\!=\!y/x$.
This dichotomy between rational and irrational spectral quotients $\lambda>0$ extends our linear remark \ref{remirr} to the linearizable polynomial case.

\subsection{Blow-up in the Siegel domain}\label{SieLin}

In our setting \eqref{fuz}, the Siegel domain consists of spectral quotients $\lambda<0$, i.e.\ of saddles $\lambda_1<0<\lambda_2$\,.
We locally normalize the stable manifold, i.e.\ the stable holomorphic leaf $\fL$, to become trivial: $\fL\subset\{z\!=\!0\}$.
This allows us to define the associated (local) \emph{nonlinear holonomy} $h$ of $u$ over $z$ analogously to the linear holonomy of $x$ over $y$ in \eqref{linhol}.

We have seen in section \ref{LinHol} how linear systems imply linear holonomy.
The following adaptation of \cite{Ilya}, theorem 22.7 to our present notation provides a partial converse, in the saddle case. 
For terminology, see definition \ref{defequi}.

\begin{thm}\label{thmsadlin} 
Assume $\lambda<0$ and linearity for the saddle holonomy $h$ of $u$ over $z$.\\
Then ODE \eqref{fuzs} is analytically orbit equivalent to its linearization at $u\!=\!z\!=\!0$, under a local biholomorphism \eqref{Psi}.
The same holds true for any holonomy which is analytically equivalent, locally, to its own linearization at the fixed point $h(0)=0$.
\end{thm}

For general nonlinear holonomy, analytic flow diagonalization in the Siegel domain of saddles $\lambda_1<0<\lambda_2$ based on nonresonance is classical.
All rational spectral quotients $0>\lambda=\lambda_1/\lambda_2\in\Q$ are resonant, and all irrational spectral quotients $\lambda\in\R\setminus\Q$ are nonresonant.
Analytic flow diagonalization, however, requires additional \emph{Diophantine conditions}, which ensure $\lambda\not\in\Q$ to be  ``sufficiently irrational'', quantitatively.
See the surveys in \cite{Milnor}, §11, in \cite{Ilya}, section I.5E, and \cite{Bern} more recently, for further details and references.
 
A sufficient condition for analytic flow diagonalization is Siegel's classical \emph{small divisor} or \emph{small denominator condition}
\begin{equation}
\label{Dio}
|\lambda_\iota-\boldsymbol{\alpha}\cdot\boldsymbol{\lambda}|\geq \delta\, |\boldsymbol{\alpha}|^{-N}\,,
\end{equation}
for some $\delta,N>0$ and uniformly in $|\boldsymbol{\alpha}|\geq2$.
See \cite{Siegel}.
Lebesgue almost all real pairs $\boldsymbol{\lambda}$ satisfy such a condition.
On the other hand, local linearizability fails for a generic set of $\lambda<0$ and rational holonomies of resonance degree at least 2 \cite{Milnor}.
Also note $|\boldsymbol{\alpha}\cdot\boldsymbol{\lambda}|< |\boldsymbol{\alpha}|^{-2}$, by continued fraction expansion of irrational spectral quotients $\lambda\in\R\setminus\Q$ \cite{ArnoldODE}.
Condition \eqref{Dio} has been celebrated in Kolmogorov-Arnold-Moser (KAM) theory on the survival of invariant tori under small non-integrable perturbations of sufficiently smooth, integrable real Hamiltonian systems \cite{Moser, ArnoldEnc}.
Symplectic discretizations can serve as such perturbations.
An optimally weakened condition which is still sufficient for analytic flow diagonalization goes back to Bryuno (aka Bruno, Brjuno) \cite{Bruno}:
\begin{equation}
\label{Bryuno}
|\lambda_\iota-\boldsymbol{\alpha}\cdot\boldsymbol{\lambda}|\geq \delta\, \exp(-|\boldsymbol{\alpha}|^{1-\eps})\,,
\end{equation}
for some $\eps,\delta>0$ and uniformly in $|\boldsymbol{\alpha}|\geq2$.
The \emph{Bryuno condition} is optimal for irrational linearization, in the sense that its violation gives rise to many counterexamples \cite{Yoccoz,P-M}.

Analytic linearization of the associated holonomy is local, in a maximal complex region called a \emph{Siegel disk}.
These are important components of the Fatou set, i.e.\ complementary to the Julia set.
See \cite{Milnor} again for a discussion of the boundary of Siegel disks.

After linearization, alas, we realize how we cannot locally reach blow-up saddles in real time, except along the stable leaf $z\!=\!0$; see section \ref{z=e}.
Still, our discussions of the linear case in remark \ref{remirr}, and of the node case $\lambda>0$ in section \ref{PoiMas}, indicate how nontrivial leafs of the foliation $\omega\!=\!0$ keep winding, with quasiperiodic holonomy multipliers $\exp(2\pi\mi\lambda)$, arbitrarily near blow-up itself.
As announced, solutions $(u(t),z(t))$ along iterated Masuda-style blow-up loops $\gamma^t$ around blow-up at $t\!=\!T$ then keep closing \hbox{\emph{almost} --} but not quite.

\section{Example: Homogeneous polynomial foliations}\label{C2m}

As a first illustration of our main results from section \ref{B}, we now study blow-up of planar ODEs \eqref{ODEF}, \eqref{ODEfg} for $m$-homogeneous complex polynomial vector fields $\mathbf{F}=(f,g)$.
We consider degrees $m\geq2$.
Following the Masuda PDE paradigm, we explore minimal blow-up loops, in complex time, and their associated real-time blow-up stars. 
See definitions \ref{defloop} and \ref{defstar}.
Theorem \ref{thmloop} reduces this task to ``linearization at infinity'', i.e.\ to linearization of the projective version \eqref{ODEuz} at equilibria $(u,z)=(0,e_j)$; see theorem \ref{thmuzmlin} below.
We assume \emph{nondegeneracy}:
\begin{equation}
\label{nondegm}
f(x,y)\!=\!g(x,y)\!=\!0\quad \Rightarrow \quad x\!=\!y\!=\!0\,.
\end{equation}
The assumption prevents radial complex lines $\sigma\cdot(x,y),\ \sigma\in\Co$, of nontrivial equilibria. 

For the linear case $m\!=\!1$ see section \ref{LinC2}.
Nondegeneracy \eqref{nondegm} amounted to $\det\neq0$, there.

For polynomial degrees $m\geq2$, the linearization at the trivial equilibrium vanishes.
However, we still obtain self-similar foliations $\omega\!=\!0$ of ODE \eqref{ODEfg} on $\C^2\setminus \{(0,0)\}$, analogously to the linear case.
Since we are interested in blow-up, we use the coordinates \eqref{xyz} of projective compactification.
From \eqref{f12}, \eqref{f21}, we obtain the scaled polynomials
\begin{equation}
\label{f12m}
\begin{aligned}
   &f_1(z) :=f(1,z)\,,\quad  &f_2(w):= f(w,1)=w^mf_1(1/w)\,,   \\
   &g_1(z) :=g(1,z)\,,\quad  &g_2(w):= g(w,1)=w^mg_1(1/w)\,.   
\end{aligned}
\end{equation}
The transformed ODEs \eqref{ODEuz}, \eqref{ODEvw} then become
\begin{equation}
\label{ODEuzm}
\begin{aligned}
    \dot u( t_1)&= -u\,f_1(z)\,,  & \\
    \dot z( t_1)&= -z\,f_1(z)+g_1(z) =: P(z)\,; &\qquad\qquad\qquad\ 
\end{aligned}
\end{equation}
\vspace{-4.5mm}
\begin{equation}
\label{ODEvwm}
\begin{aligned}
   \dot v( t_2)&=-v\,g_2(w)\,,   \\
    \dot w( t_2)&=-w\,g_2(w)+f_2(w)=: P_1(w)=w^mP(1/w)  \,.
\end{aligned}
\end{equation}
As in remark \ref{short}, we continue to distinguish rescaled times $t_1,\,t_2$ of $u,z$ and $v,w$, respectively, from original time $t$ of $x,y$, by shorthand notation.
Note how $m$-homogeneity in $(x,y)$ induces skew product structures over $z\!=\!y/x\!=\!1/w$, and linearity in $u,v$, respectively.
Our nondegeneracy assumption \eqref{nondegm} implies that $f_1, g_1$ do not share any zeros, nor do $f_1, P$.
The same nondegeneracy holds among the triple $f_2,g_2,P_1$\,.

In the linear case $m\!=\!1$ of ODE \eqref{xylin}, for example, we had $f_1=\alpha+\beta z$ and $g_1=\gamma+\delta z$.
Zeros $z\!=\!e_j$ of $P=-zf_1+g_1$ indicated the slopes $z\!=\!y/x$ of eigenvectors, and $-f_1$ their associated eigenvalues $\lambda_\iota$\,.
Nondegeneracy \eqref{nondegm}, alias $\det\neq0$, forbids kernel vectors $(1,z)$.

Let $m\geq2$.
It will be sufficient to discuss ODE \eqref{ODEuzm}.
By nondegeneracy assumption \eqref{nondegm}, the only equilibria are $u\!=\!0,\ z\!=\!e_j$ with $P(e_j)=0$.
For a global discussion of the trivial leaf $\fL\subset \{u\!=\!0\}$ at infinity with dynamics $\dot z(t_1)=P(z)$, under generic assumptions on the polynomial $P$, see section \ref{Scalar} above.
Section \ref{Global} and \cite{FiedlerShilnikov,FiedlerYamaguti} discuss further details, including finite-time blow-up of $z$ on the Riemann sphere $z\in\widehat{\C}=\CP^1$.

Over equilibria $z\!=\!e_j$ of $\dot z(t_1)=P(z)$, we obtain trivial cylindrical leafs $\fL\subset\{z\!=\!e_j\}$ with linear dynamics $\dot u(t_1)=\lambda_1u$.
Here $\lambda_1=-f_1(e_j)\neq0$, by nondegeneracy assumption \eqref{nondegm}.
In particular, $x( t_1)=1/u( t_1)$ grows exponentially, in complex sectors $\Re (\lambda_1 t_1)<0$.
This implies self-similar blow-up of $x(t)=1/u(t)$ at slope $y/x\!=\!z\!=\!e_j$\,, in finite time $t$.
Indeed, \eqref{tt1} with Euler multiplier $\rho=\rho_1=u^{m-1}$ implies
\begin{equation}
\label{ tblowup}
\frac{dt}{dt_1}=u^{m-1}=(u_0)^{m-1\,}\exp((m-1)\lambda_1 t_1)\,,
\end{equation}
for radial $ t_1\rightarrow\infty$ in sectors $\Re(\lambda_1 t_1)<0$.

Among all nontrivial leafs of the foliation \eqref{ODEuzm}, scaling by constant $\sigma\in\Co $ provides holomorphically equivalent flows.
By linearity of \eqref{ODEuzm} in $u$, the holonomy of $u$ over $z$ around $z\!=\!e_j$ is linear.
The linearization of \eqref{ODEuzm} at $u\!=\!0,\ z\!=\!e_j$ possesses semisimple eigenvalues $\lambda_1=-f_1(e_j)\neq 0$ and $\lambda_2=P'(e_j)$.
This leads to the following analytic linearization result.

\begin{thm}\label{thmuzmlin}
Assume nondegeneracy \eqref{nondegm} and consider any simple zero $z\!=\!e_j$ of $P(z)$.\\
Then the linearization of \eqref{ODEuzm} at any equilibrium $(u,z)=(0,e_j)$ is diagonal.
Explicitly, the eigenvalues are given by
\begin{align}
\label{l1m}
    \lambda_1&=-f_1(e_j)\neq0\,,   \\
\label{l2m}
    \lambda_2&=P'(e_j)\neq0 \,.
\end{align}
The holonomy of $u$ over $z$ near $z\!=\!e_j$ is also linear, with holonomy multiplier $\exp(2\pi\mi\lambda)$ determined by the spectral ratio 
\begin{equation}
\label{specm}
\lambda=\lambda_1/\lambda_2=-f_1(e_j)/P'(e_j)\,.
\end{equation}
We assume
\begin{enumerate}[(i)]
\item
either, Siegel spectrum of \emph{saddle type}, i.e.\ negative spectral ratio 
$\lambda=\lambda_1/\lambda_2<0$,
\item
 or else, nonresonant Poincaré spectrum of \emph{node type}, i.e.\ with spectral ratio 
$0<\lambda,1/\lambda\in (\Q\setminus\Z)\cup\{1\}$.
\end{enumerate}
In the nonresonant Poincaré case (ii) with real eigenvalues $\lambda_\iota$\,, we obtain local analytic \emph{flow} linearization of the polynomial ODE \eqref{ODEuzm} to diagonal form, in the sense of definition \ref{defequi}.
The complex case (ii) reduces to the real case after a complex rotation in time.
The complex nonresonant Poincaré case (ii), and the general Siegel case (i), only support local analytic \emph{orbit} diagonalization.
\end{thm}

\begin{proof}
Claims \eqref{l1m}, \eqref{l2m} follow from ODE \eqref{ODEuzm}.
Indeed, the linearization of the polynomial system \eqref{ODEuzm},  at the equilibrium $u\!=\!0,\ z\!=\!e_j$\,, is diagonal with eigenvalues $\lambda_1,\, \lambda_2$\,.

For real spectral quotients $\lambda$, we may invoke a constant complex scaling of time, if necessary, such that $\lambda_1<0\neq\lambda_2\in\R$ are both real.
Let $z(s)=\gamma(s):=\eps\exp(\mi s)\subset\C\setminus P^{-1}(0)$ denote a small circle $\gamma=\gamma^z$ around $z\!=\!e_j$ in the complex $z$-plane.
Then linearity of \eqref{ODEuzm} in $u$ implies linearity of the local holonomy in $u$, along closed curves in the homotopy class of $\gamma^z$.
In the saddle case (i), i.e.\ in the Siegel domain, theorem \ref{thmsadlin} implies analytic local linearization of the foliation, i.e.\ up to orbit equivalence.

For real eigenvalues $\lambda_\iota$ with nonresonant positive spectral quotients $\lambda$ of alternative assumption (ii), i.e.\ in the Poincaré domain, local analytic flow linearization follows from  theorem \ref{thmPoi}.
Complex time rotation only preserves orbits.
This proves the theorem.
\end{proof}

With linearization theorem \ref{thmuzmlin} at hand, the results on irrational and rational foliations of sections \ref{LinC2}, \ref{B}, and their discussions, apply to the $m$-homogeneous case.
Note how the Siegel case of saddles with irrational, and hence nonresonant, spectral quotients $0>\lambda\in\R\setminus\Q$ is linearizable by theorem \hbox{\ref{thmuzmlin} (i) --} without any small divisor conditions like \eqref{Dio} or \eqref{Bryuno}.

Consider the rational Poincaré case (ii) of nonresonant nodes $\lambda>0$ next.
For discrepancy-free minimal blow-up loops around blow-up in trivial stable leafs $\fL=(\Co)\times\{z\!=\!e_j\}$ see theorem \ref{thmMs}.
For minimal blow-up loops in nontrivial leafs $\fL$ around blow-up at nonresonant nodes $u\!=\!0,\ z\!=\!e_j$\,, see theorem \ref{thmloop}.
Note how the semisimple diagonal linearization includes the nonresonant case $\lambda\!=\!1$.

In the present radially self-similar case, the trivial stable leafs of theorem \ref{thmMs} can be treated directly, of course.
Indeed, $\dot u(t_1)=\lambda_1u=-f_1(e_j)u$ becomes $\dot u(t)=\lambda_1 u^{2-m}$, with $\lambda_1<0$.
For minimal blow-up loops $(\gamma^t,\,\gamma^u)$ in the stable leaf $z\!=\!e_j$\,, separation of variables then recovers the explicit relation
\begin{equation}
\label{tuexpli}
\lambda_1(m-1)\cdot(t-T)=u^{m-1}\,.
\end{equation}
This confirms proposition \ref{proploop}; see \eqref{wtu}, \eqref{tuexpan}.
Our very first Riccati example \eqref{Ric0} corresponded to $m\!=\!2,\ \lambda_1=-1$.

\section{Example: Polynomial Hamiltonians, and the pendulum}\label{Ham}

As a second class of examples, we now study blow-up for Hamiltonian systems in the complex plane, with polynomial Hamilton functions $H$.
Again, we aim for minimal blow-up loops and their real-time blow-up stars; see definitions \ref{defloop} and \ref{defstar}.
We discuss the second order pendulum equation, in section \ref{Hampen}, and $(m+1)$-homogeneous Hamiltonians, in section \ref{Hamm}, as extreme cases.
Hamiltonian cases, however, defy our standard approach by ``linearization at infinity'', from section \ref{B}.
Instead, section \ref{Hamalg} discusses foliations into leafs $\fL$, alias energy surfaces $\{H\!=\!c\}$ or algebraic curves, in the projective coordinates of sections \ref{Fol}, \ref{Hamhom}.
Similarly to section \ref{PoiMas} in spirit, but not in technical detail, minimal blow-up loops and blow-up stars will then result from closure $\fL\subset\cR$ as a Riemann surface which is branched over an equilibrium ``at infinity''.
See theorems \ref{thmuvwp} and \ref{thmHm}.

We recall that Hamiltonian vector fields $(f,g)$ with respect to the standard symplectic form $\Omega=dx\wedge dy$ are defined by exactness $\omega\!=\!dH$ of the associated Pfaffian foliation $\omega\!=\!0$.
See \eqref{omega}-\eqref{ODEH}.
As a consequence, any \emph{Hamiltonian level surface}, or \emph{energy surface}, $H=H(x,y)=c\in\C$ is invariant under the flow of ODE \eqref{ODEH}, in complex time $t$.
In other words, any leaf $\fL$ of the complex foliation $\omega\!=\!0$ is a connected component of some energy surface $\{H\!=\!c\}$, with the critical points of $H$ removed.
Suppose $c$ is a regular value of $H$, i.e.\ $\mathrm{grad}\,H(x,y)\neq0$ on $\{H\!=\!c\}$.
Then any leaf $\fL\subseteq\{H\!=\!c\}$ is a connected complex 1-manifold. i.e.\ a Riemann surface $\fL=\cR$.
In blow-up, however, our main interest are critical points of $H$ ``at infinity''.

\subsection{Hamiltonian leafs in projective coordinates}\label{Hamhom}

In this section we rewrite Hamiltonian systems \eqref{ODEH} in the projective coordinates of section \ref{Fol}.
To obtain polynomial vector fields \eqref{ODEH} of degree $m$, we assume $H$ is a polynomial of degree $m+1$.
We can then compactify the level sets $\{H\!=\!c\}$, individually, in projective coordinates $[\xi\!:\!\eta\!:\!\zeta]=[x\!:\!y\!:\!1]=[1\!:\!z\!:\!u]=[w\!:\!1\!:\!v]$ on $\CP^2$, much as in \eqref{xyz} and proposition \ref{propuvw}.
Specifically, $H=\sum_{j,k} H_{jk}\,x^jy^k$ with $0\leq j+k\leq m+1$ lifts to $(m{+}1)$-homogeneous $\fH{:}\ \CP^2\rightarrow\CP^1$, with
\begin{equation}
\label{fH}
\fH(\xi,\eta,\zeta):=-c\,\zeta^{m+1}+\sum_{j,k} H_{jk}\,\xi^j\eta^k\zeta^{m+1-j-k}\,.
\end{equation}
The preserved level set $H^{xy}:=H-c=0$ extends projectively as
\begin{align}
\label{Hxy}
    H^{xy}(x,y)&:=\fH(x,y,1)=H(x,y)-c=0\,,   \\
\label{Huz}
    H^{uz}(u,z)&:=\fH(1,z,u)=0\,,   \\
\label{Hvw}
    H^{vw}(v,w)&:=\fH(w,1,v)=0\,.
\end{align}
More explicitly, with $c^{xy}:=c,\ c^{uz}:=u^{m+1}c,\ c^{vw}:=v^{m+1}c$\,,
\begin{equation}
\label{cuzvw}
\begin{aligned}
 H^{xy}(x,y)&= \phantom{u^{m+1}}H(x,y)-c	^{xy}\,,\qquad\quad \  \\
 H^{uz}(u,z)&=u^{m+1}H(1/u,\,z/u)-c^{uz}\,, \\ 
 H^{vw}(v,w)&=v^{m+1}H(w/v,\,1/v)-c^{vw}\,.
\end{aligned}
\end{equation}
In other words, we have compactified the level set $\{H\!=\!c\}\subset\C^2$ to $\{\fH\!=\!0\}\subset\CP^2$.
Analogously to proposition \ref{propuvw}, we can now study the Hamiltonian structure in projective coordinates, and in the shorthand notation of remark \ref{short}.

\begin{prop}\label{propuvwH} 
The leafs $\fL\subseteq\{H\!=\!c\}$ of the foliation $\omega\!=\!0$ defined by the Hamiltonian vector field $f\!=\!H_y\,,\ g\!=\!-H_x$ of $H$ with respect to the symplectic form $\Omega^{xy}:=dx\wedge dy$, extend to leafs $\fL\subset\{\fH\!=\!0\}\subset\CP^2$ at energy level zero as follows.
With the same Euler multipliers $\rho$ and the same time transformation as in proposition \ref{propuvw}, the transformed vector fields \eqref{ODExy}-\eqref{ODEvw} remain Hamiltonian, with respect to $\fH$, appropriate symplectic forms $\Omega$, and up to Euler multipliers $\rho$\,.
Specifically, with respect to the symplectic form $\Omega^{xy}:=dx\wedge dy$ we obtain the Hamiltonian vector field
\begin{equation}
\label{ODExyH}
\begin{aligned}
    \dot x(t)&= H_y^{xy}\,,   \\
     \dot y(t)&= -H_x^{xy}\,.
\end{aligned}
\end{equation} 
In time $t_1$ and at omitted arguments $u,z$, the ODE \eqref{ODEuz} then takes the form
\begin{equation}
\label{ODEuzH}
\begin{aligned}
     \tfrac{1}{u}\,\dot u(t_1)&= -H_z^{uz}=-f_1\,, \qquad\quad  \\
     \tfrac{1}{u}\,\dot z(t_1)&= H_u^{uz}= \tfrac{1}{u}(g_1-zf_1)\,,
\end{aligned}
\end{equation}
at the energy levels $H^{xy}\!=\!0\!=\!H^{uz}$.
In particular the system is Hamiltonian with respect to the Hamiltonian function $H^{uz}$ and the symplectic form $\Omega^{uz}:=-u^{-1}du\wedge dz$.
Analogously, $\Omega^{vw}:=-v^{-1}dv\wedge dw$ in time $t_2$ recasts \eqref{ODEvw} at energy levels $H^{xy}\!=\!0\!=\!H^{vw}$ in Hamiltonian form as
\begin{equation}
\label{ODEvwH}
\begin{aligned}
     \tfrac{1}{v}\,\dot v(t_1)&= -H_w^{vw}=-g_2\,,   \\
     \tfrac{1}{v}\,\dot w(t_1)&= H_v^{vw}= \tfrac{1}{v}(f_2-wg_2)\,.
\end{aligned}
\end{equation}
\end{prop}

\begin{proof}
We only check the second ODE of \eqref{ODEuzH}.
All other claims follow analogously.
By definition \eqref{FH} we obtain the Hamiltonian ODE for $\dot z(t_1)$, associated to the Hamiltonian function $H^{uz}(u,z)$ and with respect to the symplectic form $\Omega^{uz}$ as follows: 
\begin{equation}
\label{Huzu}
\begin{aligned}
H_u^{uz}(u,z)&:=\tfrac{\partial}{\partial u}\big(u^{m+1}\left(-c+H(1/u,\,z/u)\right)\!\big)=\\
&\phantom{:}=(m+1)u^m\left(-c+H(x,y)\right)-u^{m-1}\left(H_x^{xy}(x,y)+zH_y^{xy}(x,y)\right)=\\
&\phantom{:}=\tfrac{1}{u}\left(u^mg(1/u,\,z/u)-zu^mf(1/u,\,z/u)\right)=\\
&\phantom{:}=\tfrac{1}{u}\left(g_1(u,z)-zf_1(u,z)\right)=\tfrac{1}{u}\,\dot z(t_1)\,.
\end{aligned}
\end{equation}
on the invariant leaf $\fL\subseteq\{H\!=\!c\}=\{H^{xy}\!=\!0\}=\{H^{uz}\!=\!0\}$.
We have also used \eqref{ODEH} and, transiently, $x\!=\!1/u,\ y\!=\!z/u$.
This proves the second ODE of \eqref{ODEuzH}.
Analogous checking establishes the Hamiltonian structure \eqref{ODEuzH}, \eqref{ODEvwH} of the blow-up systems \eqref{ODEuz}, \eqref{ODEvw} from proposition \ref{propuvw}, as claimed.
\end{proof}

\begin{rem}
We have deviated from the standard approach to coordinate transformations by pull-back of $\Omega^{xy}=dx\wedge dy$ in $\Omega(\mathbf{F}^H,\cdot)=dH$ of \eqref{FH}, e.g. under $(x,y)=(1/u,\,z/u)=:\Psi(u,z)$.
For example, $\Psi^*\Omega^{xy}=u^{-3}dz\wedge du=u^{-2}\,\Omega^{uz}$ and $\Psi^*H^{xy}=H\circ\Psi=u^{-m-1}H^{uz}$.
In particular, $\Psi$ is not symplectic.
A deeper reason is that homogeneously projectivized $\fH,\fF$ do not define vector fields on $\CP^2$; see \eqref{ODEfF}, \eqref{fxfy}.
There are two reasons for our choices: 
\begin{enumerate}[(i)]
  \item we preserve the polynomial character of the Hamiltonian function, and 
  \item we choose the prefactors $1/u,\,1/v$ of the symplectic forms to generate Euler multipliers which minimize the degrees of the resulting polynomial ODEs \eqref{ODEuzH}, \eqref{ODEvwH}.
\end{enumerate}
\end{rem}

\subsection{Complex algebraic curves}\label{Hamalg}

In this section we briefly relate foliations by algebraic energy surfaces $\{H\!=\!c\}$ to the more standard language of complex algebraic curves.
As in \eqref{Hxy}, consider polynomial Hamiltonian functions $H=H(x,y)$ of degree $m\geq2$.
We now assume the polynomial equation $H\!=\!c$ is algebraically irreducible over $\C$.
If $c$ is a regular value of $H$, then the level set $H\!=\!c$ defines a complex 1-dimensional manifold of $(x,y)\in\C^2$, i.e.\ a Riemann surface $\cR$.
The connected Riemann surface constitutes the leaf $\fL=\{H\!=\!c\}=\cR$, in this regular case.
Finite branching may occur, over $x$ or over $y$, but not over time $t$.
We have seen how to compactify $\cR$ as $\fH=0$ in $\CP^2$.
In fact, any compact Riemann surface arises by such a construction.
See for example \cite{Bries, Jost, Lamotke} for further details.

Our particular interest, however, is in critical values $c$ of $H$, i.e.\ in equilibria of ODE \eqref{ODExyH}, and its projective variants, within the level set $\{H\!=\!c\}$, which is now an algebraic variety.
In the nondegenerate Morse case $\det D^2H(e)\neq0$, in particular, the algebraic variety $\{H\!=\!c\}$ contains the stable and unstable leafs.
In general, the details of the local ODE flow, in original complex time, depend on the details of the resolution by Puiseux expansions at equilibria; see \cite{Bries}.
A simple example arises in the next section.

\subsection{The pendulum and hyperelliptic curves}\label{Hampen}

Our first class of Hamiltonian examples are second order pendulum equations.
This leads to classical elliptic and hyperelliptic curves.
The hyperelliptic case is of interest because it features totally degenerate ``linearizations at infinity''.
In particular, this fails all results of section \ref{B}, and motivates our approach via foliations by energy surfaces $\{H\!=\!c\}$, alias complex algebraic curves.
See theorem \ref{thmuvwp}.
For comparison of language and elementary illustration, we include the classical torus case \eqref{penW}, \eqref{loopW} of Weierstrass elliptic functions, and even the cylinder case \eqref{ODExyp1}, \eqref{Hlin} of the linear pendulum.

For polynomial nonlinearities $g$, the standard pendulum $\ddot x-g(x)=0$ with Newtonian force law $g(x)$ at location $x$ and velocity $y$ reads
\begin{equation}
\label{pen}
\begin{aligned}
    \dot x &=y\,,   \\
    \dot y &=g(x)\,.  
\end{aligned}
\end{equation}
In other words, the pendulum ODE is Hamiltonian with respect to the polynomial Hamiltonian function
\begin{equation}
\label{Hpen}
H(x,y):=\tfrac{1}{2}y^2-G(x)\,.
\end{equation}
Here $G$ is a primitive function of $g$, i.e.\ $G_x=g$.
In coefficients: $G=\sum_{j=0}^{m+1}G_{m+1-j}\,x^j$ with $G_0\neq0$, and $g=G_x=\sum_{j=0}^m g_{m-j}\,x^j$ with $g_{m-j}=(j+1)G_{m-j}$\,.
Note the critical points $g(e)=0,\ y=0$ of $H$, alias equilibria $(x,y)=(e,0)$ of the pendulum \eqref{pen}.
For polynomial degrees $m=\deg\,g\geq4$, the compact Riemann surfaces $\cR$ associated to regular values $H\!=\!0$ possess genus $[m/2]\geq2$ and are called \emph{hyperelliptic curves}.
We comment on the torus case $m\!=\!2$ of quadratic $g$ and Weierstrass elliptic functions below.

Consider any equilibrium $H(e,0)=c$.
Without loss of generality, we may assume $c\!=\!0$, subsuming $c$ into the potential $G$.
Locally, we obtain branched energy surfaces $\cR=\{H\!=\!0\}$.
Let $r$ denote the multiplicity of the zero $e$ of $g$.
Shifting $e$ to $e\!=\!0$ we obtain an expansion
\begin{equation}
\label{y(x)}
\begin{aligned}
y^2&=2G(x)=a^2x^{r+1}(1+\ldots)\,, \ i.e.\\
y\phantom{^2}&=\pm a\,x^{(r+1)/2}\,\sqrt{1+\ldots}\ ,
\end{aligned}
\end{equation}
with some coefficient $a\neq0$.
For odd multiplicities $r$, we obtain two unbranched leafs of $y$ over $x$.
At simple zeros $r\!=\!1$, for example, these constitute the stable and unstable manifold of the saddle $(e,0)$, with resonant spectral ratio $\lambda=-1$; compare section \ref{SieLin}.
For even multiplicities $r$, we obtain a single leaf $y$ which is branched at $x\!=\!e$ of ramification index 2 over $x$.

Towards infinity $[\xi\!:\!\eta\!:\!0]$, we invoke propositions \ref{propuvw} and \ref{propuvwH} for projective coordinates $[x\!:\!y\!:\!1]=[1\!:\!z\!:\!u]=[w\!:\!1\!:\!v]$.
In view of \eqref{Hxy}-\eqref{Hvw}, again on energy level $H\!=\!c\!=\!0$, we therefore work in energy surfaces $H=H^{xy}=H^{uz}=H^{vw}=0$, i.e.\ with energy relations
\begin{align}
\label{xyp}
    y^2&=2G(x)\,,   \\
\label{uzp}
    u^{m-1}z^2&=2u^{m+1}G(1/u)=2\sum_{j=0}^{m+1} G_{m+1-j}\,u^{m+1-j}\,,\\
\label{vwp}
    v^{m-1}&=2v^{m+1}G(w/v)=2\sum_{j=0}^{m+1} G_{m+1-j}\,w^jv^{m+1-j}\,.
\end{align}

\begin{thm}\label{thmuvwp}
Consider the energy level $c\!=\!0$ of pendulum Hamiltonians $H^{xy}=H=\tfrac{1}{2}y^2-G(x)$, with $f\!=\!y$ and polynomial force law $g\!=\!G_x$ of degree $m\geq2$, i.e.\ $G_0\neq0$. \\
Then the following claims hold true.
\begin{enumerate}[(i)]
  \item The ODEs \eqref{ODEuzH} and \eqref{ODEvwH} become
\begin{equation}
\label{ODEuzp}
\begin{aligned}
     \dot u(t_1)&= -u\,H_z^{uz}=-u^mz\,,   \\
     \dot z(t_1)&= \phantom{-}u\,H_u^{uz}= u^mg(1/u)-u^{m-1}z^2\,;
\end{aligned}
\end{equation}
\begin{equation}
\label{ODEvwp}
\begin{aligned}
     \dot v(t_2)&= -v\,H_w^{vw}=-v\,v^mg(w/v)\,,  \qquad\ \  \\
     \dot w(t_2)&=  \phantom{-}v\,H_v^{vw}= v^{m-1}-w\,v^mg(w/v)\,.
\end{aligned}
\end{equation}
  \item The only equilibrium at infinity, i.e.\ for projective coordinates $[\xi\!:\!\eta\!:\!0]$, occurs at $v\!=\!w\!=\!0$, i.e.\ in the closure of the local leaf or leafs $\fL\subset\{H^{vw}=0\}$.
For degrees $m\geq3$, the linearization of \eqref{ODEvwp} vanishes at $v\!=\!w\!=\!0$.
In particular, the linearization results of section \ref{B} do not apply.
  \item Local expansions of $\fL$ in a complex parameter $\theta$ take the analytic form
\begin{equation}
\label{vth}
v=a\,\theta^{m+1}(1+\ldots)\,,\quad w=:\theta^{m-1}\,, \quad a:=\pm(2G_0)^{1/(m-1)}\,,
\end{equation}
where we omit higher order terms in $\theta$.
Original time $t$ of \eqref{pen} is related to $\theta$ and $w$ by 
\begin{equation}
\label{dtdw}
dt=-\tfrac{2}{m-1}(1+\ldots)\,dw\,.
\end{equation}
  \item In particular, small loops of $\theta\in\Co$ around $0$ define discrepancy-free minimal blow-up loops $(t,v,w)\in(\gamma^t,\,\gamma^{vw})$ around blow-up time $t\!=\!T\!=\!0$ and $v\!=\!w\!=\!0$ such that
\begin{equation}
\label{loopH}
\begin{aligned}
     t(\theta)&\phantom{:}=-\tfrac{2}{m-1}\,\theta^{m-1}(1+\ldots)\,,\\
    v(\theta)&\phantom{:}= \quad\,\phantom{-}a\,\theta^{m+1}(1+\ldots) \,,  \\
    w(\theta)&:= \quad\,\phantom{-a\,}\theta^{m-1}\,.
\end{aligned}
\end{equation}
  \item For even polynomial degrees $m\geq2$, \eqref{vth} defines a single leaf.
The blow-up loop $(\gamma^t,\,\gamma^{uz})$ closes up, first, after one cycle of $\theta$, i.e.\ with winding numbers $\fw_t=m-1=\fw_w$ and $\fw_v=m+1$ of $\gamma^t,\,\gamma^w$ and $\gamma^v$, respectively.\\
For odd $m$, where \eqref{vth} defines two leafs $\pm v$, half a cycle of $\theta$ is sufficient to close the blow-up loop of $(t,v,w)\in(\gamma^t,\,\gamma^{vw})$, on each leaf, with half those winding numbers.
  \item For either parity of $m$, the closure $\cR:=\fL\cup\{v\!=\!w\!=\!0\}$ of any local leaf $\fL$ is a Riemann surface with branch point $v\!=\!w\!=\!0$.
The ramification indices over $v$ and $w$ are given by the winding numbers $\fw_v$ and $\fw_w$\,, respectively.
  \item For either parity of $m$ in (v), the blow-up stars in the local leafs $\fL\subset\cR$ are parametrized, alternatingly, by the $\fw_t=\fw_w$  real analytic blow-up and blow-down branches of $\theta$ in \eqref{loopH}, such that $t\leq T$ and $t\geq T$, respectively.
\end{enumerate}
\end{thm}

\begin{proof}
Claim (i) follows from proposition \ref{propuvwH}, with the substitutions \eqref{f12}, \eqref{f21} of proposition \ref{propuvw} applied to $f\!=\!y$ and $g\!=\!g(x)$.

To prove claim (ii), we first note that $u\!=\!z\!=\!0$ is not an equilibrium of \eqref{ODEuzp}.
Indeed, $u\!=\!z\!=\!0$ implies $\dot z=g_0=(m+1)G_0\neq0$.
Therefore we only have to discuss equilibria $v\!=\!0,\,w\!=\!e_j$ of ODE \eqref{ODEvwp}.
Then $0=-\dot w=w\,v^mg(w/v)=g_0\,e_j^{m+1}=0$ implies $e_j=0$, which also satisfies energy relation \eqref{vwp}.
For $m\geq3$, the linearization of \eqref{ODEvwp} at $v\!=\!w\!=\!0$ is the zero matrix, i.e.\ totally degenerate.

Concerning claim (iii), we first prove \eqref{vth}.
We proceed via the Newton polygon of energy relation \eqref{vwp} at $v\!=\!w\!=\!0$.
Puiseux expansion identifies the leading terms $v^{m-1}=2\,G_0\,w^{m+1}+\ldots$ of
\begin{equation}
\label{vwpui}
v^{m-1}=2\,G_0\,w^{m+1}+\ldots+2\,G_{m+1}v^{m+1}\,.
\end{equation}
For $G_0\neq0$, in general, this defines a local leaf, or leafs, with 
$(m+1):(m-1)$ torus knots, up to a reducing factor 2 in case of odd degrees $m$; see also section \ref{RatFol}.
We introduce the complex covering variable $\theta^{m-1}\!=\!w$, i.e.\ a root of $w$.
Substitution into \eqref{vwpui} and taking roots yields the local expansion \eqref{vth}, but with a polynomial of higher order terms still involving not only $\theta$, but also $v$.
The negative sign of $a$ only appears for odd $m$.
For even $m$, the exponents $m\pm1$ are coprime, and we may pick any root $a$, at positive sign.
Solving \eqref{vth} for $v=v(\theta)$ by the implicit function theorem eliminates any higher order terms in $v$.
The omitted terms now depend on $\theta$, only, and have become an analytic power series, as claimed in \eqref{vth}.

Expansion \eqref{vth} implies the following time transformation:
\begin{equation}
\label{dtdwpf}
\begin{aligned}
dt &= v^{m-1}dt_2 = \left(v^{m-1}/\dot w(t_2)\right)\, dw = \left(v^{m-1}/( v^{m-1}-w\,v^mg(w/v))\right) dw =\\
    &= \left(1- vw\, g(1/(a \theta^2+\ldots))\right)^{-1} dw =\\
    &= \left(1- a\theta^{2m}(1+\ldots)\, g\left((a \theta^2)^{-1}(1+\ldots)\right)\right)^{-1} dw =\\
    &= (1-a^{1-m}g_0+\ldots)^{-1}\,dw = (1-g_0/(2G_0)+\ldots)^{-1}\,dw = -\tfrac{2}{m-1}(1+\ldots)\,dw\,.
\end{aligned}
\end{equation}
Here we have successively used \eqref{tt2}, \eqref{ODEvwp}, and expansion \eqref{vth} with $w=\theta^{m-1}$.
Re-inspection shows the validity for either parity of $m$.
This proves expansion \eqref{dtdw} and completes claim (iii).

Claim (iv) follows from expansion (iii).
The omitted terms of higher order in $\theta$ are borrowed from locally analytic expansions \eqref{dtdw} and \eqref{vth}.
By Cauchy's theorem, the curves $(t,v,w)=(t,v,w)(\theta)$ defined by \eqref{loopH} form a closed minimal blow-up loop.
This proves claim (iv).

To prove claim (v) consider even degrees $m$, first.
Then $m\pm1$ are coprime, and it is sufficient to consider $a:=+(2G_0)^{1/(m-1)}$ in \eqref{vth}.
Indeed, the negative sign can be absorbed into $\theta$, and we obtain a single leaf.
The winding numbers $\fw_\iota$ follow from expansions \eqref{loopH} of the closed loops $\gamma^t,\,\gamma^v,\,\gamma^w$ over the loop $\gamma^\theta\subset\Co$.

On the other hand, suppose $m$ is odd.
Then $(m\pm1)/2$ are coprime, the negative sign of $a$ cannot be absorbed into $\theta$, and we obtain two leafs $\pm v(\theta)$ in \eqref{vth}.
Since $m-1$ is even, we are allowed to replace $m\pm1$ by coprime $(m\pm1)/2$, from regularization \eqref{vth} onward.
In other words, we may introduce $\tilde\theta:=\theta^2$ as a new regularizing parameter.
In particular, half a cycle of $\theta$ corresponds to one full cycle of $\tilde\theta$.
With these substitutions, the previous proof then shows how the previous winding numbers $\fw=m\pm1$ are cut in half, correspondingly, for odd $m$.

Claim (vi) on blow-up stars at blow-up time $T\!=\!0$ follows from time expansion \eqref{loopH}, with the above modification for odd $m$, just as theorem \ref{thmloop}(iv) followed from time expansion \eqref{at}.
This proves the theorem.
\end{proof}

Note that the winding numbers in theorem \ref{thmuvwp}(v) contradict the expectation \eqref{wtu}, here $\fw_t=(m-1)\,\fw_v$\,, of winding numbers suggested by proposition \ref{proploop}.
The vanishing linearization in (ii) is to blame for that, of course.

For an explicit illustration of theorem \ref{thmuvwp}, let us consider the cubic case of \emph{classical elliptic functions}, i.e.\ the Hamiltonian pendulum \eqref{pen}, \eqref{Hpen} with lowest admissible degree $m\!=\!2$ and quadratic nonlinearity
\begin{equation}
\label{penW}
g(x)=6x^2-1\,.
\end{equation}
This corresponds to the Weierstrass elliptic case with $H=y^2-(4x^3-\tilde g_2x-\tilde g_3)$ and $\tilde g_2=2$.
The discriminant $d$ of the polynomial $G=2x^3-\tfrac{1}{2}\tilde g_2x-\tfrac{1}{2}\tilde g_3$ is given by $d=\tilde g_2^3-27\tilde g_3^2$\,.
Note $d\neq0$, for $\tilde g_2=2,\ \tilde g_3\neq \pm(2/3)^{3/2}$.
Explicit solutions in complex time are then given by $(x,y)=(\wp(t),\dot\wp(t))$.
For nonzero discriminant $d$, the \emph{Weierstrass elliptic function} $\wp$ is doubly periodic with respect to a suitable period lattice $t\in\Lambda=\langle1,\tau\rangle_\Z\subset \C$, for some lattice parameter $\Im\tau>0$ which turns out purely imaginary for real coefficients $\tilde g_2,\,\tilde g_3$\,.
Blow-up occurs via the cubic pole $(1/t^2,\,-2/t^3)$ of $(x,y)$, at blow-up time $t\!=\!T\!=\!0\ \mod\Lambda$.
In particular any simple small loop $\gamma^t\subset\Co$ around $T\!=\!0$ defines a minimal blow-up loop $(v(t),w(t))\in\gamma^{vw}$ given by
\begin{equation}
\label{loopW}
\begin{aligned}
v=1/y&=-\tfrac{1}{2}t^3+\ldots \,,\\
w=x/y&=-\tfrac{1}{2}t+\ldots \,.
\end{aligned}
\end{equation}
The loop closes up after a single cycle $\fw_t\!=\!1$ of $\gamma^t$ with $\fw_v\!=\!3,\ \fw_w\!=\!1$.
With $G_0\!=\!2,\ a\!=\!4,\ \theta\!=\!-t/2+\ldots$\,, this confirms \eqref{loopH}.
In contrast to \eqref{loopW}, however, our derivation of \eqref{loopH} was only based on $G_0\neq0$.
We did not require any discriminant condition $d\neq0$ which, for example, excludes blow-up leafs associated to finite equilibria.
For comments on the Duffing pendulum of cubic force law $g$, viz.\ the $m\!=\!3$ variant of \emph{Jacobi elliptic functions}, see \eqref{Duf} below.

We have excluded the linear pendulum $m\!=\!1$ in theorem \ref{thmuvwp}.
The linear case 
\begin{equation}
\label{ODExyp1}
\begin{aligned}
     \dot x(t)&= y \,,  \\
     \dot y(t)&=x \,,
\end{aligned}
\end{equation}
is associated to the homogeneously quadratic Hamiltonian
\begin{equation}
\label{Hlin}
H=\tfrac{1}{2}y^2-\tfrac{1}{2}x^2\,.
\end{equation}
Since the elementary transformations of theorem \ref{thmuvwp}(i) remain valid, we may restrict attention to ODE \eqref{ODEuzp} which now reads
\begin{equation}
\label{ODEuzp1}
\begin{aligned}
     \dot u(t)&= -uz\,,   \\
     \dot z(t)&= 1-z^2\,,
\end{aligned}
\end{equation}
in original time.
Indeed, $t_1=t$ for $m=1$.
The driving ODE of directions $z:=y/x$ is just the Riccati equation of figure \ref{figRiccati} in reversed time.
Note and compare the resonances $-1\!:\!1$ of \eqref{ODExyp1} at the saddle $x\!=\!y\!=\!0$, versus $\mp1\!:\!\mp2$ of \eqref{ODEuzp1} at the sink and source nodes $(u,z)=(0,e_\pm)$ with $e_\pm=\pm1$.

Of course, the linear pendulum \eqref{ODExyp1} does not experience blow-up in finite time.
As in section \ref{LinC2}, it is still instructive to look at the associated complex leafs $H\!=\!0,\,1/2$.
The trivial energy level $H\!=\!0$ consists of the stable and unstable eigenspaces $y=\mp x$ associated to the $-1\!:\!1$ resonant eigenvalues $\lambda_1=-1$ and $\lambda_2=+1$ at the saddle $x\!=\!y\!=\!0$.

The energy levels $H\!=\!c\neq0$ all rescale to $H\!=\!1/2$, biholomorphically.
The linear analogue to the Weierstrass function $\wp$ is the hyperbolic sine. 
Indeed $(x,y)=(\sinh t, \cosh t)$, with imaginary period $2\pi\mi$, parametrizes the unique leaf $\fL\subset\{H\!=\!1/2\}$.
Note the explicit $\mp1\!:\!\mp2$ winding of the nontrivial leaf $\fL$ in coordinates $u(t)=1/\sinh t,\ z(t)=\coth t$ at the nodes $(u,z)=(0,\pm 1)$, for $t\rightarrow \pm\infty$, due to the minimal period $\pi\mi$ of the hyperbolic cotangent.
In conclusion, all leafs $\fL$ of ODE \eqref{ODExyp1} are biholomorphically equivalent to the cylinder $\C/2\pi\mi\Z\cong \C\setminus\{0\}$.

A diagonal variant of the reciprocally linear ODE \eqref{xyr} in section \ref{RecC2}, with $\alpha\!=\!\delta\!=\!0$ and $\beta\!=\!\gamma\!=\!1$, reads
\begin{equation}
\label{xyrp}
\begin{aligned}
    \dot x(t) &=1/x\,, \\
    \dot y(t) &=1/y\,.  
\end{aligned}
\end{equation}
Under an Euler multiplier $dt=\rho\,dt_1$\,, with $\rho:=xy$, this leads back to the linear pendulum \eqref{ODExyp1}.
We can write explicit solutions as
\begin{equation}
\label{xyrexpli}
    \tfrac{1}{2}y^2(t) =t= \tfrac{1}{2}x^2(t) +t_0\,.
\end{equation}
The time shift $t_0\in\C$ between $x(t)$ and $y(t)$ constitutes the constant energy level
\begin{equation}
\label{Hlinr}
H=\tfrac{1}{2}y^2-\tfrac{1}{2}x^2=t_0\,.
\end{equation}
Direct blow-up or blow-down involving the singularity saddle $x\!=\!y\!=\!0$ can only occur along the trivial eigenspaces $y=\pm x$, characterized by synchrony $t_0=0$; see\eqref{xyrexpli}.
The minimal blow-up loops $\gamma^t,\gamma^x,\gamma^y$ around finite-time blow-up at $x\!=\!y\!=\!t\!=\!T\!=\!0$  are obvious, with $\gamma^y=\pm\gamma^x$, winding numbers $\fw_x=\fw_y=1$, and $\fw_t=2$.
In particular we observe real-time blow-up stars with two incoming and outgoing branches, each.

The general asynchronous energy level $H\!=\!t_0\neq0$ provides blow-up with associated real-time blow-up stars, separately, at $y\!=\!0$ for $t\!=\!T_y\!=\!0$, and at $x\!=\!0$ for $t\!=\!T_x\!=\!t_0$\,.
The resulting winding numbers $\fw_\tau=2$ for $\tau:=t-T$ and, respectively, $\fw_y=1$ at $T_y=0$ and $\fw_x=1$ at $T_x=t_0$\,, remain unchanged.

\subsection{Nondegenerate (m+1)-homogeneous Hamiltonians}\label{Hamm}

The results on hyperelliptic curves, in the previous section \ref{Hampen}, were based on an extreme imbalance between the degrees $m+1$ in $x$, and 2 in $y$, of the polynomial Hamiltonian $H(x,y)=\tfrac{1}{2}y^2-G(x)$.
At the other extreme, we now address $(m{+1})$-homogeneous Hamiltonians
\begin{equation}
\label{Hm}
H(x,y)=\sum_{j=0}^{m+1}\, h_j x^jy^{(m+1-j)}\,,
\end{equation}
for $m\geq2$.
We first resort to section \ref{B} and attempt to establish minimal blow-up loops and real-time blow-up stars based on projective coordinates and ``linearization at infinity''.
Resonance, however, will fail our attempt to just combine the results of section \ref{C2m} on nonresonant $m$-homogeneous ODEs \eqref{ODEfg}, with section \ref{Hamhom} on general Hamiltonian systems \eqref{FH}.
See propositions \ref{propuvw}, \ref{proploop}, and \ref{propuvwH}.
Instead, we will again apply algebraic methods to the complex foliation by energy levels.

By \eqref{fH}, the $(m{+}1)$-homogeneous Hamiltonian $H$ lifts to $(m{+}1)$-homogeneous
\begin{equation}
\label{fHm}
\fH(\xi,\eta,\zeta):=-c\,\zeta^{m+1}+H(\xi,\eta)\,.
\end{equation}
With the abbreviations $H_1(z):=H(1,z),\ H_2(w):=H(w,1)=w^{m+1}H_1(1/w)$, \eqref{cuzvw} becomes
\begin{equation}
\label{cuzvwm}
\begin{aligned}
 H^{xy}(x,y)&= H(x,y)-c\,, \qquad\quad\\
 H^{uz}(u,z)&=H_1(z)\,-c\,u^{m+1}\,, \\ 
 H^{vw}(v,w)&=H_2(w)-c\,v^{m+1}\,.
\end{aligned}
\end{equation}
We assume the Hamiltonian $H$ to be \emph{nondegenerate}, i.e.
\begin{align}
\label{nondegH0}
h_0,\,h_{m+1}\quad&\neq\quad\  0\,,\\
\label{nondegH1}
H_1(e)=0\quad &\Rightarrow \quad H_1'(e)\neq0\,,\\
\label{nondegH2}
H_2(e)=0\quad &\Rightarrow \quad H_2'(e)\neq0\,.
\end{align}
We first note that assumption \eqref{nondegH0} is equivalent to $H_1(0),\,H_2(0)\neq0$.
Simplicity of zeros of $H_1, H_2$ is assumed in \eqref{nondegH1}, \eqref{nondegH2}.
Note $e\neq0$ by \eqref{nondegH0}.
Since $H_2(w)\!=\!w^{m+1}H_1(1/w)$, condition \eqref{nondegH1} and its variant \eqref{nondegH2} are equivalent, and hence mutually redundant.
Indeed, we may just swap $H_1$ with $H_2$\,, and $e\neq0$ with $1/e$.

In projective blow-up coordinates, the $m$-homogeneous Hamiltonian ODE \eqref{ODExyH} reads
\begin{align}
\label{ODEuHm}
     \dot u(t_1)&= -uH_1'(z)=-uf_1(z)\,,   \\
\label{ODEzHm}
     \dot z(t_1)&= -(m+1)cu^{m+1} = -(m+1)H_1(z)\,;\\[2mm]
\label{ODEvHm}
     \dot v(t_1)&= -vH_2'(w)=-vg_2(w)\,,   \\
\label{ODEwHm}
     \dot w(t_1)&= -(m+1)cv^{m+1} = -(m+1)H_2(w)\,;
\end{align}
on the trivial energy levels $\fH\!=\!H^{xy}\!=\!H^{uz}\!=\!H^{vw}\!=\!0$.
We have already substituted these energy levels in \eqref{ODEzHm}, \eqref{ODEwHm}.
Equilibria of \eqref{ODEuHm}, \eqref{ODEzHm} are $(u,z)\!=\!(0,e)$ at the simple zeros $z\!=\!e\neq0$ of $H_1$.
In fact, $H_2(0)\neq0$ allows us to restrict attention to local blow-up coordinates $(u,z)$.
Simplicity \eqref{nondegH1}  of zeros $z\!=\!e$ of $H_1$ is also equivalent to nondegeneracy assumption \eqref{nondegm} on their associated $m$-homogeneous Hamiltonian ODEs \eqref{ODExyH}.
Indeed, nontrivial equilibria $(x,y)\neq(0,0)$ are equivalent to equilibria $(u,z)$ with $u\neq0$ and, by \eqref{ODEuHm}, \eqref{ODEzHm}, to zeros $H_1(e)\!=\!0$ of higher multiplicity.

The spectral quotients $\lambda=1/(m+1)>0$ of the linearization identify all equilibria $(0,e)$ as nodes with integer Poincaré resonance $1\,{:}\,(m{+}1)$.
Therefore, the linearization theorem \ref{thmuzmlin}(ii) for general $m$-homogeneous ODEs fails to apply in the Hamiltonian case.
Since we are in the Poincaré domain of nodes, and not in the Siegel domain of saddles, linearization theorem \ref{thmsadlin} also fails, even though \eqref{ODEuHm}, \eqref{ODEvHm} are linear in $u,v$, respectively, with linear holonomy.
In theorem \ref{thmHm} we will establish minimal blow-up loops and stars, nonetheless.

To replace linearization, we invoke invariance of the Hamiltonian energy levels.
Indeed, the Hamiltonian $H$ is constant on leafs $\fL\subset \{H\!=\!c\}$.
Biholomorphic scaling by $\sigma\in\Co$ conjugates the flows on energy levels $c\neq0$ to the reference complex algebraic curve $c\!=\!1$.
The only other case is the trivial energy level $H\!=\!c\!=\!0$.
We treat both cases separately.

The case $H\!=\!c\!=\!0$ further simplifies \eqref{ODEuHm} to
\begin{equation}
\label{ODEuzHm0}
\dot u(t_1)= -uH_1'(e)
\end{equation}
on the leaf $\fL$ associated to invariant slope $z\!=\!e$.
In \eqref{nondegH1}, we have assumed $H_1'(e)\neq0$.
In original time $dt=u^{m-1}dt_1=-(H_1'(e))^{-1}u^{m-2}du$, this establishes 
$t$-loops $\gamma^t$ with winding number  $\fw_t=m-1$ around blow-up at $t\!=\!T$,
over simple blow-up loops $\fw_u$ in $u$.

For $H_2'(e)\neq0$ in ODE \eqref{ODEvHm}, analogously, simple blow-up loops in $v$ establish blow-up loops $\gamma^t$ of the same winding.
As we saw in \eqref{Ric0}, already, the Riccati case $m\!=\!2$ presents a biholomorphism between $t$ and $u=1/x\in\C$.
This settles the trivial energy level $H\!=\!0$.

To study the nontrivial energy level $H\!=\!c\!=\!1$, we invoke global invariance of $H^{uz}\!=\!0$ in \eqref{cuzvwm}. 
This identifies that energy level as a branched compact Riemann surface $\cR$, alias a complex algebraic curve
\begin{equation}
\label{cR}
\cR\,{:}\quad u^{m+1}=H_1(z)
\end{equation}
of $u$ over $z\in\CP^1=\widehat{\C}$.
Branching points of ramification index $m+1$ are the simple zeros $z\!=\!e$ of $H_1$ at $u\!=\!0$\,; see nondegeneracy assumption \eqref{nondegH1}.
Note that $1/w\!=\!z\!=\!\infty$ is not a branching point, by variant \eqref{nondegH2} of \eqref{nondegH1}.

In time $t_1$\,, the skew product \eqref{ODEuHm}, \eqref{ODEzHm} of $u$ over $z$ is linear in the holonomy fiber $u$.
The driving flow in $z$ on $\cR$ is the compactification on the Riemann sphere of the polynomial ODE
\begin{equation}
\label{ODEH1}
\dot z(t_1) = -(m+1)H_1(z)\,.
\end{equation}
See section \ref{Scalar} again, with $P:=-(m+1)H_1$\,, and note the raised degree $m+1$ of $H_1(z)=H(1,z)$.

\begin{thm}\label{thmHm}
For $m\geq2$, consider Hamiltonian vector fields \eqref{ODExyH}-\eqref{ODEvwH} at energy level $H\!=\!c\!=\!1$ of $(m{+}1)$-homogeneous Hamiltonians $H\!=\!H(x,y)$.
Assume nondegeneracies \eqref{nondegH0}-\eqref{nondegH2}.\\
Then blow-up or blow-down, say at time $t\!=\!T\!=\!0$, on the closure of any nontrivial leaf $\fL\subset\{H\!=\!1\}$ only occurs at $(u,z)=(0,e)$, where $z\!=\!e$ are  simple zeros of the polynomial $H_1(z)=H(1,z)$.
Any blow-up or blow-down possesses a local minimal blow-up loop $(\gamma^t,\,\gamma^{uz})$ based on any small cycle $\gamma^t\subset\Co$ of winding number $\fw_t=m-1$, with a real-time blow-up star of $\fw_t$ incoming and outgoing branches.
The winding numbers of the loop $\gamma^{uz}$, in agreement with $1\,{:}\,(m{+}1)$ resonance, are  $\fw_u=1$ and $\fw_z=m+1$.
\end{thm}

\begin{proof}
We use $u$ as a regular parametrization of the Riemann surface $\cR$ in \eqref{cR}, locally at any simple zero $z\!=\!e$ of the polynomial $H_1(z)$.
Indeed, nondegeneracy $H_1'(e)\neq0$ and the implicit function theorem locally imply that 
\begin{equation}
\label{z(u)Hm}
z=z(u)=e+\tfrac{1}{H_1'(e)}u^{m+1}\,(1+\ldots)
\end{equation} 
is a graph, analytically over $u$.
Therefore small simple loops $\gamma^u\subset\Co$ lift to loops $\gamma^z\subset\C\setminus\{e\}$ with winding numbers $\fw_z=m+1$.

It remains to determine the winding number $\fw_t$\,.
We proceed similarly to \eqref{dtdwpf}, i.e.
\begin{equation}
\label{dtdupf}
dt = u^{m-1}dt_1 = u^{m-1}\frac{du}{\dot u(t_1)} = -\frac{u^{m-2}}{H_1'(z(u))}\,du\,.
\end{equation}
Here we have used time rescaling \eqref{tt1} and ODE \eqref{ODEuHm}.
Since $H_1'(e)\neq0$, integration of \eqref{dtdupf} yields
\begin{equation}
\label{tuH1}
t=-\tfrac{1}{(m-1)H_1'(e)}\,u^{m-1}(1+\ldots)\,.
\end{equation}
In particular, small simple loops $\gamma^u\subset\Co$ define loops $\gamma^t\subset\C\setminus\{e\}$ with winding numbers $\fw_t=m-1$, and vice versa.
Blow-up stars at blow-up time $T\!=\!0$ follow from time expansion \eqref{tuH1}, just as theorem \ref{thmloop}(iv) followed from time expansion \eqref{at}.
This proves the theorem.
\end{proof}

\section{Case study: Galerkin caricatures of quadratic heat equations}\label{Fuji}

We now return to our main motivation: 
the crude two-mode Galerkin caricature \eqref{ODEb} of blow-up in the spatially inhomogeneous PDE \eqref{PDEb}; see section \ref{C2}.
The system is quadratic, $m\!=\!2$, but neither $m$-homogeneous nor Hamiltonian.
The main results from section \ref{B} on minimal blow-up loops and real-time blow-up stars apply.
They reduce the blow-up analysis to nonresonant analytic ``linearization at infinity'', viz.\ at equilibria in the projective coordinates of sections \ref{Fol} and \ref{Car}.
See theorems \ref{thmMs}, \ref{thmPoi}, and \ref{thmloop}.

In section \ref{Carsy} we address symmetric spatial coefficients $\mathbf{b}(\bxi)=1+(a-2)\cos(2\bxi)$. The Masuda choice $a\!=\!2$ refers to the original Fujita PDE \eqref{PDEw}.
For antisymmetric modifications $\mathbf{b}(\bxi)$ of $a\!=\!0$ see section \ref{Caras}.
A less technical summary is provided in section \ref{Con}.

\subsection{Projective caricatures}\label{Car}

For the quadratic Galerkin caricature \eqref{ODEb}, proposition \ref{propuvw} on projective coordinates  $[\xi\!:\!\eta\!:\!\zeta]=[x\!:\!y\!:\!1]=[1\!:\!z\!:\!u]=[w\!:\!1\!:\!v]\in\CP^2$ generates the triplet of ODE systems
\begin{align}
\label{ODExG}
 \dot x(t) &=f(x,y):=  \qquad \quad x^2 + b_1xy + \tfrac{1}{4}(2+b_2) y^2\,,  \\
\label{ODEyG}
\qquad\qquad\quad \dot y(t) &=g(x,y):=-y + b_1x^2 + (2+b_2)xy  + \tfrac{1}{4}(3b_1+b_3)y^2\,; \\[2mm]
\label{ODEuG}
 \dot u(t_1) &= \quad\ \  -u\left(1+b_1z+\tfrac{1}{4}(2+b_2)z^2\right)\,,  \\
\label{ODEzG}
 \dot z(t_1) &=b_1+(1+b_2-u)z+\tfrac{1}{4}(b_3-b_1)z^2-\tfrac{1}{4}(2+b_2)z^3\,;\\[2mm]
\label{ODEvG}
 \dot v(t_2) &=\quad-v\left(\tfrac{1}{4}(3b_1+b_3) -v+(2+b_2)w+b_1w^2  \right)\,,\\
\label{ODEwG}
 \dot w(t_2) &=\tfrac{1}{4}(2+b_2)+\left(\tfrac{1}{4}(b_1-b_3)+v\right)w-(1+b_2)w^2-b_1w^3\,;
\end{align}
with real coefficients $b_1,b_2,b_3$\,.
The associated Euler multipliers are
\begin{equation}
\label{tt1t2}
dt=u\,dt_1=v\,dt_2\,.
\end{equation}
Consider \emph{symmetric} $\mathbf{b}(\bxi)=1+(a-2)\cos (2\bxi)$, i.e.\ $b_1\!=\!b_3\!=\!0$, and $a:=2+b_2$\,.
We assume $a\neq 0,1$.
The projectivized ODEs \eqref{ODExG}-\eqref{ODEwG} simplify to
\begin{align}
\label{ODExG1}
 \dot x(t) &=\ x^2 + \tfrac{1}{4}a y^2\,,  \qquad\qquad\qquad \qquad\qquad\qquad\\
\label{ODEyG1}
 \dot y(t) &=y(-1 + a x)  \,; \\[2mm]
\label{ODEuG1}
 \dot u(t_1) &= -u\,(1+\tfrac{1}{4}az^2)\,,  \\
 \label{ODEzG1}
 \dot z(t_1) &=z\,(a-1-u-\tfrac{1}{4}az^2)\,; \\[2mm]
 \label{ODEvG1}
 \dot v(t_2) &=\quad v\left( v-aw  \right)\,,\\
\label{ODEwG1}
 \dot w(t_2) &=\tfrac{1}{4}a+vw-(a-1)w^2\,.
\end{align}
Note how only variant \eqref{ODEuG1}, \eqref{ODEzG1} involves cubic terms.

For \emph{antisymmetric modifications} $\mathbf{b}(\bxi)=1-2\cos (2\bxi)+b_1\cos \bxi+b_3\cos(3\bxi)$ of $a=b_2+2=0$, we obtain the simplified triplet
\begin{align}
\label{ODExG2}
 \dot x(t) &=  \qquad \quad x(x + b_1y)\,,  \\
\label{ODEyG2}
 \dot y(t) &=-y + b_1x^2 + \tfrac{1}{4}(3b_1+b_3)y^2\,; \qquad\qquad\quad\   \\[2mm]
\label{ODEuG2}
 \dot u(t_1) &= \qquad -u\,(1+b_1z)\,,  \\
\label{ODEzG2}
 \dot z(t_1) &=b_1-(1+u)z+\tfrac{1}{4}(b_3-b_1)z^2\,; \\[2mm]
 \label{ODEvG2}
 \dot v(t_2) &=-v\left(\tfrac{1}{4}(3b_1+b_3) -v+b_1w^2  \right)\,,\\
\label{ODEwG2}
 \dot w(t_2) &=w\left(\tfrac{1}{4}(b_1-b_3)+v+w-b_1w^2\right)\,. 
\end{align}
Here the only cubic variant is \eqref{ODEvG2}, \eqref{ODEwG2}.
For simplicity, we only address these two simplifications, separately.

\subsection{Spatially symmetric coefficients}\label{Carsy}
Consider the symmetric Galerkin caricature \eqref{ODExG1}-\eqref{ODEwG1} for spatially symmetric coefficients $\mathbf{b}(\bxi)=1+(a-2)\cos(2\bxi)$, first.
Then the space of spatially symmetric PDE solutions $\mathbf{w}(t,\pi-\bxi)=\mathbf{w}(t,\bxi)$ is time-invariant.
In the Galerkin caricature, only the $x$-axis of spatially homogeneous $\mathbf{w}=\mathbf{w}(t)=x(t)$ is spatially symmetric.
Therefore the $x$-axis $\{y\!=\!0\}$ and the $u$-axis $\{z\!=\!y/x\!=\!0\}$ are invariant, in the Galerkin caricature.
On the $x$-axis $\{y\!=\!0\}$, not surprisingly, we recover the purely quadratic Riccati equation \eqref{Ric0} with explicit blow-up.
On the $u$-axis $\{z\!=\!0\}$, \eqref{ODEuG1} and the associated Euler multiplier \eqref{tt1t2} imply $\dot u(t)= \dot u(t_1)/u=-1$.
Our assumption $a\neq0$ implies absence of a trivial blow-up equilibrium $v\!=\!w\!=\!0$ of \eqref{ODEvG1}, \eqref{ODEwG1}.
It is therefore sufficient to study the blow-up coordinates $u,z$ of \eqref{ODEuG1}, \eqref{ODEzG1}.

\subsubsection{The trivial blow-up equilibrium $u=z=0$}\label{Carsy0}

Linearization at $u\!=\!z\!=\!0$ provides the semisimple eigenvalues and spectral quotient
\begin{align}
\label{s0l1}
\lambda_1&=-1\,, \\
\label{s0l2}
\lambda_2&=a-1\,, \\
\label{s0l}
\lambda&=1/(1-a)\,.
\end{align}
The free parameter $a\neq1$ of $\mathbf{b}(\bxi)=1+(a-2)\cos(2\bxi)$ in the spectral quotient $\lambda=1/(a-1)$ covers all real linearization options studied in section \ref{B}, with polynomial degree $m\!=\!2$.

\begin{description}
\item[Case 0: The trivial spatially homogeneous leaf $\fL\,{:}\ z=0$.]\hfill\\[2mm]
In the trivial spatially homogeneous leaf  of the spatially homogeneous blow-up equilibrium $u\!=\!z\!=\!0$, the above arguments constitute a minimal blow-up loop $\gamma^t=-\gamma^u,\ \gamma^z=0$ of winding number $\fw(\gamma^t)=\fw(\gamma^u)=m-1=1$ and a non-spectacular one-in, one-out real-time blow-up star.
See figure \ref{figm234}, left.

Alternatively, we could have applied theorem \ref{thmMs} to the stable leaf $\fL\,{:}\ z=0$ of the equilibrium $u\!=\!z\!=\!0$ with eigenvalue $\lambda_1=-1$; see \eqref{s0l1}.
Or we could have invoked the explicit solution \eqref{Ric0} of \eqref{ODExG1} for $y\!=\!0$, directly.

\item[Case 1: The Poincaré domain $a<1$.]\hfill\\[2mm]
In the \emph{Poincaré domain} $a<1$, we encounter a real stable node $\lambda_2=a-1<0$ and arbitrary spectral quotients $\lambda>0$.
See section \ref{PoiLin}.

\begin{description}
  \item[Case 1.1: Irrational spectral quotients $\lambda>0$.]\hfill\\
Irrational $\lambda$ are automatically nonresonant.
By theorem \ref{thmPoi}, the flow possesses an analytic flow linearization near $u\!=\!z\!=\!0$.
The associated irrational winding of the nontrivial leafs, however, prevents discrepancy-free blow-up loops.
Instead, iterated loops $\gamma^t$ produce almost-loops $\gamma^{uz}$, due to their quasiperiodic winding.
See remark \ref{remirr}.

\item[Case 1.2: Rational spectral quotients $\lambda>0$.]\hfill\\[2mm]
Let $\lambda=n_1/n_2>0$ with coprime positive integers $n_1,n_2$\,.
We rescale time such that $\lambda_1=-n_1,\ \lambda_2=-n_2$\,.
In other words, $a=1-1/\lambda=1-n_2/n_1$\,.
\emph{Nonresonance} only forbids $n_1=1\neq n_2$\,, as well as $n_2=1\neq n_1$\,.
Since the eigenvalues are semisimple, theorem \ref{thmloop} applies to the stable node case, with degree $m\!=\!2$.
We obtain minimal blow-up loops $(\gamma^t,\,\gamma^{uz})$, in original time $t$ and variables $u\!=\!1/x,\ z\!=\!y/x$ of \eqref{ODExG1}, \eqref{ODEyG1}, around blow-up at $t\!=\!T\!=\!0$ and $u\!=\!z\!=\!0$.
The associated winding numbers are $\fw_t=\fw_u=n_1>0$\,,
We obtain real-time blow-up stars of $n_1$ blow-up branches, alternating with $n_1$ blow-down branches.
In particular, a single-loop Masuda detour near homogeneity generates discrepancies, even for rational spectral quotients, in case $n_1>1$.
The only nonresonant case which produces single-loop, discrepancy-free Masuda detours, along nontrivial leafs, is the semisimple degeneracy $n_1=n_2=1$; i.e.\ $\lambda=1,\ a=0$.

Even \emph{at resonance} $n_1=1<n_2$\,, however, not all is lost.
The formal normal form then features resonant monomials $(0,u^{n_2k})$, for integer $k\geq1$.
Even though such resonant monomials are absent in \eqref{ODEzG1}, the formal normal form might fail to diagonalize \eqref{ODEuG1}, \eqref{ODEzG1} up to any finite order in $u,z$, producing resonant monomials instead.
The closure question of single Masuda detours $\gamma^t,\,\gamma^{uz}$ therefore requires, first, an analysis of that formal normal form.
In case Masuda closure succeeds, on that formal level, we can conclude that the discrepancy \eqref{discrepancyuz} of $\gamma^{uz}$ will be small beyond any upper estimates of finite order in $u,z$.
For quantitative estimates leading to exponentially small discrepancies, see \cite{Neishtadt, ArnoldEnc}.
At present, we do not pursue this intriguing resonance problem any further.
\end{description}

\item[Case 2: The Siegel domain $a>1$.]\hfill\\[2mm]
We obtain a saddle $u\!=\!z\!=\!0$ with $\lambda_1=-1<0<\lambda_2=a-1$ and arbitrary spectral quotients $\lambda=1/(1-a)<0$.
By \eqref{ODEuG1}, \eqref{ODEzG1}, the holonomy $h(u)$ of $u$ is linear, locally near the saddle and over loops $\gamma^z$ around the invariant axis $\{z\!=\!0\}$.
Therefore the dynamics can be analytically orbit diagonalized, locally, by theorem \ref{thmsadlin}.
In particular, neither nonresonance, nor Siegel or Bryuno small divisor conditions \eqref{Dio}, \eqref{Bryuno} are required.
For rational spectral quotients $\lambda=-n_1/n_2$\,, the leafs $\fL$ close upon $n_1\geq1$ iterations of simple loops $\gamma^t$, because $h^{n_1}=\mathrm{Id}$.
This includes the Fujita equation \eqref{PDEw}, i.e.\ $a=2,\ \mathbf{b}=1$ with $n_2=-1=-n_1$\,.
The original Masuda detour on the trivial leaf then closes upon small antisymmetric perturbations of initial conditions, in the Galerkin caricature.
The perturbed solutions, however, do not blow up at $u\!=\!z\!=\!0$, locally.
Therefore they do not qualify as local blow-up loops.
For irrational loops, quasiperiodic holonomy describes how leafs keep closing with better and better approximation, but not quite, upon higher and higher iterations of $\gamma^t$ and $h$.
See remark \ref{remirr} again.
\end{description}

\subsubsection{Nontrivial blow-up equilibria $(u,z)=(0,e_\pm)$}\label{Carsye}

To further explore minimal blow-up loops and blow-up stars, according to section \ref{B}, it remains to linearize \eqref{ODEuG1}, \eqref{ODEzG1} at the two nontrivial blow-up equilibria $z\!=\!e_\pm$ on the compactifying line $[\xi\!:\!\eta\!:\!\zeta]=[\xi\!:\!\eta\!:\!0]\in\CP^1$ at infinity.
Explicitly
\begin{align}
\label{sepm}
e_\pm&=\pm 2\sqrt{1-1/a}\,; \\
\label{l1epm}
\lambda_1&=-a\,, \\
\label{l2epm}
\lambda_2&=2(1-a)\,, \\
\label{lepm}
\lambda&=\tfrac{1}{2}\,a/(a-1)\,, \quad \mathrm{i.e.} \quad a=\,\lambda/(\lambda-\tfrac{1}{2})\,.
\end{align}
Note the \emph{non-semisimple} degeneracy $\lambda\!=\!1$ at the Masuda choice $a\!=\!2$.
Again, the spectral quotients $\lambda$ cover all options for blow-up loops and blow-up stars with polynomial degree $m\!=\!2$ addressed in section \ref{B}.
We distinguish the three cases 
$a<0,\ 0<a<1,$ and $a>1$.

\begin{description}
\item[Case 1: The Poincaré domain $a<0$.]\hfill\\[2mm]
The real unstable nodes $z\!=\!e_\pm$ realize any spectral quotient $0<\lambda<\tfrac{1}{2}$\,.
In backward real time, we obtain analytic flow linearization at nonresonance, and the discussion follows case 1 of section \ref{Carsy0}, verbatim.

\item[Case 2: The Poincaré domain $a>1$.]\hfill\\[2mm]
The real stable nodes $z\!=\!e_\pm$ realize any spectral quotient $\lambda>\tfrac{1}{2}$\,.
Case 1 of section \ref{Carsy0} applies to spectral quotients $\lambda>\tfrac{1}{2}$ at the stable nodes $z\!=\!e_\pm$\,, as well, with one caveat.
The degeneracy $\lambda\!=\!1$ at the Masuda choice $a\!=\!2$ of the original Fujita PDE \eqref{PDEw}, with equal eigenvalues $\lambda_1=\lambda_2=-2$, is not semisimple.
In fact, a term $-e_\pm\neq0$ appears in the lower left of the lower triangular linearizing Jordan block; see \eqref{ODEzG1}.
By remark \ref{remnonss}, \eqref{Jordan} we therefore have to exclude $\lambda\!=\!1$.
In our Galerkin caricature, this is a source of discrepancies for Masuda detours which lead to blow-up at finite size ratios $y/x=e_\pm=\pm\sqrt{2}$.

\item[Case 3: The Siegel domain $0<a<1$.]\hfill\\[2mm]
We obtain purely imaginary saddles $z\!=\!e_\pm$ with arbitrary real spectral quotients $\lambda={\tfrac{1}{2}a/(1-a)}<0$.
Due to absence of a trivial invariant vertical axis, however, the holonomy $h(u)$ of $u$ may not remain linear, after analytic trivialization of the vertically tangent stable manifold. 
Rational $\lambda$, in particular, are resonant and cannot be linearized in general.
Therefore, the linear branching analysis of closed loops in section \ref{RatFol} does not apply.
For irrational $\lambda$, local analytic flow linearization can be achieved under classical Siegel or under optimal Bryuno small divisor conditions \eqref{Dio}, \eqref{Bryuno}.
Quasiperiodic holonomy then describes how leafs only keep closing ``almost'': with better and better approximation, but not quite, upon higher and higher iterations of $\gamma^t$ and $h$.
See remark \ref{remirr}.

Suppose an irrational saddle spectral quotient $\lambda<0$ violates the small divisor condition \eqref{Bryuno}.
Although fully analytic diagonalization may fail, mere absence of resonant terms still estimates deviations from smooth diagonalization to be small beyond any finite order.
For exponential smallness estimates in the analytic case see \cite{Neishtadt, ArnoldEnc}.
This identifies two sources of failing discrepancy closure: the quasiperiodicity induced by $0>\lambda\not\in\Q$, and further perturbations which are exponentially small in the size of the almost-closing, nearly quasiperiodic loops.
\end{description}

\subsection{Spatially antisymmetric modifications}\label{Caras}
Consider the Galerkin caricature \eqref{ODExG2}-\eqref{ODEwG2} for $\mathbf{b}(\bxi)=1-2\cos (2\bxi)+b_1\cos \bxi+b_3\cos(3\bxi)$, next.
The coefficients $b_1,b_3$ generate antisymmetric modifications of the degenerate spatially symmetric case $a=2+b_2=0$; see case 1.2 of section \ref{Carsy0}.
Replacing $\bxi$ by $\pi-\bxi$, we may assume $b_1>0$ without loss of generality.
Let $\beta:=b_3/b_1$\,.
We also assume 
\begin{equation}
\label{bas}
\beta\,\not\in\,\{-3,\, 1,\, 1{+}1/b_1^2\}\,.
\end{equation}
Our assumption $a\!=\!0$ makes the $y$-axis $\{x\!=\!0\}$ in \eqref{ODExG2} invariant.
Since $b_1\neq0$ prevents a blow-up equilibrium $u\!=\!z\!=\!0$, it will now be sufficient to consider blow-up coordinates $v\!=\!1/y,\ w\!=\!x/y$\,; see \eqref{ODEvG2}, \eqref{ODEwG2}.
In particular, both axes in blow-up coordinates $v,w$ are invariant.

\subsubsection{The trivial blow-up equilibrium $v=w=0$}\label{Caras0}

On the trivial spatially antisymmetric leaf $\fL\,{:}\ w=x/y=0$, the dynamics of $v\!=\!1/y$ is equivalent to the Riccati dynamics of $v$ on the $v$-axis.
Indeed, $3b_1+b_3\neq0$, and \eqref{ODEvG2} recovers the quadratic blow-up of Riccati section \ref{ODERic}.
It remains to discuss the nontrivial leafs.

The semisimple spectrum along the invariant axes is 
\begin{align}
\label{a0l1}
b_1\lambda_1&=-\tfrac{1}{4}(\beta+3)\,, \\
\label{a0l2}
b_1\lambda_2&=-\tfrac{1}{4}(\beta-1)\,, \\
\label{a0l}
\lambda&=(\beta+3)/(\beta-1)\,, \quad \mathrm{i.e.} \quad\beta =\,(\lambda+3)/(\lambda-1)\,.
\end{align}
Again, the free parameter $\beta$ in the spectral quotient $\lambda$ covers all options for minimal blow-up loops and blow-up stars with polynomial degree $m\!=\!2$ studied in section \ref{B}.
With $b_1>0$, we distinguish the three cases $\beta<-3,\ -3<\beta<1$, and $\beta>1$.

\begin{description}
\item[Case 1: The Poincaré domains $\beta<-3$ \textbf{and} $\beta>1$.]\hfill\\[2mm]
In the Poincaré domains $\beta<-3$ and $\beta>1$, we encounter a real unstable node and a real stable node, respectively, with arbitrary spectral quotients $0<\lambda\neq1$.
Note how this excludes $1{:}1$ nonresonant semisimple double eigenvalues.
The only resonances are other integer ratios $\lambda,\,1/\lambda\in\N\setminus\{1\}$.
Therefore the Poincaré linearization theorem \ref{thmPoi} applies, with all the consequences for minimal blow-up loops, blow-up stars, and quasiperiodicity already discussed in case 1 of section \ref{Carsy0}.

\item[Case 2: The Siegel domain $-3<\beta<1$.]\hfill\\[2mm]
We obtain a saddle $v\!=\!w\!=\!0$ with $\lambda_1<0<\lambda_2$ and arbitrary spectral quotients $\lambda=(\beta+3)/(\beta-1)<0$.
The associated holonomy $h(v)$ of $v$ is not linear, \emph{a priori}, locally near the saddle and over loops $\gamma^w$ around the invariant axis $\{v\!=\!0\}$; see \eqref{ODEvG1}.
Analytic diagonalization may therefore fail, at resonant rational spectral quotients $\lambda$.
Irrational $\lambda$, although nonresonant with smooth diagonalization beyond finite order, require additional small divisor properties \eqref{Dio} or \eqref{Bryuno} to qualify for analytic diagonalization.
See section \ref{SieLin} and case 3 of section \ref{Carsye} for further discussion.
\end{description}

\subsubsection{Nontrivial blow-up equilibria $(v,w)=(0,e_\pm)$}\label{Carase}

To further explore minimal blow-up loops and blow-up stars, via section \ref{B}, it remains to consider the two nontrivial blow-up equilibria $w\!=\!e_\pm\neq0$ on the compactifying line $[\xi\!:\!\eta\!:\!\zeta]=[\xi\!:\!\eta\!:\!0]\in\CP^1$ at infinity.
Explicitly, \eqref{ODEwG2} and assumption \eqref{bas} imply
\begin{equation}
\label{w2}
b_1w^2-w+\tfrac{1}{4}(b_3-b_1)=0\,,\quad \textrm{with discriminant } d=1+b_1^2(1-\beta)\neq0\,;
\end{equation}
For the equilibria and their linearizations, we obtain
\begin{align}
\label{sasepm}
e_\pm&=\tfrac{1}{2}(1\pm\sqrt{d})/b_1\,; \\
\label{l1asepm}
\lambda_1&=-e_\pm-b_1\,, \\
\label{l2asepm}
\lambda_2&=-e_\pm-b_1(1-\beta)/2\,
\\
\label{lasepm}
\lambda&=\left(e_\pm+b_1\right)\,/\,\left(e_\pm+b_1(1-\beta)/2\right)
\end{align}
For positive discriminant $d>0$, the interested reader may explore the resulting real equilibria and real eigenvalues along the lines of our previous discussions.
We address negative discriminants $d<0$ instead, i.e.\ $\beta>1+b_1^{-2}$.
This leads to equilibria $w\!=\!e_\pm$\,, eigenvalues $\lambda_1,\lambda_2$ and spectral quotients $\lambda$, which are all nonreal complex.
In this case we are in the Poincaré domain by proposition \ref{propPoires}(i), automatically, hence nonresonant by proposition \ref{propPoires}(ii), and analytically diagonalizable by theorem \ref{thmPoi}.
However, diagonalization also implies that the only way to obtain any blow-up loops $(\gamma^t,\,\gamma^{vw})$ for $\lambda\in\C\setminus\R$ is within the eigenspaces of the diagonalization.
The eigenspace of $\lambda_2$ is the nonlinearly invariant $w$-axis $v\!=\!0$, which resides in the compactification $[w\,{:}\,1\,{:}\,0]\in\CP^1$ at infinity.
Therefore the only eligible eigenspace for blow-up loops is the complex eigenspace of $\lambda_1=-e_\pm-b_1$ with complex linearly dependent, nonvanishing components $(v,w)$.
In diagonalized coordinates, note the linear holonomy \eqref{linhol} in terms of the spectral quotient $\lambda=\lambda_1/\lambda_2$\,.
Upon a complex rotation of time, which rotates $\lambda_1$ to $-|\lambda_1|<0$, our previous comments on the trivial stable leafs of case 0 in section \ref{Carsy0} apply, albeit with nonlinearly related loops $\gamma^v,\,\gamma^w$ in the original variables.

\section{Conclusions on real Galerkin blow-up}\label{Con}
We summarize the results of section \ref{Fuji} with a focus on two aspects: spatial symmetry of $\mathbf{w}(t,\bxi)$ in $\bxi$, and real-time blow-up of real $\mathbf{w}$.
We compare the blow-up effects of symmetric spatial inhomogeneities, versus antisymmetric modifications, caused by the spatially inhomogeneous coefficient
\begin{equation}
\label{bxi}
\mathbf{b}(\bxi)=1+b_2\cos(2\bxi)+b_1\cos\bxi+b_3\cos(3\bxi)
\end{equation}
with real Fourier coefficients $b_k$\,; compare \eqref{bj}.
In particular, real-valued subspaces are invariant, in real time.
For the crucial role of rational versus irrational spectral quotients $\lambda$ of analytic ``linearizations at infinity'', we urge the reader to recall remark \ref{remirr}.

We have studied the effects in the Galerkin caricature $\mathbf{w}=x+y\cos\bxi$ of the quadratic Fujita PDE $\mathbf{w}_t=\mathbf{w}_{\bxi\bxi}+\mathbf{b}(\bxi)\mathbf{w}^2$; see \eqref{PDEb}-\eqref{ODEa} and \eqref{ODExG}, \eqref{ODEyG}.
We call $\mathbf{w}$ \emph{symmetric} and \emph{antisymmetric}, respectively, if $\mathbf{w}(\pi-\bxi)=\pm\mathbf{w}(\bxi)$; similarly for Fourier components of $\mathbf{b}(\bxi)$.
In our Galerkin caricature, symmetric $\mathbf{w}$ correspond to the $x$-axis $y\!=\!0$.
Somewhat exaggerating mere ``symmetry'', the $x$-axis only consists of spatially homogeneous $\mathbf{w}=x$.
The $y$-axis represents antisymmetric $\mathbf{w}$, of spatial average $x\!=\!0$.
Somewhat exaggerating ``antisymmetry'', as well, zero average $x$ implies antisymmetry, in our crude caricature.

We have explored asymptotics towards blow-up time $t\nearrow T=0$, and along complex time loops $\gamma^t $ around $T\!=\!0$.
Necessarily, this required analysis of ODE solutions $(x,y)\in\C^2$, and their associated complex foliations, by maximal local solution leafs $\fL$.
Blow-up was studied via projective compactification in rescaled time; see \eqref{ODEuG}-\eqref{ODEwG}.
The blow-up singularities of complex foliations, of complex codimension two, arise as equilibria ``at infinity'', i.e.\ for projective coordinates $u\!=\!1/x\!=\!0$, or $v\!=\!1/y\!=\!0$.
The complementary coordinates $z\!=\!y/x$, or $w\!=\!x/y\!=\!1/z$, keep track of the relative proportions of the spatially homogeneous, i.e.\ symmetric, $\mathbf{w}$-components $x,u$ versus the spatially antisymmetric components $y,v$.
Our analysis was based on the spectral quotient $\lambda:=\lambda_1/\lambda_2$ of eigenvalues of the linearization at such blow-up equilibria.
Our main interest were \emph{minimal blow-up loops} $(\gamma^t,\gamma^{uz})$ and $(\gamma^t,\gamma^{vw})$, respectively, around blow-up equilibria, together with their associated \emph{real-time blow-up stars} of real or complex continuations of solutions.
See definitions \ref{defloop} and \ref{defstar}.

In blow-up coordinates, the spatially homogeneous $x$-axis and the zero average antisymmetric $y$-axis are represented by the axes $z\!=\!y/x\!=\!0$ of $u\!=\!1/x$, and $w\!=\!x/y\!=\!0$ of $v\!=\!1/y$, respectively.
Symmetric inhomogeneities $\mathbf{b}(\bxi)$, or actually $b_1=0$ alone, imply invariance of the symmetric $x$-axis and $u$-axis.
Note the coefficient $a=2+b_2$ in \eqref{ODExG1}, \eqref{ODEwG1} which breaks antisymmetric invariance of the $y$- and $v$-axes in the ODEs \eqref{ODExG1}-\eqref{ODEwG1} of symmetric $\mathbf{b}(\bxi)$.
Curiously, invariance of these axes does hold for $a=b_2+2=0$ in \eqref{ODExG2}-\eqref{ODEwG2}, even though $\mathbf{b}(\bxi)$ itself is not antisymmetric and, therefore, the ODEs are not equivariant under the action of antisymmetry $\kappa\mathbf{w}(t,\pi-\bxi):=-\mathbf{w}(t,\bxi)$, alias the reflection $\kappa(x,y):=(-x,y)$.

Unrelated to symmetry, invariance of the $z$-axis at infinity, $u\!=\!0$, and of the $w$-axis at infinity, $v\!=\!0$, is caused by the Euler multipliers of our time regularizations $t_1,t_2$ at blow-up.
They generate invisible dynamics of $z,w$ within the projectively compactifying invariant Riemann sphere at infinity.
By  \eqref{ODEzG}, \eqref{ODEwG}, we obtain scalar cubic ODEs for $z,w$.
We refer to \cite{FiedlerShilnikov, FiedlerYamaguti} for a detailed discussion of the global dynamics for generic polynomial ODEs on the Riemann sphere; see also figure \ref{fig14trees}.

The \emph{symmetric case} $\mathbf{b}(\bxi)=1+(a-2)\cos (2\bxi)$ of \eqref{bxi} of \eqref{ODExG1}-\eqref{ODEwG1} has been studied in section \ref{Carsy}.
Minimal blow-up loops and real-time blow-up stars at the \emph{trivial equilibrium} $u\!=\!z\!=\!0$ address small relative perturbations $z\!=\!y/x$ of the purely homogeneous Riccati case $y\!=\!0\!=\!z$; see section \ref{Carsy0}.

Small relative perturbations $z\!=\!y/x$ of homogeneity already feature rescaled profiles $\mathbf{w}(t,\bxi)$ with extrema at the Neumann boundary, but still dominated by the large spatial average $x\!=\!1/u$.
See \eqref{s0l1}-\eqref{s0l} for the semisimple real eigenvalues and their spectral quotient $\lambda=1/(1-a)$.
Nonresonant rational spectral quotients $0<\lambda, 1/\lambda\in\Q\setminus\N$ did establish discrepancy-free minimal blow-up loops, and real-time blow-up stars, in the Poincaré domain $a<1$, where the trivial blow-up equilibrium is a stable node, in real time.
See definitions \ref{defloop}, \ref{defstar}, sections \ref{PoiLin}, \ref{PoiMas}, and case 1.2 in section \ref{Carsy0}.
Mostly, minimal blow-up loops required multiple winding of $\gamma^t$ around the blow-up time $T\!=\!0$.
A single winding was sufficient only for the semisimple $1{:}1$ resonance $\lambda\!=\!1$ at $a\!=\!0$.
See remark \ref{remnonss}.
The technical tool was analytic flow linearization in the Poincaré domain.
Irrational spectral quotients led to quasiperiodic holonomy of the foliations, at blow-up.
See section \ref{Carsy0}, case 1.1 and remark \ref{remirr}.
This prevents discrepancy-free blow-up loops.
The resulting non-closure of leafs, however, is alleviated by quasiperiodic discrepancies in $u,z$ tending to zero, for certain sufficiently high iterates of the complex time-loop $\gamma^t$ around blow-up time $T\!=\!0$.

Due to linear holonomy, the Siegel domain $a>1$ of saddle spectral quotients $\lambda<0$ is analytically orbit linearizable at $u\!=\!z\!=\!0$.
Rational $\lambda$ lead to closed nontrivial loops, near blow-up.
This includes our caricature of the Fujita case \eqref{PDEw} with $a\!=\!2,\ \mathbf{b}\!=\!1$.
Irrational spectral quotients $0>\lambda\not\in\Q$ lead to quasiperiodic holonomy again, but do require small divisor conditions like \eqref{Dio}, \eqref{Bryuno} for the requisite analytic diagonalization of the blow-up system.
See section \ref{Carsy0}, case 1.2 again.

\emph{Nontrivial blow-up equilibria} $u\!=\!0, \ z\!=\!e_\pm\neq0$ for symmetric $\mathbf{b}(\bxi)$ describe asymptotic blow-up profiles at fixed ratios $z\!=\!y/x$ between the antisymmetric perturbation component $y$ and the spatial average $x$.
See section \ref{Carsye}.
This addresses asymptotic real blow-up profiles with leading blow-up located at the left or right Neumann boundary of the spatial interval $0<\bxi<\pi$.
See \eqref{sepm}-\eqref{lepm} for $e_\pm$\,, real eigenvalues, and their real spectral quotient $\lambda=\tfrac{1}{2}a/(a-1)$.
In the Poincaré domains $a<0$ and for $a>1$, the equilibria are also real.
See case 2, $a>1$, of section \ref{Carsye} for minimal blow-up loops and real-time blow-up stars at nonresonant spectral quotients $0<\lambda\neq1$.
Although the $1{:}1$ resonant original Fujita case $a\!=\!2,\ \mathbf{b}\!=\!1$ with $\lambda_1=\lambda_2=-2$ qualifies for linearization, however, the eigenvalues are not semisimple.
Therefore discrepancies of Masuda detours arise, as described in remark \ref{remnonss}.
For blow-down variants see case 1 of $a<0$.
Sign changing real blow-down and blow-up asymptotics $\mathbf{w}$ occur for $|e_\pm|>1$, i.e.\ in the Poincaré domain $a<0$ and for $1<a<4/3$.

All blow-up equilibria have been real, so far, and accessible by real blow-up solutions.
The symmetric Siegel domain $0<a<1$, however, first features conjugate complex blow-up equilibria with nonreal ratios $y/x\!=\!z\!=\!e_\pm$\,, here purely imaginary.
By invariance of the real subspace $(x,y)\in\R^2$, these cannot be reached by any blow-up of real initial conditions in real time.
From a real (viz. limited) point of view, therefore, nontrivial nonreal blow-up equilibria are of mathematical interest, ``only''.
See section \ref{Carsye}, case 3, anyway.

The \emph{antisymmetric modifications} $\mathbf{b}(\bxi)=1-2\cos (2\bxi)+b_1(\cos\bxi+\beta\cos(3\bxi))$ with $\beta:=b_3/b_1$\,, alias ODEs \eqref{ODExG2}-\eqref{ODEwG2}, have been studied in section \ref{Caras}, choosing $b_1>0$.
Somewhat atypically, zero average $x\!=\!0$, which also represents spatial antisymmetry $\mathbf{w}(\pi-\bxi)=-\mathbf{w}(\bxi)$ in our Galerkin caricature, is preserved by our choice of $a=b_2+2=0$.
See \eqref{ODExG2}.
The local dynamics around the \emph{trivial blow-up equilibrium} $v\!=\!w\!=\!0$ represents almost antisymmetric dynamics.
The stable/unstable blow-up node of the Poincaré domains $\beta\in\R\setminus[-3,1]$, as well as the resulting Siegel saddle for $\beta\in(-3,1)$ feature real spectral quotients $\lambda=(\beta+3)/(\beta-1)$ and have been discussed as before; see section \ref{Caras0}.

For negative discriminant $d=1+b_1^2(1-\beta)$, we again encounter \emph{nontrivial blow-up equilibria} $e_\pm=\tfrac{1}{2}(1\pm\sqrt{d})/b_1$ which are nonreal conjugate complex; see \eqref{sasepm}.
Nonreal asymptotic proportions $w\!=\!x/y\!=\!e_\pm$ cannot be reached from real initial conditions, in real time.
Therefore we just refer back to section \ref{Carase}, for mathematical details.
We have left the case of positive discriminant $d$, which follows the previous scheme, to the reader.

In conclusion, we have seen how real-time blow-up of $(x,y)\in\R^2$ depends, explicitly but in quantitatively subtle ways, on the actual parameters and symmetry aspects of the spatially inhomogeneous quadratic coefficient $\mathbf{b}(\bxi)$ in \eqref{bxi}.
And analogously for blow-down.
We have encountered open domains of blow-up towards nodes, in Poincaré domains, versus separatrix blow-up towards saddles, in Siegel domains, as well as the existence of minimal complex-time blow-up loops and real-time blow-up stars.
It will be exciting to explore this rich phenomenology, in complex time, for the higher-dimensional Galerkin ODE context, or even for the full PDE.

\section{Galerkin embedding of planar ODEs in parabolic PDEs}\label{PDElift}

In this section we sketch a complementary approach to replace our crude Galerkin caricatures of PDE blow-up.
Instead of exploring inadequate 2-mode Galerkin caricatures of given, prescribed PDEs, theorem \ref{thmLift} \emph{constructs parabolic PDEs} with any exactly prescribed, given ODE dynamics in an invariant, complex two-dimensional Galerkin subspace.
Any ODE blow-up then embeds into the resulting PDE dynamics, verbatim.

Such an approach is reminiscent of backwards error analysis in numerics, where we do not ask for the numerical error with respect to an exact solution of a given system.
Instead, we ask how large a  perturbation of the original system is required, to render the numerical solution exact in the perturbed system.
See for example \cite{Backward}, and for complex methods \cite{Matthiesback, OliverWulff, WulffOliver}.

Specifically, consider any given complex ODE flow 
\begin{equation}
\label{ODEfgid}
\begin{aligned}
   \dot x&=-x+f(x,y)\,,\\ 
   \dot y&=-y+g(x,y)
\end{aligned}
\end{equation}
on $\mathbf{x}=(x,y)\in\C^2$.
Here we have just slightly shifted the notation of nonlinearities $f,g$, in our original ODE \eqref{ODEfg}.
We aim at a linear embedding of \eqref{ODEfgid} into a PDE
\begin{equation}
\label{PDEf}
\mathbf{w}_t=\mathbf{w}_{\bxi\bxi}+\mathbf{f}(\bxi,\mathbf{w},\mathbf{w}_{\bxi})
\end{equation}
under periodic boundary conditions $\bxi\in\S^1=\R/2\pi\Z$\,.
More precisely, we require the complex Galerkin (alias Fourier) subspace $X_1:=\langle \sin\bxi,\,\cos\bxi\rangle_\C\cong\C^2$ of the phase space $X$ of \eqref{PDEf} to be invariant under \eqref{PDEf}, with the prescribed semiflow \eqref{ODEfgid} on it.
In complex coordinates $\mathbf{w}(t,\bxi)=x(t)\cos\bxi+y(t)\sin\bxi$\,, this is equivalent to the condition
\begin{equation}
\label{ODElift}
   \mathbf{f}(\bxi,x\cos\bxi+y\sin\bxi,-x\sin\bxi+y\cos\bxi)\ :=\ f(x,y)\cos\bxi+g(x,y)\sin\bxi
\end{equation}
for the function $\mathbf{f}{:}\ \S^1\times\C^2\rightarrow\C$, on all arguments $(\bxi,x,y)\in\S^1\times\C^2$.
Let $\mathbf{R}(\bxi)\mathbf{x}$ denote the standard rotation of a vector $\mathbf{x}=(x,y)\in\R^2$ or $\C^2$ by an angle $\bxi\in\S^1$.
Then \eqref{ODElift} abbreviates to
\begin{equation}
\label{ODEliftR}
   \mathbf{f}(\bxi,\mathbf{x})\ :=\ 
   f(\mathbf{R}(\bxi)\mathbf{x})\cos\bxi+g(\mathbf{R}(\bxi)\mathbf{x})\sin\bxi\,.
\end{equation}
Note how the map $(f,g)\mapsto\mathbf{f}$ defined by \eqref{ODEliftR} preserves complex analyticity in $\mathbf{x}$.
It maps entire functions $(f,g)$ to entire functions $\mathbf{f}$, and polynomials to polynomials with preserved polynomial degrees.
Conversely, any nonlinearity $\mathbf{f}{:}\ \S^1\times\C^2\rightarrow\C$ in \eqref{PDEf}, for which the Galerkin subspace $X_1$ is flow-invariant, defines the vector field $(f,g)$ of the ODE \eqref{ODEfgid} on it. This proves the following theorem.

\begin{thm}\label{thmLift}
Define the polynomial nonlinearity $\mathbf{f}$ by \eqref{ODEliftR}.\\
This constructs PDE examples \eqref{PDEf} on the unit circle $\bxi\in\S^1$ \emph{for each and every blow-up result} on complex planar polynomial ODEs \eqref{ODEfgid}, viz.\  \eqref{ODEfg}, of the present paper.
The ODEs are given by exact, error-free Galerkin projection onto the PDE-invariant Galerkin (alias Fourier) subspace $X_1=\langle \sin\bxi,\,\cos\bxi\rangle_\C$\,.
\end{thm}

In the real case, this idea goes back to Heidelberg days; see theorem 2 in \cite{FSandstede}.
For the milder requirement of \emph{jet embeddings} in flow-invariant center manifolds, i.e.\ prescribed finite reduced Taylor expansions, see \cite{FPolacik}.
For delay equations, the idea had been pioneered by \cite{HaleEmbed}.
For a broader survey of real ODE flow embeddings in PDEs see section 7 in \cite{PolacikEmbed}. 

At least in the real parabolic case, we have to caution the reader that the above embedding is dynamically unstable.
This is due to exponential dichotomies associated with real PDEs \eqref{PDEf}; see \cite{RochaDichotomy} and, for delay equations, \cite{MPDichotomy}.

\section{Comments}\label{Com}

We comment on some previous literature and discuss further prospects.
Section \ref{Ent} recalls a few results which force \emph{any} nonstationary solution of polynomial and complex entire ODEs to blow up in finite time. 
Section \ref{Global} collects some global results in low complex dimensions 1 and 2, including nongeneric variants of section \ref{Scalar}.
Already the global classification of quadratic ODEs in complex dimension 2 has met very substantial obstacles.
For higher dimensions, we are essentially left to the local linearization devices of section \ref{B}; see section \ref{Dim}.
Real-time heteroclinic ODE orbits $\mathbf{x}\,{:}\ \mathbf{e}_-\leadsto \mathbf{e}_+$\,, however, are often accompanied by blow-up in imaginary time, even in the PDE settings of section \ref{PDE}.
The final section \ref{App} advocates a revival of the fascinating subject of ODE and PDE dynamics in complex phase spaces and complex time, towards applications.
We include a 1,000\;\euro\ prize question.
The question asks for nonconstant homoclinic orbits $\mathbf{x}\,{:}\ \mathbf{e}\leadsto \mathbf{e}$ which are entire, i.e.\ globally analytic in complex time.
This is motivated by spurious splitting effects attached to discretization of real ODEs and PDEs in real time.

\subsection{Blow-up versus entire solutions}\label{Ent}

For polynomial nonlinearities $\mathbf{F}$ of ODE \eqref{ODEF} in complex dimensions 1 and 2, we have encountered meromorphic branching versus essential singularities; see remark \ref{remirr} and section \ref{B}.
Already Rellich  considered the second order pendulum \eqref{pen} with general nonlinear entire force laws $g$ \cite{Rellich}.
He proved the absence of \emph{any} entire, i.e.\ globally holomorphic, nonequilibrium solutions $t\mapsto \mathbf{w}(t)$.

For entire unstable manifolds in complex dimensions 2 and higher, see \cite{Ushiki1, Ushiki2}.
As a corollary, on any saddle-saddle heteroclinic orbits via one-dimensional unstable \emph{and} one-dimensional stable manifolds of entire ODEs \eqref{ODEF}, blow-up has to occur.

For entire, nonpolynomial nonlinearities, the description of separatrix basin boundaries of finite un/stable nodes or centers is difficult.
See \cite{MuciRiem, MuciEss} and, most recently, some examples, descriptions, and detailed references in \cite{LebEss}.

\subsection{Global aspects}\label{Global}

Blow-up in real time, for \emph{generic polynomial or rational ODEs} in scalar $x\in\C$, has been discussed in \cite{FiedlerShilnikov, FiedlerYamaguti} from the viewpoint of global $C^0$ orbit equivalence, as defined in \ref{defequi}.
The polynomial case $\dot x=P(x)$ was reviewed in section \ref{Scalar} above and reappeared, e.g., in separatrices at blow-up.
See sections \ref{C2m}, \ref{Hamm}, \ref{Con}.
Blow-up and blow-down separatrices where given by stable and unstable leafs of separatrix saddles in the Poincaré compactification at $x=\infty\in\CP^1=\widehat{\C}$.
Regularized generic rational vector fields 
\begin{equation}
\label{rat}
\dot x=P(x)/Q(x)
\end{equation}
are $C^0$ orbit equivalent to \emph{all} Morse flows on the Riemann 2-sphere $x\in\widehat{\C}$; see \eqref{hybrid}. 

The general rational case \eqref{rat} is not necessarily generic, and therefore much more complicated.
See \cite{DiPol} for a complete characterization of the global dynamics in the polynomial case $f\!=\!P$, and \cite{DiRat} for the general rational case $f\!=\!P/Q$.
The characterization involves nine necessary properties of global flows on $\widehat{\C}$.

Violation of hyperbolic nondegeneracy, for example, already happens at Lyapunov centers in families $f=f(\lambda,\cdot)$ with a scalar real parameter $\lambda$.
Global descriptions of their unfoldings are a wide open field, in the complex context.
See \cite{Rousseaud=4} for interesting sample bifurcations, including a 2-parameter unfolding.
For the polynomial degrees of the numerator and denominator polynomials $P$ and $Q$ they chose 4 and 2, respectively, to satisfy the difference constraint by Euler characteristic 2 and avoid equilibria and poles at infinity.
For higher degrees, even the global transitions between generic phase portraits, induced by Lyapunov centers in generic 1-parameter families of scalar ODEs, seem not to be completely understood.
The difficulties include the global real-time consequences of a mere complex rotation in time for, say, the connection graphs enumerated in figure \ref{fig14trees}.

An obvious real Euler multiplier $\rho=|Q(x)|^2$, which preserves real and imaginary time directions, converts \eqref{rat} to the hybrid scalar ODE 
\begin{equation}
\label{hybrid}
\dot x=P(x)\overline{Q(x)} 
\end{equation}
with a holomorphic factor $P$ and an \emph{antiholomorphic factor} $\overline{Q}$.
Consider the purely antiholomorphic, reciprocal case $P\!=\!1$ first.
Then the resulting antiholomorphic flow of $x\in\C\cong\R^2$ is, simultaneously, a gradient flow \emph{and} Hamiltonian with respect to the real and imaginary parts of the integrated denominator $\int^x Q\,dx$.
As in the polynomial case $Q\!=\!1$, a complete global $C^0$ classification of the generic cases ensues.
This leads to a description of compactified connection graphs by noncrossing circle chords, and an associated explicit combinatorial case count, much in the spirit of \eqref{countpol1}.
See section 2.3 in \cite{FiedlerYamaguti} for details.

Classification of rational ODEs \eqref{rat}, from the hybrid point of view, relies on the factorization  \eqref{hybrid}, after Euler regularization.
Without such additional assumptions, polynomial hybrid ODEs $\dot x = P(x,\bar{x})$ are just another notation for arbitrary polynomial systems of degree $m$ in $\R^2$.
The unsolved Hilbert problem XVI concerning a bound $L=L(m)<\infty$ on the number of limit cycle periodic orbits arises in that class.
See \cite{Ilya}, section IV.24 for a survey with references.
In particular, any global classification of \eqref{ODEF} in two complex dimensions $\mathbf{x}=(x_1,x_2)\in\C^2$ is notoriously difficult.

For a $C^0$ classification in the real quadratic case $m\!=\!2$, and under a simplifying reduction which collapses limit cycles nested around the same single equilibrium, see \cite{Llibre1}.
They identify 44 distinct generic phase portraits, in the Poincaré compactification. 
For generic degeneracies of real codimension one, that count increases up to 211, under pending realizability of 7 cases.
A main difficulty is the classification of saddle-saddle heteroclinic orbits and homoclinic loops, which are not accessible to purely algebraic methods.
For a purely algebraic classification of singularities, only, see also \cite{Llibre2}.

\subsection{Higher dimensions}\label{Dim}

Classical studies by Wittich \cite{Wittich1, Wittich2} extended \cite{Rellich} to large classes of scalar $N$-th order versions of ODE \eqref{ODEF} which are nonpolynomial entire in $x\!=\!x_1$\,, and polynomial in $t$ and the derivatives $\dot x,\ldots, x^{(N)}$.

The projective compactification of sections  \ref{Fol} and \ref{B} applies to $\C^N\subset\CP^N$, in principle, via the higher-dimensional variants of analytic diagonalization in the nonresonant Poincaré and Siegel domains.
\emph{Resonance} for higher-dimensional spectral and integer coefficient vectors $\boldsymbol{\lambda}, \boldsymbol{\alpha}$ is defined by \eqref{res}, verbatim. 
For simplicity of presentation, we refrained from indulging in such generalizations here, say for higher-dimensional Galerkin approximations of the spatially heterogeneous quadratic Fujita equation \eqref{PDEb}. 
More importantly, it is not clear how well direct Galerkin approximations represent blow-up phenomena in complex time, at all.                 

Another approach to blow-up singularities of ODEs \eqref{ODEF} for higher-dimensional complex entire vector fields $\mathbf{F}(\mathbf{x})\in \mathbf{X}=\mathbb{C}^N$, in complex time, starts from real heteroclinic orbits $\mathbf{x}(t)\rightarrow \mathbf{e}_\pm$ in real time $t\rightarrow\pm\infty$.
In symbols: 
\begin{equation}
\label{het2}
\mathbf{x}\,{:}\quad \mathbf{e}_-\leadsto \mathbf{e}_+ \,.
\end{equation}
See also \eqref{het} for this notation.
Assume separately Poincaré nonresonant, real unstable and real stable spectra at the hyperbolic saddle equilibria $\mathbf{e}_-$ and $\mathbf{e}_+$\,, respectively.
Then the maximal complex time extension of the real-time heteroclinic orbit $\mathbf{x}(t)$ experiences finite-time blow-up, in the imaginary time direction.

This result is based on Poincaré linearization, as in section \ref{PoiLin}, applied inside the complex analytic, finite-dimensional, unstable and stable manifolds of the heteroclinic source $\mathbf{e}_-$ and target $\mathbf{e}_+$\,, separately and respectively.
Suppose, indirectly, that an entire heteroclinic orbit \eqref{het} does exist.
Global incompatibility of the coefficients of quasiperiodicity near $\mathbf{e}_\pm$\,, in imaginary time, then leads to a contradiction.
Such an indirect proof, unfortunately, does not reveal the geometry of singularities in complex time.
It does not even prove that singularities are isolated in time.
Nor does it address issues like higher-mode Galerkin approximations of spatio-temporal singularities in PDEs.
See \cite{FiedlerClaudia} for details and further discussion.
See \cite{JaquetteStuke, JaquetteMasuda, Jaquetteqp}, however, for computer assisted proofs of blow-up in certain quadratic PDEs.

\subsection{Partial differential equations}\label{PDE}

Subsequently, we have explored consequences of real heteroclinicity for blow-up in a PDE context; see \cite{FiedlerFila}.
In real time $r$, consider real analytic power series solutions $r\mapsto \mathbf{w}(r,\bxi)$ of a reaction-diffusion equation
\begin{equation}
\label{w6}
\mathbf{w}_r=\mathbf{w}_{\bxi\bxi}+6\mathbf{w}^2-\ell
\end{equation}
with a quadratic nonlinearity of Riccati type. 
Here $0<\bxi<1/2$\,, with Neumann boundary condition, and $\ell>0$ serves as a bifurcation parameter.
Then $\mathbf{w}=\mathbf{w}(r,\bxi)$ extends to complex time $t=r+\mi s$, locally, and $\boldsymbol{\psi}(s,\bxi):=\mathbf{w}(r_0-\mi s)$ satisfies the \emph{nonconservative quadratic Schrödinger equation}
\begin{equation}
\label{psi6}
\mathrm{i}\boldsymbol{\psi}_s=\boldsymbol{\psi}_{\bxi\bxi}+6\boldsymbol{\psi}^2-\ell
\end{equation}
with real initial condition $\boldsymbol{\psi}(0,\bxi)=\mathbf{w}(r_0,\bxi)$.
In this way, complex time extension provides families of analytic Schrödinger solutions $\boldsymbol{\psi}$, from a single heat solution $\mathbf{w}$, and vice versa. 
Moreover, solution families of one type are related by a (semi)flow of the other type; see \eqref{flow}.
The simplest case of spatially homogeneous solutions $\mathbf{w}=\mathbf{w}(t)$ leads back to the Riccati ODEs of section \ref{Ric}.
For some previous literature in such PDE context, see for example \cite{COS,Yanagida,Stukediss,Stukearxiv,JaquetteHet,JaquetteStuke,JaquetteMasuda, Jaquetteqp, Fasondini24}.
For a quick summary and further references see again \cite{FiedlerFila}.

One PDE feature of the real quadratic heat equation \eqref{w6} is the abundance of nonconstant heteroclinic orbits $\mathbf{w}\,{:}\ \mathbf{e}_-\leadsto \mathbf{e}_+$ in real time $r$ and in suitable Banach spaces $\mathbf{w}(t)=\mathbf{w}(t,\cdot)\in X$. See also \eqref{het}, \eqref{het2}.
Here $\mathbf{e}_\pm$ denote equilibria of \eqref{w6} which, up to reversed sign, spatially satisfy the pendulum equation \eqref{pen}, \eqref{penW} discussed above.
\emph{Connecting orbits}, \emph{traveling fronts} (or backs), or \emph{solitons} are other names for $\mathbf{w}$, depending on context.
In fact, any globally bounded, real-valued solution $\mathbf{\mathbf{w}}(r,x), \ r\in\mathbb{R}$ of \eqref{PDEb} is heteroclinic between distinct equilibria $\mathbf{e}_\pm$\,.
This is due to a decreasing Lyapunov functional on real solutions; see \cite{Zelenyak, Matano, Lappicy}.
For further developments in the real case, see also \cite{MatanoLap,brfi88,brfi89,Raugelattr,Galaktionov,firoSFB,firoFusco} and the many references there.

The main result in \cite{FiedlerFila} asserts that real-time real heteroclinic orbits $\mathbf{w}\,{:}\ \mathbf{e}_-\leadsto \mathbf{e}_+$ between hyperbolic equilibria $\mathbf{e}_\pm$ then experience blow-up in imaginary time $\pm s$, i.e.\ for some Schrödinger solution $\boldsymbol{\psi}(s,\bxi):=\mathbf{w}(r_0-\mi s)$ and a suitable choice of $r_0$\,.

Limiting assumptions are Poincaré nonresonant unstable spectrum of the Sturm-Liouville spectrum at the source $\mathbf{e}_-$\,, and linear asymptotic stability of the target, here $\mathbf{e}_+=-\sqrt{\ell/6}$\,.
The stability requirement remedies the absence of any infinite-dimensional versions of Poincaré linearization.
Unstable spectral nonresonance at spatially inhomogeneous source equilibria $\mathbf{e}_-$ is a very delicate issue which limits the result further, to (most) source equilibria of unstable dimension not exceeding 22, or to sufficiently fast unstable manifolds of $\mathbf{e}_-$\,.
See theorems 1.2 and 1.3 in \cite{FiedlerFila}.
As in the ODE paradigm of \cite{FiedlerClaudia}, the precise geometry of the Schrödinger blow-up remains unexplored, save vague statements like $L^\infty$ blow-up.

Beyond the spatial inhomogeneity of the quadratic coefficient $\mathbf{b}=\mathbf{b}(\bxi)$ in \eqref{PDEb}, thin domain analysis  motivates variants of the Fujita PDE \eqref{PDEw} with spatially inhomogeneous diffusion, like
\begin{equation}
\label{PDEa}
\mathbf{w}_t=\mathbf{a}(\bxi)^{-1}(\mathbf{a}(\bxi)\, \mathbf{w}_{\bxi})_{\bxi}+\mathbf{w}^2 \,.
\end{equation}
See for example the survey \cite{Raugelthin}.
Other variants, popular in homogenization, prefer purely variational form and omit the term $\mathbf{a}(\bxi)^{-1}$.
One feature of \eqref{PDEa} and such variants is invariance of the subspace of spatially homogeneous solutions.
In the projective Galerkin caricatures of proposition \ref{propuvw}, this implies $g(x,0)=0$ and hence $p(u,0)=0$.
This provides a blow-up equilibrium $u\!=\!z\!=\!0$, \textit{a priori}, along with automatic invariance of its $u$-axis $z\!=\!0$.
It also defines a meaningful local holonomy of $u$ over $z$, without previous analytic diagonalization.
See section \ref{Carsy0},
The interested reader may pursue the details of the complex blow-up analysis, in analogy to our section \ref{Con}.
We also recall the exact embedding of polynomial ODEs on $\C^2$ into scalar parabolic polynomial PDEs on $\S^1$, as outlined in section \ref{PDElift}.

\subsection{Real applications}\label{App}

\emph{``The differences between pure and applied mathematics are social, rather than scientific''}.
Thus spake Vladimir I. Arnold, opening a congress on Industrial and Applied Mathematics \cite{ICIAM}.
He continued with a political definition: \emph{``Mathematics is called \emph{applied}, if it has been supported by, both, \emph{CIA} and \emph{KGB}.''}
His paradigm, thereafter, was quite ``pure'' catastrophe theory.

To some, applications of ODEs and PDEs in complex variables may still seem a bit of a stretch.
Well, are they?

After a brief reminder of the long tradition on complex variables, we indicate several applied directions beyond complex-time blow-up loops around, and real-time blow-up stars towards, \emph{real-time} blow-up.
An obvious application of complex-valued differential equations are real-valued systems, which can be interpreted as coupled systems of real equations for their real and imaginary parts.
Passage from real to imaginary time then relates superficially ``unrelated'' equations, via commuting (semi)flows \eqref{flow}.
And at least Schrödinger equations have been around for a century, later augmented by complex Ginzburg-Landau equations and other variants.
We conclude with the role of complex time extensions and their relevance for spurious solutions in real-time discretizations of homoclinic orbits -- a topic related to adiabatic averaging, exponentially small ``spurious'' separatrix splitting, and ``invisible'' chaos.

Complex numbers, once upon a time, had been anathema. 
They were considered ``imaginary'' -- a terminology much protested by C.F. Gauss, though in vain.
Thanks to mathematical progress, quite real positive capacitors and inductances have found wide interpretation as ``imaginary'' resistors in classical and contemporary electrical engineering.

\emph{Real} parts $\mathbf{u}=\Re \mathbf{w}$ of holomorphic functions $\mathbf{w}=\mathbf{w}(\bz)$ of $\bz=\bxi+\mi \boldsymbol{\eta}$ are harmonic, i.e.\ they satisfy the planar Laplace equation $\Delta_{\bxi,\boldsymbol{\eta}}\, \mathbf{u}=0$.
In particular, harmonicity is therefore invariant under biholomorphic equivalences in $\bz$.
The Riemann mapping theorem and its Teichmüller variants then open the door to general geometries.
The resulting velocity fields $\nabla \mathbf{u}$, for example, may describe particle dynamics in divergence-free, irrotational fluid flows.
See also \eqref{Ricuv} below.

Semilinear elliptic real PDEs like $\Delta \mathbf{u}+f(\mathbf{u})=0$ in two space dimensions $(\bxi,\boldsymbol{\eta})\in\mathbb{R}^2$, although not in complex form, lead to mathematical interpretations with a two-dimensional ``time'' group $(\mathbb{R}^2,+)$ given by spatial shifts in $(\bxi,\boldsymbol{\eta})$.
See for example \cite{Vishik} and the references there.
Similarly to complex time $t=r+\mi s\in(\mathbb{C},+)$, spatial shifts in the two unbounded or periodic real ``time'' directions $\bxi$ and $\boldsymbol{\eta}$ commute; compare \eqref{flow}.
They are not Cauchy-Riemann related, however.
Alternatively, real analyticity of solutions defines local extensions to shifts $(\bxi,\boldsymbol{\eta})\in\mathbb{C}^2$, perhaps including spatial ``blow-up''.

Pseudoholomorphic, alias $J$-holomorphic, curves can be viewed as  first-order systems of such elliptic PDEs.
They involve the complex Cauchy-Riemann PDE operator $\partial_{\,\bar t}$\,, instead of the Laplacian, and describe the gradient-``flow'' associated to the action integral of real Hamiltonian systems, in compact symplectic settings.
Their investigation by Conley and Zehnder famously led to the solution of the Arnold conjecture on the minimal number of periodic orbits in real Hamiltonian mechanics, under time-periodic forcing.
See for example \cite{Hofer, Donaldson, Salamon}.
See \cite{Vishik} and the references there, for a related Conley index approach to solutions of semilinear elliptic PDEs in strip-shaped domains.

The periodic trigonometric functions, and their hyperbolic cousins, are solutions to the linear pendulum equations $\ddot x\pm x=0$, which are related by a sign change in the force law of physics, alias a purely mathematical passage from real to imaginary time.
The groundbreaking work of K. Weierstrass on analogous doubly periodic elliptic functions, of course, traces back to elliptic integrals as in the \emph{nonlinear pendulum} setting of section \ref{Hampen}.
Less common lore is the relation between Duffing pendulum equations of single-well and double-well potentials \cite{Duffing}:
\begin{equation}
\label{Duf}
\ddot x \pm x(1-x^2)=0\,.
\end{equation}
Superficially opposites, we may invoke the elementary relation $t\mapsto \mi t$ between real and imaginary time directions $t\in\C$, and recognize them as real cousins in a complex world.

Related examples are the transition from the ``soft spring'' mathematical pendulum \eqref{pen} with force law $g(x)=-\sin x$ to the ``hard spring'' pendulum $g(x)=-\sinh v$.
Alternatively, we may view the two ``pendula'' as equilibria of real reaction-diffusion equations.
For suitable classes of force laws $g$, the passage from real to imaginary space variables $x$ then corresponds to a transition from Dirichlet to Neumann boundary conditions -- a quite ominous transition within conventional physics.
See \cite{FiedlerYamaguti} for a more detailed discussion.

In reverse real time, we may also split our very first Riccati example \eqref{Ric0} into real and imaginary parts $x=u+\mi v$ as
\begin{equation}
\label{Ricuv}
\begin{aligned}
\dot u &= v^2-u^2\,, \\
\dot v &= -2\,uv\,.
\end{aligned}
\end{equation}
We can then reinterpret that ODE in terms of \emph{population dynamics} or abstract \emph{chemical reaction networks}, at least for $u,v\geq0$.
Note the global attractor $u\!=\!v\!=\!0$ there.
The discussion of the full complex Riccati equation in section \ref{Ric} is classical, of course.
Once upon a time, figure \ref{figRiccati} was motivated by field lines and potentials of electric and magnetic dipoles \cite{Marsden}.
Similarly, it illustrates flow lines and potential of a steady, incompressible, irrotational planar fluid flow with source and sink at $w=\pm1$.

Decomposing the original Fujita PDE \eqref{PDEw} into real and imaginary parts $\mathbf{w}=\mathbf{u}+\mi \mathbf{v}$, in the same spirit, leads to the \emph{real reaction-diffusion system}
\begin{equation}
\label{PDEuv}
\begin{aligned}
   \mathbf{u}_t=&\mathbf{u}_{\bxi\bxi}+\mathbf{u}^2-\mathbf{v}^2\,,   \\
   \mathbf{v}_t=&\mathbf{v}_{\bxi\bxi}+2\mathbf{uv}\,,  
\end{aligned}
\end{equation}
with equal diffusion coefficients for $\mathbf{u}$ and $\mathbf{v}$.
See \cite{Yanagida}, where \eqref{PDEuv} has been studied in real time as a viscous version of the nonlocal \emph{Constantin-Lax-Majda equation} \cite{Lax}, a spatially one-dimensional model for the vorticity equation of fluid dynamics.

The linear \emph{Schrödinger equation} makes use of complex state spaces $L^2$ for the quantum mechanical wave function $\boldsymbol{\psi}$.
Standard nonlinear versions of the Schrödinger equation require gauge invariance, like
\begin{equation}
\label{NLS}
\mathrm{i}\boldsymbol{\psi}_s=\boldsymbol{\psi}_{\bxi\bxi}+|\boldsymbol{\psi}|^2\boldsymbol{\psi}\,.
\end{equation}
Therefore they are at least cubic and involve hybrid ``probability'' coefficients like $|\boldsymbol{\psi}|^2$.
See for example \cite{Sulem}.

We have mentioned the nonconservative Schrödinger equation \eqref{psi6} as a cousin to the heat equation, in imaginary time.
It arises from the nonlinear Schrödinger equation as a toy model involving external electrical fields and other variants \cite{Jaquetteqp}.
For a discussion of several related PDEs involving $\boldsymbol{\psi}$ and its complex conjugate, see \cite{JaquetteHet}.
The complex Ginzburg-Landau equation features similar hybrid terms and is yet another rich source of real physics \cite{CGL} in real time.
The global comments of section \ref{Global} on hybrid scalar ODEs \eqref{hybrid}, however, do not apply here because we achieve neither fully holomorphic nor antiholomorphic structure by any factorization with an Euler multiplier.
At least well-posedness of some hybrid quadratic Schrödinger PDEs, including the antiholomorphic variant $\overline{w}^2$ of the purely quadratic Fujita nonlinearity $w^2$, has recently been surveyed by \cite{LiuR}.

Intriguingly, the phase space of the \emph{3d Navier-Stokes equation} has been complexified in \cite{LiSinai}.
By a purely mathematical tour de force, blow-up in finite real time has been claimed, for an open set of complex-valued initial conditions.
In our ODE setting, such blow-up for open sets of initial conditions occurs in the Poincaré domain of stable blow-up nodes; see section \ref{PoiMas}.
The Siegel domain of section \ref{SieLin}, in contrast, was characterized by saddles at infinity, and blow-up along their separatrices, only.

In ODEs \eqref{ODEF}, real one-step \emph{discretizations} of nonconstant real homoclinic orbits $\mathbf{x}\,{:}\ \mathbf{e}\leadsto \mathbf{e}$\,, alias saddle-loops, mandate yet another excursion to complex time.
Yes, first-order explicit Euler discretization already qualifies.
Once advertised as ``spurious solutions'', they defy real treatment by Melnikov functions.
Indeed their separatrix splitting is exponentially small, beyond finite order, in the discretization step $\eps$.
In the guise of rapid forcing with period $\eps$, great Poincaré had famously fallen into that trap in his flawed -- but accepted -- original submission for the Swedish King Oscar Prize.
First exponential smallness  estimates go back to the adiabatic elimination method of exponential averaging by Neishtadt \cite{Neishtadt, ArnoldEnc}.
For very fine asymptotic expansions, which are based on complex poles $t=T\not\in\R$ of $\mathbf{x}(t)$ and resolved the Poincaré difficulty after a century, in passing, see work by Lazutkin and Gelfreich \cite{Gelfreich01, Gelfreich02} and their references.
This is yet another body of work which has hardly been celebrated enough.
In \cite{FiedlerScheurle} we have contributed upper exponential estimates of order $\exp(-c/\eps)$ for this real discretization effect of ``invisible'' chaos, under the assumption that the homoclinic orbit  $\mathbf{x}(t)$ remains analytic in a strip $|\Im\,t|\leq c$ of complex time $t$.
In Gevrey spaces of high spatial regularity, these results have been extended to suitable classes of PDEs \cite{MatthiesDiss,Matthieshom,Matthiesell,MatthiesScheel}.
For numerically oriented PDE applications in terms of \emph{backward error analysis}, see for example \cite{Matthiesback,OliverWulff, WulffOliver} and their references.

In \cite{FiedlerClaudia}, we have raised a 1,000\;\euro\ \emph{prize question} concerning the existence of complex entire homoclinic loops $\mathbf{x}$ for entire vector fields $\mathbf{F}$, in several complex dimensions.
See sections 1.8 and 7 there, for further discussion.
A positive answer would establish elusive super-exponentially small separatrix splittings under rapidly time-periodic forcing or time-discretization.
As so many issues of differential equations in the complex domain, the question remains open.

\newpage



\begin{thebibliography}{9999)999}

{\footnotesize{

\bibitem[AbMa85]{AM}
R.~Abraham and J.E.~Marsden.
\emph{Foundations of Mechanics}, 2nd ed.
Addison-Wesley, 1985. 

\bibitem[APMR22]{MuciEss}
A.~Alvarez-Parrilla and J.\ Muciño-Raymundo.
Dynamics of singular complex analytic vector fields with essential singularities. I and II.
\emph{Conform.\ Geom.\ Dyn.} \textbf{21} (2017), 126--224, and \emph{J.\ Sing.} \textbf{24} (2022), 1--78.

\bibitem[Ama90]{Amann} 
H.~Amann. 
\emph{Ordinary Differential Equations. An Introduction to Nonlinear Analysis.}
De Gruyter, Berlin 1990.

\bibitem[Arn88]{ArnoldODE}
V.I.~Arnold.                       
\emph{Geometrical Methods in the Theory of Ordinary Differential Equations.}
Springer-Verlag, Berlin 1988.

\bibitem[Arn95]{ICIAM}
V.I.~Arnold. Opening talk at ICIAM 95.
K.~Kirchgässner, O.~Mahrenholtz, R.~Mennicken, Reinhard (eds.)
Third International Congress on Industrial and Applied Mathematics,
Hamburg 1995.
Akademie Verlag Berlin, 1996. Personal notes.

\bibitem[AKN06]{ArnoldEnc}
V.I.~Arnold, V.V.~Kozlov, A.I.~Neishtadt.
\emph{Mathematical Aspects of Classical and Celestial Mechanics.} 3rd ed.
Springer-Verlag, Berlin 2006;
\url{https://doi.org/10.1007/978-3-540-48926-9}

\bibitem[ALR18]{Llibre1}
J.C.~Artés, J.~Llibre, A.C.~Rezende.
\emph{Structurally Unstable Quadratic Vector Fields of Codimension One.}
Birkhäuser, Cham 2018. 

\bibitem[ALSV21]{Llibre2}
J.C.~Artés, J.~Llibre, D.~Schlomiuk, N.~Vulpe.
\emph{Geometric Configurations of Singularities of Planar Polynomial Differential Systems. 
A global classification in the quadratic case.} 
Birkhäuser, Cham 2021. 

\bibitem[Ber23]{Bern}
P.~Bernard.
The Siegel-Bruno linearization theorem. 
\emph{Regul.\ Chaot.\ Dyn.} \textbf{28} (2023), 756--762; 
\url{https://doi.org/10.1134/S1560354723040147}

\bibitem[BraDi10]{DiPol}
B.~Branner and K.~Dias.
Classification of complex polynomial vector fields in one complex variable.
\emph{J.\ Difference Eqs.\ Appl.} \textbf{16} (2010), 463--517;
\url{https://doi.org/10.1080/10236190903251746} 

\bibitem[BrKn86]{Bries}
E.~Brieskorn and H.~Knörrer.
\emph{Plane Algebraic Curves.}
Birkhäuser, Basel 1986. 

\bibitem[BrFie88]{brfi88}
P.~Brunovsk\'y and B.~Fiedler.
 Connecting orbits in scalar reaction diffusion equations.
 \emph{Dynamics Reported} \textbf{1} (1988), 57--89.

\bibitem[BrFie89]{brfi89}
P.~Brunovsk\'y and B.~Fiedler.
Connecting orbits in scalar reaction diffusion equations {II}: The complete solution.
 \emph{J.\ Diff.\ Eqs.} \textbf{81} (1989), 106--135.
 
\bibitem[Bryu71]{Bruno}
A.D.~Bryuno. 
Analytic form of differential equations. 
\emph{Trans.\ Moscow Math.\ Soc.} \textbf{25} (1971), 131--288.

\bibitem[COS16]{COS}
C.-H.~Cho, H.~Okamoto, M.~Sh\={o}ji. 
A blow-up problem for a nonlinear heat equation in the complex plane of time. 
\emph{Japan J.\ Ind.\ Appl.\ Math.} \textbf{33} (2016), 145--166.

\bibitem[CLM85]{Lax}
P.~Constantin, P.D.~Lax, A.J.~Majda. 
A simple one-dimensional model for the three-dimensional vorticity equation.
\emph{Comm.\ Pure Appl.\ Math.} \textbf{38} (1985), 715--724.

\bibitem[DiGa21]{DiRat}
K.~Dias and A.~Garijo.
On the separatrix graph of a rational vector field on the Riemann sphere.
\emph{J.\ Diff.\ Eqs.} \textbf{282} (2021), 541--565;
\url{https://doi.org/10.1016/j.jde.2021.02.021}

\bibitem[Don05]{Donaldson}
S.K.~Donaldson.
What is … a pseudoholomorphic curve?
\emph{Notices Am. Math. Soc.} \textbf{52} (2005), 1026--1027. 

\bibitem[Duf1918]{Duffing}
G.~Duffing.
Erzwungene Schwingungen bei veränderlicher Eigenfrequenz und ihre technische Bedeutung. 
Vieweg, Braunschweig 1918.
 
\bibitem[DNZ23]{CGL}
G.K.~Duong, N.~Nouaili, H.~Zaag.
\emph{Construction of Blowup Solutions for the Complex Ginzburg-Landau Equation with Critical Parameters.} 
Mem.\ Am.\ Math.\ Soc.\ \textbf{1411}, Providence RI 2023. 

\bibitem[FKW24]{Fasondini24}
M.~Fasondini, J.R.~King, J.A.C.~Weideman.
Complex-plane singularity dynamics for blow up in a nonlinear heat equation: analysis and computation. 
\emph{Nonlinearity} \textbf{37} (2024);
\url{https://doi.org/10.1088/1361-6544/ad700b}

\bibitem[Fie02]{fi02}
B.~Fiedler (ed.).  \emph{Handbook of Dynamical
Systems} \textbf{2}. Elsevier, Amsterdam 2002.

\bibitem[Fie23]{FiedlerClaudia}
B.~Fiedler.
Real chaos and complex time. (2023); 
\url{https://arxiv.org/abs/2310.08136}

\bibitem[Fie25a]{FiedlerShilnikov}
B.~Fiedler.
Scalar polynomial vector fields in real and complex time. 
\emph{Reg.\ Chaotic Dyn.} \textbf{30} (2025), 188--225;
\url{https://doi.org/10.1134/S1560354725020030}

\bibitem[Fie25b]{FiedlerYamaguti}
B.~Fiedler.
Real-time blow-up and connection graphs of rational vector fields on the Riemann sphere. 
\emph{Jap. J.~Industr. Appl. Math.} (2025);
\url{https://doi.org/10.1007/s13160-025-00749-8}

\bibitem[FieMan00]{FiedlerMantel}
B.~Fiedler and R.-M.~Mantel.
Crossover collisions of scroll wave filaments. 
\emph{Doc.\ Math.\ J.\ DMV} {\textbf 5} (2000), 695--731;
\url{https://doi.org/10.4171/dm/92}

\bibitem[FieMat07]{FiedlerMatano}
B.~Fiedler and H.~Matano.
Blow-up shapes on fast unstable manifolds of one-dimensional reaction-diffusion equations. 
\emph{J.~Dyn.\ Differ.\ Eqs.} \textbf{19} (2007), 867--893.

\bibitem[FiePol90]{FPolacik} 
B.~Fiedler and P.~Pol\'{a}\v{c}ik.
Complicated dynamics of one-dimensional reaction diffusion equations with a non-local term.
\emph{Proc.\ R.\ Soc.\ Edinb.} \textbf{115A} (1990), 167--192.

\bibitem[FieRo23]{firoSFB}
B.~Fiedler and C.~Rocha.
Design of Sturm global attractors 1: Meanders with three noses, and reversibility.
\emph{Chaos} \textbf{33}, 083127 (2023); 
\url{https://doi.org/10.1063/5.0147634}

\bibitem[FieRo24]{firoFusco}
B.~Fiedler and C.~Rocha.
Design of Sturm global attractors 2: Time-reversible Chafee-Infante lattices of 3-nose meanders.
\emph{São Paulo J.\ Math.\ Sciences.} (2024);
\url{https://doi.org/10.1007/s40863-023-00385-5}

\bibitem[FieSch02]{FiedlerScheel}
B.~Fiedler and A.~Scheel. 
Dynamics of Reaction-Diffusion Patterns.
In \emph{Trends in Nonlinear Analysis. Festschrift dedicated to Willi Jäger
for his 60th birthday.}
M.~Kirkilionis, R.~Rannacher, F.~Tomi (eds.), 23--152.
Springer-Verlag, Heidelberg 2002.

\bibitem[FSV98]{Vishik}
B.~Fiedler, A.~Scheel, M.I.~Vishik. 
Large patterns of elliptic systems in infinite cylinders. 
\emph{J.\ Math.\ Pures Appl.} \textbf{77} (1998), 879--907.

\bibitem[FieSch96]{FiedlerScheurle}
B.~Fiedler and J.~Scheurle.
\emph{Discretization of Homoclinic Orbits, Rapid Forcing and “Invisible” Chaos.}
Mem.\ Am.\ Math.\ Soc.\ \textbf{570}, Providence RI 1996. 

\bibitem[FieStu25]{FiedlerFila}
B.~Fiedler and H.~Stuke.
Real eternal PDE solutions are not complex entire: a quadratic parabolic example.
\emph{J.\ Ell.\ Par.\ Eqs.} (2025), 53pp.
\url{doi.org/10.1007/s41808-024-00309-0}

\bibitem[For81]{Forster}
O.~Forster.
\emph{Lectures on Riemann Surfaces.}
Springer-Verlag, New York 1981.

\bibitem[FuRo24]{RochaDichotomy}
G.~Fusco and C.~Rocha.
On the structure of the infinitesimal generators of scalar one-dimensional semigroups with discrete Lyapunov functionals. 
\emph{São Paulo J.\ Math.\ Sci.} \textbf{18} (2024), 1026--1054. 

\bibitem[Gala04]{Galaktionov}
V.A.~Galaktionov.
\emph{Geometric Sturmian Theory of Nonlinear Parabolic Equations and Applications.} 
Chapman\&Hall/CRC, Boca Raton FL 2004. 

\bibitem[GL01]{Gelfreich01}
V.G.~Gelfreich and V.F.~Lazutkin. 
Splitting of separatrices: perturbation theory and exponential smallness. 
\emph{Russ.\ Math.\ Surv.} \textbf{56} (2001), 499--558.

\bibitem[Gel02]{Gelfreich02}
V.G.~Gelfreich.
Numerics and exponential smallness. 
In \cite{fi02} (2002), 265--312. 

\bibitem[GNSY13]{Yanagida}
J.-S.~Guo, H.~Ninomiya, M.~Shimojo, E.~Yanagida.
Convergence and blow-up of solutions for a complex-valued heat equation with a quadratic nonlinearity. 
\emph{Trans.\ Am.\ Math.\ Soc.} \textbf{365} (2013), 2447--2467.

\bibitem[HLW02]{Backward}
E.\ Hairer, Ch.\ Lubich, G.\ Wanner.
\emph{Geometric Numerical Integration. Structure-preserving Algorithms for Ordinary Differential Equations. }
Springer-Verlag, Berlin 2002. 

\bibitem[Hale85]{HaleEmbed}
J.K.~Hale.
Flows on centre manifolds for scalar functional differential equations. 
\emph{Proc. R. Soc. Edinb.} \textbf{101A} (1985), 193--201.  

\bibitem[Har02]{Hartman}
Ph.~Hartman.
\emph{Ordinary Differential Equations.}
SIAM, Providence RI 2002.

\bibitem[HoZe94]{Hofer}
H.~Hofer and E.~Zehnder.
\emph{Symplectic Invariants and Hamiltonian Dynamics.}
Birkhäuser, Basel 1994. 

\bibitem[IlYa08]{Ilya}
Y.~Ilyashenko and S.~Yakovenko.
\emph{Lectures on Analytic Differential Equations.}
AMS, Providence RI 2008.

\bibitem[Jaq24]{Jaquetteqp}
J.~Jaquette.
Quasiperiodicity and blowup in integrable subsystems of nonconservative nonlinear Schrödinger equations. 
\emph{J.\ Dyn.\ Differ.\ Eqs.} \textbf{36} (2024), 1--25; 
\url{https://doi.org/10.1007/s10884-021-10112-3}

\bibitem[JLT22a]{JaquetteHet}
J.~Jaquette, J.-P.~Lessard, A.~Takayasu.
Global dynamics in nonconservative nonlinear Schrödinger equations. 
\emph{Adv.\ Math.} \textbf{398} (2022), 108234. 

\bibitem[JLT22b]{JaquetteStuke}
J.~Jaquette, J.-P.~Lessard, A.~Takayasu.
Singularities and heteroclinic connections in complex-valued evolutionary equations with a quadratic nonlinearity. 
\emph{Comm.\ Nonlinear Sci.\ Numer.\ Simul.} \textbf{107} (2022), 106188.

\bibitem[Jost06]{Jost}
J.~Jost.
\emph{Compact Riemann Surfaces. An Introduction to Contemporary Mathematics.}
Springer-Verlag, Berlin 2006.

\bibitem[KaiLe25]{LebEss}
N.~Kainz and D.~Lebiedz.
Separatrix configurations in holomorphic flows.
(2025);
\url{https://doi.org/10.48550/arXiv.2505.14594}

\bibitem[KliRou21]{Rousseaud=4}
M.~Klimeš and Ch.~Rousseau.
Remarks on rational vector fields on $\mathbb{CP}^1$.
\emph{J.\ Dyn.\ Control Syst.} \textbf{27} (2021), 293--320. 

\bibitem[Kup13]{Kupitz}
D.~Kupitz.
\emph{Wechselwirkungen von Scrollwellen.} 
Dissertation, Otto-von-Guericke-Universität Magdeburg 2013;
\url{http://hydra.nat.uni-magdeburg.de/fnw/pdf/DKde.pdf}

\bibitem[KupH12]{KupitzHauser}
D.~Kupitz and M.J.B.~Hauser.
Interaction of a pair of parallel scroll waves.
\emph{J.\ Phys.\ Chem.} A \textbf{117} (2013), 12711--12718;
\url{https://doi.org/10.1021/jp409269u}

\bibitem[Lam09]{Lamotke}
K.~Lamotke.
\emph{Riemannsche Flächen.} 
Springer-Verlag, Heidelberg 2009.

\bibitem[Lang95]{Lang}
S.~Lang. 
\emph{Differential and Riemannian Manifolds.}
Springer-Verlag, New York 1995.

\bibitem[Lap23]{Lappicy}
Ph.~Lappicy and E.~Beatriz. 
An energy formula for fully nonlinear degenerate parabolic equations in one spatial dimension.
\emph{Math.\ Ann.} (2023); 
\url{https://doi.org/10.1007/s00208-023-02740-5}

\bibitem[LiSin08]{LiSinai}
D.~Li and Y.G.~Sinai.
Blow ups of complex solutions of the 3D Navier-Stokes system and renormalization group method.
\emph{J.\ Eur.\ Math.\ Soc.} \textbf{10} (2008), 267--313;
\url{https://doi.org/10.4171/JEMS/111}

\bibitem[Liu25]{LiuR}
R.~Liu.
Local well-posedness of the periodic nonlinear Schrödinger equation with a quadratic nonlinearity $\bar{u}^2$ in negative Sobolev spaces.
\emph{J.\ Dyn.\ Differential Eqs.} \textbf{37} (2025), 509--538;
\url{https://doi.org/10.1007/s10884-023-10295-x}

\bibitem[MPS96]{MPDichotomy}
J.~Mallet-Paret and G.~Sell. 
Systems of differential delay equations: Floquet multipliers and discrete Lyapunov functions. 
\emph{J.\ Diff.\ Eqs.} \textbf{125}
(1996), 385--440.

\bibitem[MH99]{Marsden}
J.E.~Marsden and M.J.~Hoffman.
\emph{Basic Complex Analysis.}
Freeman, New York 1999.

\bibitem[Mas82]{Masuda1}
K.~Masuda.
Blow-up of solutions of some nonlinear diffusion equations.  
\emph{North-Holland Math.\ Stud.} \textbf{81} (1982), 119--131. 

\bibitem[Mas84]{Masuda2}
K.~Masuda.
Analytic solutions of some nonlinear diffusion equations. 
\emph{Math.\ Z.} \textbf{187} (1984), 61--73. 

\bibitem[Mat82]{MatanoLap}
H.~Matano.
Nonincrease of the lap-number of a solution for a one-dimensional semilinear parabolic equation. 
\emph{J.\ Fac.\ Sci.\ Univ.\ Tokyo I A} \textbf{29} (1982), 401--441. 

\bibitem[Mat88]{Matano}
H.~Matano.
Asymptotic behavior of solutions of semilinear heat equations on
$\S^1$.
In \emph{Nonlinear Diffusion Equations and their Equilibrium States
  {II}}. W.-M.~Ni, L.A.~Peletier, J.~Serrin (eds.), Springer-Verlag,
New York 1988.

\bibitem[Mat01]{MatthiesDiss}
K.~Matthies.
Time-averaging under fast periodic forcing of parabolic partial differential equations: Exponential estimates.
\emph{J.\ Differ.\ Eqs.} \textbf{174} (2001), 133--180. 

\bibitem[Mat03a]{Matthieshom}
K.~Matthies.
Exponentially small splitting of homoclinic orbits of parabolic differential equations under periodic forcing.
\emph{Discr.\ Contin.\ Dyn.\ Syst.} \textbf{9}  (2003), 585--602. 

\bibitem[Mat03b]{Matthiesback}
K.~Matthies. 
Backward error analysis of a full discretisation scheme for a class of parabolic
partial differential equations. \emph{Nonlin.\ Analysis TMA} \textbf{52} (2003), 805–826.

\bibitem[Mat05]{Matthiesell}
K.~Matthies.
Homogenisation of exponential order for elliptic systems in infinite cylinders.
\emph{Asympt.\ Analysis} \textbf{43} (2005), 205--232. 

\bibitem[MS03]{MatthiesScheel}
K.~Matthies and A.~Scheel.
Exponential averaging for Hamiltonian evolution equations.
\emph{Trans.\ Am.\ Math.\ Soc.} \textbf{355} (2003), 747--773. 

\bibitem[McDSa12]{Salamon}
D.~McDuff and D.~Salamon.
\emph{$J$-holomorphic Curves and Symplectic Topology.} 2nd ed.
Am. Math. Soc., Providence, RI 2012. 

\bibitem[Mil06]{Milnor}
J.~Milnor.
\emph{Dynamics in One Complex Variable.} 3rd ed.
Princeton, NJ 2006. 

\bibitem[Mos73]{Moser}
J.~Moser.
\emph{Stable and Random Motions in Dynamical Systems. With Special Emphasis on Celestial Mechanics.}
Princeton, N.J.\ 1973.

\bibitem[Muc02]{MuciRiem}
J.~Muciño-Raymundo.
Complex structures adapted to smooth vector fields. 
\emph{Math.\ Ann.} \textbf{322} (2002), 229--265. 

\bibitem[Nei84]{Neishtadt}
A.I.~Neishtadt. 
On the separation of motions in systems with rapidly rotating phase. 
\emph{J.\ Appl.\ Math.\ Mech.} \textbf{48} (1984), 134--139.

\bibitem[OWW04]{OliverWulff}
M.~Oliver, M.~West, C.~Wulff.
Approximate momentum conservation for spatial semidiscretizations of semilinear wave equations. 
\emph{Numer.\ Math.} \textbf{97} (2004), 493--535. 

\bibitem[PM01]{P-M}
R.~Pérez-Marco. 
Total convergence or general divergence in small divisors.
\emph{Comm.\ Math.\ Phys.} \textbf{223} (2001), 451--464.

\bibitem[Pol02]{PolacikEmbed}
P.~Poláčik.
Parabolic equations: Asymptotic behavior and dynamics on invariant manifolds.
In \cite{fi02}, 2002, 835--883.

\bibitem[QS19]{Quittner}
P.~Quittner and Ph.~Souplet.
\emph{Superlinear Parabolic Problems. Blow-Up, Global Existence and Steady States.} 2nd ed., Birkhäuser, Cham 2019. 

\bibitem[Ra95]{Raugelthin}
G.~Raugel.
Dynamics of partial differential equations on thin domains.
In \emph{Dynamical Systems.}
R.~Johnson (ed.). CIME, Montecatini Terme, Italy 1994. Springer-Verlag, Berlin 1995, 208--315.

\bibitem[Ra02]{Raugelattr}
G.~Raugel.
Global attractors in partial differential equations.
In \cite{fi02}, 2002, 885--982.

\bibitem[Rel40]{Rellich}
F.~Rellich.
Elliptische Funktionen und die ganzen Lösungen von $y''=f(y)$.
\emph{Math.\ Z.} \textbf{47} (1940), 153--160;
\url{https://doi.org/10.1007/bf01180954} 

\bibitem[SaFie92]{FSandstede}
 B.~Sandstede and B.~Fiedler.
 Dynamics of periodically forced parabolic equations on the circle. 
 \emph{J.~Erg.\ Th.\ Dyn.\ Syst.} \textbf{12} (1992), 559--571.


\bibitem[Sie42]{Siegel}
C.L.~Siegel. 
Iteration of analytic functions. 
\emph{Ann.\ Math.\ II Ser.} \textbf{43} (1942), 607--612.

\bibitem[Slo24a]{oeispol}
L.P.A.~Sloane (ed.).
\emph{The Online Encyclopedia of Integer Sequences.}
A002995 (2024); \url{https://oeis.org/A002995}

\bibitem[Sto00]{countpol2}
A.~Stoimenow.
On the number of chord diagrams.
\emph{Disc.\ Math.} \textbf{218} (2000), 209--233.

\bibitem[Stu17]{Stukediss}
H.~Stuke.
\emph{Blow-up in Complex Time.}
Dissertation, Freie Universität Berlin 2017;
\url{http://dx.doi.org/10.17169/refubium-11743}

\bibitem[Stu18]{Stukearxiv}
H.~Stuke.
\emph{Complex time blow-up of the nonlinear heat equation.}
(2018); 
\url{https://arxiv.org/abs/1812.10707} 

\bibitem[SuSu99]{Sulem}
C.~and P.-L.~Sulem.
\emph{The Nonlinear Schrödinger Equation. Self-Focusing and Wave Collapse.}
Springer-Verlag, New York 1999. 

\bibitem[TLJO22]{JaquetteMasuda}
A.~Takayasu, J.-P.~Lessard, J.~Jaquette, H.~Okamoto.
Rigorous numerics for nonlinear heat equations in the complex plane of time.
\emph{Numer.\ Math.} \textbf{151} (2022), 693--752. 

\bibitem[Ush80]{Ushiki1}
S.~Ushiki. 
On unstable manifolds of analytic diffeomorphisms of the plane.
\emph{RIMS Kyoto Kokyuroku} \textbf{403} (1980), 1--7.

\bibitem[Ush81]{Ushiki2}
S.~Ushiki. 
Unstable manifolds of analytic dynamical systems.
\emph{J.\ Math.\ Kyoto Univ.} \textbf{21} (1981), 763--785.

\bibitem[Van89]{Vdb}
A.~Vanderbauwhede.
Centre manifolds, normal forms and elementary bifurcations.
\emph{Dynamics Reported} \textbf{2} (1989), 89--169. 

\bibitem[Wal72]{countpol1} 
D.W.~Walkup.
The number of plane trees.
\emph{Mathematika} \textbf{19} (1972), 200--204.

\bibitem[Wit41]{Wittich1}
H.~Wittich.
Ganze Lösungen der Differentialgleichung $w''=f(w)$.
\emph{Math.\ Z.} \textbf{47} (1941), 422--426;
\url{https://doi.org/10.1007/BF01180973}

\bibitem[Wit50]{Wittich2}
H.~Wittich.
Ganze transzendente Lösungen algebraischer Differentialgleichungen.
\emph{Math.\ Ann.} \textbf{122} (1950), 37--46;
\url{https://doi.org/10.1007/BF01342967}

\bibitem[WO16]{WulffOliver}
C.~Wulff and M.~Oliver.
Exponentially accurate Hamiltonian embeddings of symplectic A-stable Runge-Kutta methods for Hamiltonian semilinear evolution equations. 
\emph{Proc.\ R.\ Soc.\ Edinb.\ A, Math.} \textbf{146} (2016), 1265--1301.

\bibitem[Yoc95]{Yoccoz}
J.-C.~Yoccoz.
Théorème de Siegel, nombres de Bruno et polynômes quadratiques.
\emph{Astérisque} \textbf{231} (1995), 3--88.

\bibitem[Zel68]{Zelenyak}
T.I.~Zelenyak.
Stabilization of solutions of boundary value problems for a second order parabolic equation with one space variable.
\emph{Diff.\ Eqs.} \textbf{4} (1968), 17--22.
}}

\end{thebibliography}
\end{document}